\documentclass{amsart}

\usepackage{amssymb}
\usepackage{url} 
\usepackage{graphicx}
\usepackage[all]{xy}
\usepackage[all]{xy}

\newtheorem{theorem}{Theorem}[section]
\newtheorem{corollary}[theorem]{Corollary}
\newtheorem{proposition}[theorem]{Proposition}
\newtheorem{lemma}[theorem]{Lemma}
\newtheorem{claim}[theorem]{Claim}
\newtheorem{notation}[theorem]{Notation}
\newtheorem{caution}[theorem]{Caution}

\theoremstyle{remark}
\newtheorem{remark}[theorem]{Remark}

\theoremstyle{definition}
\newtheorem{definition}[theorem]{Definition}
\newtheorem{example}[theorem]{Example}

\numberwithin{equation}{section}
\numberwithin{theorem}{section}

\begin{document}

\title{Cohomology theory in 2-categories}

\author{Hiroyuki NAKAOKA}
\address{Graduate School of Mathematical Sciences, The University of Tokyo 
3-8-1 Komaba, Meguro, Tokyo, 153-8914 Japan}
\email{deutsche@ms.u-tokyo.ac.jp}

\thanks{The author wishes to thank Professor Toshiyuki Katsura for his encouragement}
\thanks{The author is supported by JSPS}
\thanks{This paper was submitted to {\it Theory and Applications of Categories} in 25 November 2007}

\begin{abstract}
Recently, symmetric categorical groups are used for the study of the Brauer groups of symmetric monoidal categories. As a part of these efforts, some algebraic structures of the 2-category of symmetric categorical groups $\mathrm{SCG}$ are being investigated. In this paper, we consider a 2-categorical analogue of an abelian category, in such a way that it contains $\mathrm{SCG}$ as an example. As a main theorem, we construct a long cohomology 2-exact sequence from any extension of complexes in such a 2-category. Our axiomatic and self-dual definition will enable us to simplify the proofs, by analogy with abelian categories.

\end{abstract}

\maketitle

\section{Introduction}\label{Introduction}

In 1970s, B. Pareigis started his study on the Brauer groups of symmetric monoidal categories in \cite{P4}. Around 2000, the notion of symmetric categorical groups are introduced to this study by E. M. Vitale in \cite{VPicBr} (see also \cite{VBrBT}). By definition, a symmetric categorical group is a categorification of an abelian group, and in this sense the 2-category of symmetric categorical groups $\mathrm{SCG}$ can be regarded as a 2-dimensional analogue of the category $\mathrm{Ab}$ of abelian groups. As such, $\mathrm{SCG}$ and its variants (e.g. 2-category $\mathbb{G}\text{-}\mathrm{SMod}$ of symmetric categorical groups with $\mathbb{G}$-action where $\mathbb{G}$ is a fixed categorical group) admit a 2-dimensional analogue of the homological algebra in $\mathrm{Ab}$.

For example,  E. M. Vitale constructed for any monoidal functor $F:\mathbf{C}\rightarrow\mathbf{D}$ between symmetric monoidal categories $\mathbf{C}$ and $\mathbf{D}$, a 2-exact sequence of Picard and Brauer categorical groups
\[
\mathcal{P}(\mathbf{C})\rightarrow\mathcal{P}(\mathbf{D})\rightarrow\overline{\mathcal{F}}\rightarrow\mathcal{B}(\mathbf{C})\rightarrow\mathcal{B}(\mathbf{C}).
\]
By taking $\pi_0$ and $\pi_1$, we can induce the well-known Picard-Brauer and Unit-Picard exact sequences of abelian groups respectively. In \cite{RMV}, A. del R\'{\i}o, J. Mart\'{\i}nez-Moreno and E. M. Vitale defined a more subtle notion of the relative 2-exactness, and succeeded in constructing a cohomology long 2-exact sequence from any short relatively 2-exact sequence of complexes in $\mathrm{SCG}$. In this paper, we consider a 2-categorical analogue of an abelian category, in such a way that it contains $\mathrm{SCG}$ as an example, so as to treat $\mathrm{SCG}$ and their variants in a more abstract, unified way.

In section \ref{Preliminaries}, we review general definitions in a 2-category and properties of $\mathrm{SCG}$, with simple comments. In section \ref{Definition of a RelEx2Cat}, we define the notion of a relatively exact 2-category as a generalization of $\mathrm{SCG}$, also as a 2-dimensional analogue of an abelian category. We try to make the homological algebra in $\mathrm{SCG}$ (\cite{RMV}) work well in this general 2-category. It will be worthy to note that our definition of a relatively exact 2-category is self-dual.

\begin{tabular}
[c]{|c|c|c|}\hline
& category & 2-category\\\hline
general theory & abelian category & relatively exact 2-category\\\hline
example & $\mathrm{Ab}$ & $\mathrm{SCG}$\\\hline
\end{tabular}

In section \ref{Existence of PFS}, we show the existence of {\it proper factorization systems} in any relatively exact 2-category, which will make several diagram lemmas more easy to handle. In any abelian category, any morphism $f$ can be written in the form $f=e\circ m$ (uniquely up to an isomorphism), where $e$ is epimorphic and $m$ is monomorphic. As a 2-dimensional analogue, we show that any 1-cell $f$ in a relatively exact 2-category $\mathbf{S}$ admits the following two ways of factorization:

{\rm (1)} $i\circ m\Longrightarrow f$ where $i$ is fully cofaithful and $m$ is faithful.

{\rm (2)} $e\circ j\Longrightarrow f$ where $e$ is cofaithful and $j$ is fully faithful.\\
(In the case of $\mathrm{SCG}$, see \cite{KV}.)
In section \ref{Definition of Rel2Ex}, complexes in $\mathbf{S}$ and the relative 2-exactness are defined, generalizing those in $\mathrm{SCG}$ (\cite{RMV}). Since we start from the self-dual definition, we can make good use of duality in the proofs. In section \ref{Long Cohomology Sequence}, as a main theorem, we construct a long cohomology 2-exact sequence from any short relatively 2-exact sequence (i.e. an extension) of complexes. Our proof is purely diagrammatic, and is an analogy of that for an abelian category. In section \ref{Definition of Rel2Ex} and \ref{Long Cohomology Sequence}, several 2-dimensional diagram lemmas are shown. Most of them have 1-dimensional analogues in an abelian category, so we only have to be careful about the compatibility of 2-cells.

Since $\mathrm{SCG}$ is an example of a relatively exact 2-category, we expect some other 2-categories constructed from $\mathrm{SCG}$ will be a relatively exact 2-category. For example, $\mathbb{G}\text{-}\mathrm{SMod}$, $\mathrm{SCG}\times\mathrm{SCG}$ and the 2-category of bifunctors from $\mathrm{SCG}$ are candidates. We will examine such examples in forthcoming papers.

\section{Preliminaries}\label{Preliminaries}

\paragraph{Definitions in a 2-category}\quad\\
%
%
\begin{notation}
Throughout this paper, $\mathbf{S}$ denotes a 2-category $($in the strict sense$)$. We use the following notation.\\
$\mathbf{S}^0$, $\mathbf{S}^1$, $\mathbf{S}^2:$ class of 0-cells, 1-cells, and 2-cells in $\mathbf{S}$, respectively.\\
$\mathbf{S}^1(A,B):$ 1-cells from $A$ to $B$, where $A,B\in\mathbf{S}^0$.\\
$\mathbf{S}^2(f,g):$ 2-cells from $f$ to $g$, where $f,g\in\mathbf{S}^1(A,B)$ for certain $A,B\in\mathbf{S}^0$.\\
$\mathbf{S}(A,B):$ $\mathrm{Hom}$-category between $A$ and $B$ $($i.e. $\mathrm{Ob}(\mathbf{S}(A,B))=\mathbf{S}^1(A,B)$, $\mathbf{S}(A,B)(f,g)=\mathbf{S}^2(f,g))$.

In diagrams, $\longrightarrow$ represents a 1-cell, $\Longrightarrow$ represents a 2-cell, $\circ$ represents a horizontal composition, and $\cdot$ represents a vertical composition. We use capital letters $A,B,\ldots$ for 0-cells, small letters $f,g,\ldots$ for 1-cells, and Greek symbols $\alpha,\beta,\ldots$ for 2-cells.
\end{notation}

For example, one of the conditions in the definition of a 2-category can be
written as follows (see for example \cite{Mac}):

\begin{remark}
For any diagram in $\mathbf{S}$
\[
\xy(-16,0)*+{A}="0";
(0,0)*+{B}="2";
(16,0)*+{C}="4";
{\ar@/^1.00pc/^{f_1} "0";"2"};
{\ar@/_1.00pc/_{f_2} "0";"2"};
{\ar@/^1.00pc/^{g_1} "2";"4"};
{\ar@/_1.00pc/_{g_2} "2";"4"};
{\ar@{=>}_{\alpha} (-8,2);(-8,-2)};
{\ar@{=>}_{\beta} (8,2);(8,-2)};
\endxy
,
\]
we have
\begin{equation}
(f_1\circ\beta)\cdot(\alpha\circ g_2)=(\alpha\circ g_1)\cdot(f_2\circ\beta).
\label{2-cell}
\end{equation}
(Note: composition is always written diagrammatically.)
\end{remark}
This equality is frequently used in later arguments.

Products, pullbacks, difference kernels and their duals are defined by the universality.

\begin{definition}
For any $A_1$ and $A_2\in\mathbf{S}^0$, their product $(A_1\times A_2,p_1,p_2)$ is defined as follows$:$\\
{\rm (a)} $A_1\times A_2\in\mathbf{S}^0$, $p_i\in\mathbf{S}^1(A_1\times A_2,A_i)$ $(i=1,2)$.\\
{\rm (b1)} $($existence of a factorization$)$

For any $X\in\mathbf{S}^0$ and $q_i\in\mathbf{S}^1(X,A_i)$ $(i=1,2)$, there exist $q\in\mathbf{S}^1(X,A_1\times A_2)$ and $\xi_i\in\mathbf{S}^2(q\circ p_i,q_i)\ (i=1,2)$.
\[
\xy
(0,6)*+{X}="0";
(-18,-10)*+{A_1}="2";
(0,-10)*+{A_1\times A_2}="4";
(18,-10)*+{A_2}="6";
{\ar_{q_1} "0";"2"};
{\ar_{q} "0";"4"};
{\ar^{q_2} "0";"6"};
{\ar^{p_1} "4";"2"};
{\ar_{p_2} "4";"6"};
{\ar@{=>}^{\xi_1} (-3,-7);(-7.5,-2.5)};
{\ar@{=>}_{\xi_2} (3,-7);(7.5,-2.5)};
\endxy
\]
\noindent{\rm (b2)} $($uniqueness of the factorization$)$

For any factorizations $(q,\xi_1,\xi_2)$ and $(q^{\prime},\xi_1^{\prime},\xi_2^{\prime})$ which satisfy {\rm (b1)}, there exists a unique 2-cell $\eta\in\mathbf{S}^2(q,q^{\prime})$ such that $(\eta\circ p_i)\cdot\xi_i^{\prime}=\xi_i\ (i=1,2)$.
\[
\xy
(-10,7)*+{q\circ p_i}="0";
(10,7)*+{q^{\prime}\circ p_i}="2";
(0,-8)*+{q_i}="4";
(0,8)*+{}="6";
{\ar@{=>}^{\eta\circ p_i} "0";"2"};
{\ar@{=>}_{\xi_i} "0";"4"};
{\ar@{=>}^{\xi_i^{\prime}} "2";"4"};
{\ar@{}|\circlearrowright"6";"4"};
\endxy
\]
\end{definition}
The coproduct of $A_1$and $A_2$ is defined dually.

\begin{definition}\label{DefPullback}
For any $A_1,A_2,B\in\mathbf{S}^0$ and $f_i\in\mathbf{S}^1(A_i,B)\ (i=1,2)$, the pullback $(A_1\times_BA_2,f_1^{\prime},f_2^{\prime},\xi)$ of $f_1$ and $f_2$ is defined as follows$:$

\noindent{\rm (a)} $A_1\times_BA_2\in\mathbf{S}^0$, $f_1^{\prime}\in\mathbf{S}^1(A_1\times_BA_2,A_2)$, $f_2^{\prime}\in\mathbf{S}^1(A_1\times_BA_2,A_1)$, $\xi\in\mathbf{S}^2(f_1^{\prime}\circ f_2,f_2^{\prime}\circ f_1)$.
\[
\xy
(-12,0)*+{A_1\times_B A_2}="0";
(0,9)*+{A_2}="2";
(0,-9)*+{A_1}="4";
(12,0)*+{B}="6";
{\ar^{f^{\prime}_1} "0";"2"};
{\ar_{f^{\prime}_2} "0";"4"};
{\ar_{f_1} "4";"6"};
{\ar^{f_2} "2";"6"};
{\ar@{=>}^{\xi} (0,3);(0,-3)};
\endxy
\]
\noindent{\rm (b1)} $($existence of a factorization$)$

For any $X\in\mathbf{S}^0$, $g_1\in\mathbf{S}^1(X,A_2)$, $g_2\in\mathbf{S}^1(X,A_1)$ and $\eta\in\mathbf{S}^2(g_1\circ f_2,g_2\circ f_1)$, there exist $g\in\mathbf{S}^1(X,A_1\times A_2),\xi_i\in\mathbf{S}^2(g\circ f_i^{\prime},g_i)\ (i=1,2)$ such that $(\xi_1\circ f_2)\cdot\eta=(g\circ\xi)\cdot(\xi_2\circ f_1)$.
\[
\xy
(-10,0)*+{A_1\times_B A_2}="0";
(0,10)*+{A_2}="2";
(0,-10)*+{A_1}="4";
(10,)*+{B}="6";
(-32,0)*+{X}="8";
{\ar^{f^{\prime}_1} "0";"2"};
{\ar_{f^{\prime}_2} "0";"4"};
{\ar_{f_1} "4";"6"};
{\ar^{f_2} "2";"6"};
{\ar^{g} "8";"0"};
{\ar@/^1.20pc/^{g_1} "8";"2"};
{\ar@/_1.20pc/_{g_2} "8";"4"};
{\ar@{=>}^{\xi} (0,5);(0,-5)};
{\ar@{=>}_{\xi_1} (-12,2);(-16,8)};
{\ar@{=>}^{\xi_2} (-12,-2);(-16,-8)};
\endxy
\quad
\xy
(-14,6)*+{g\circ f^{\prime}_1\circ f_2}="0";
(14,6)*+{g_1\circ f_2}="2";
(-14,-6)*+{g\circ f^{\prime}_2\circ f_1}="4";
(14,-6)*+{g_2\circ f_1}="6";
{\ar@{=>}^{\xi_1\circ f_2} "0";"2"};
{\ar@{=>}_{g\circ\xi} "0";"4"};
{\ar@{=>}^{\eta} "2";"6"};
{\ar@{=>}_{\xi_2\circ f_1} "4";"6"};
{\ar@{}|\circlearrowright"0";"6"};
\endxy
\]
\noindent{\rm (b2)} $($uniqueness of the factorization$)$

For any factorizations $(g,\xi_1,\xi_2)$ and $(g^{\prime},\xi_1^{\prime},\xi_2^{\prime})$ which satisfy {\rm (b1)}, there exists a unique 2-cell $\zeta\in\mathbf{S}^2(g,g^{\prime})$ such that $(\zeta\circ f_i^{\prime})\cdot\xi_i^{\prime}=\xi_i\ (i=1,2)$.
\end{definition}
The pushout of $f_i\in\mathbf{S}^1(A,B_i)\ (i=1,2)$ is defined dually.

\begin{definition}
For any $A,B\in\mathbf{S}^0$ and $f,g\in\mathbf{S}^1(A,B)$, the difference kernel 
\[ (\mathrm{DK}(f,g),d_{(f,g)},\varphi_{(f,g)}) \]
of $f$ and $g$ is defined as follows$:$

\noindent{\rm (a)} $\mathrm{DK}(f,g)\in\mathbf{S}^0$, $d_{(f,g)}\in\mathbf{S}^1(\mathrm{DK}(f,g),A)$, $\varphi_{(f,g)}\in\mathbf{S}^2(d_{(f,g)}\circ f,d_{(f,g)}\circ g)$.
\[
\xy
(-28,0)*+{\mathrm{DK}(f,g)}="0";
(0,0)*+{A}="2";
(20,0)*+{B}="4";
{\ar_{d_{(f,g)}} "0";"2"};
{\ar@/^0.50pc/^{f} "2";"4"};
{\ar@/_0.50pc/_{g} "2";"4"};
\endxy
,\ 
\xy
(-20,0)*+{\mathrm{DK}(f,g)}="0";
(16,0)*+{B}="4";
{\ar@/^1.10pc/^{d_{(f,g)}\circ f} "0";"4"};
{\ar@/_1.10pc/_{d_{(f,g)}\circ g} "0";"4"};
{\ar@{=>}_{\varphi_{(f,g)}} (0,2)*{};(0,-2)*{}} ;
\endxy
\]

\noindent{\rm (b1)} $($existence of a factorization$)$

For any $X\in\mathbf{S}^0$, $d\in\mathbf{S}^1(X,A)$, $\varphi\in\mathbf{S}^2(d\circ f,d\circ g)$, there exist $\underline{d}\in\mathbf{S}^1(X,\mathrm{DK}(f,g)),\underline{\varphi}\in\mathbf{S}^2(\underline{d}\circ d_{(f,g)},d)$ such that $(\underline{d}\circ\varphi_{(f,g)})\cdot(\underline{\varphi}\circ g)=(\underline{\varphi}\circ f)\cdot\varphi.$

\noindent{\rm (b2)} $($uniqueness of the factorization$)$

For any factorizations $(\underline{d},\underline{\varphi})$ and $(\underline{d^{\prime}},\underline{\varphi^{\prime}})$ which satisfy {\rm (b1)}, there exists a unique 2-cell $\eta\in\mathbf{S}^2(\underline{d},\underline{d^{\prime}})$ such that $(\eta\circ d_{(f,g)})\cdot\underline{\varphi^{\prime}}=\underline{\varphi}$.
\end{definition}
The difference cokernel of $f$ and $g$ is defined dually.

The following definition is from \cite{DV}.
\begin{definition}
Let $f\in\mathbf{S}^1(A,B)$.\\
{\rm (1)} $f$ is said to be faithful if $f^{\flat}:=-\circ f:\mathbf{S}^1(C,A)\rightarrow\mathbf{S}^1(C,B)$ is faithful for any $C\in\mathbf{S}^0$.\\
{\rm (2)} $f$ is said to be fully faithful if $f^{\flat}$ is fully faithful for any $C\in\mathbf{S}^0$.\\
{\rm (3)} $f$ is said to be cofaithful if $f^{\sharp}:=f\circ-:\mathbf{S}^1(B,C)\rightarrow\mathbf{S}^1(A,C)$ is faithful for any $C\in\mathbf{S}^0$.\\
{\rm (4)} $f$ is said to be fully cofaithful if $f^{\sharp}$ is fully faithful for any $C\in\mathbf{S}^0$.
\end{definition}

\paragraph{Properties of $\mathrm{SCG}$}\quad\\

By definition, a symmetric categorical group is a symmetric monoidal category $(\mathbb{G},\otimes,0)$, in which each arrow is an isomorphism and each object has an inverse up to an equivalence with respect to the tensor $\otimes$. More precisely;
\begin{definition}
A symmetric categorical group $(\mathbb{G},\otimes,0)$ consists of\\
{\rm (a1)} a category $\mathbb{G}$\\
{\rm (a2)} a tensor functor $\otimes:\mathbb{G}\times\mathbb{G}\rightarrow\mathbb{G}$\\
{\rm (a3)} a unit object $0\in \mathrm{Ob}(\mathbb{G)}$\\
{\rm (a4)} natural isomorphisms
\begin{eqnarray*}
&\alpha_{A,B,C}:A\otimes(B\otimes C)\rightarrow(A\otimes B)\otimes C,&\\
&\lambda_A:0\otimes A\rightarrow A,\ \rho_A:A\otimes0\rightarrow A,\ \gamma_{A,B}:A\otimes B\rightarrow B\otimes A&
\end{eqnarray*}
which satisfy certain compatibility conditions $($cf. {\rm \cite{P1}}$)$, and the following two conditions are satisfied:\\
{\rm (b1)} For any $A,B\in \mathrm{Ob}(\mathbb{G)}$ and $f\in\mathbb{G}(A,B)$, there exists $g\in\mathbb{G}(B,A)$ such that $f\circ g=\mathrm{id}_A$, $g\circ f=\mathrm{id}_B$.\\
{\rm (b2)} For any $A\in \mathrm{Ob}(\mathbb{G})$, there exist $A^{\ast}\in \mathrm{Ob}(\mathbb{G})$ and $\eta_A\in\mathbb{G}(0,A\otimes A^{\ast})$.
\end{definition}

In particular, there is a \lq zero categorical group' $0$, which consists of only one object $0$ and one morphism $\mathrm{id}_0$.

\begin{definition}
For symmetric categorical groups $\mathbb{G}$ and $\mathbb{H}$, a monoidal functor $F$ from $\mathbb{G}$ to $\mathbb{H}$ consists of\\
{\rm (a1)} a functor $F:\mathbb{G}\rightarrow\mathbb{H}$\\
{\rm (a2)} natural isomorphisms
\[ F_{A,B}:F(A\otimes B)\rightarrow F(A)\otimes F(B)\ \text{and}\ F_I:F(0)\rightarrow 0 \]
which satisfy certain compatibilities with $\alpha$, $\lambda$, $\rho$, $\gamma$. $($cf. {\rm \cite{P1}}$)$
\end{definition}

\begin{remark}\label{AddedRemark}
For any monoidal functors $F:\mathbb{G}\rightarrow\mathbb{H}$ and $G:\mathbb{H}\rightarrow\mathbb{K}$, their composition $F\circ G:\mathbb{G}\rightarrow\mathbb{K}$ is defined by
\begin{eqnarray}
(F\circ G)_{A,B}&:=&G(F_{A,B})\circ G_{F(A),F(B)}\\
(F\circ G)_I&:=&G(F_I)\circ G_I. \label{1-4}
\end{eqnarray}

\end{remark}
In particular, there is a \lq zero monoidal functor' $0_{\mathbb{G},\mathbb{H}}:\mathbb{G}\rightarrow\mathbb{H}$ for each $\mathbb{G}$ and $\mathbb{H}$, which sends every object in $\mathbb{G}$ to $0_{\mathbb{H}}$, every arrow in $\mathbb{G}$ to $\mathrm{id}_{0_{\mathbb{H}}}$, and $(0_{\mathbb{G},\mathbb{H}})_{A,B}=\lambda_0^{-1}=\rho_0^{-1}$, $(0_{\mathbb{G},\mathbb{H}})_I=\mathrm{id}_0$. It is easy to see that $0_{\mathbb{G},\mathbb{H}}\circ0_{\mathbb{H},\mathbb{K}}=0_{\mathbb{G},\mathbb{K}}$ ($\forall\mathbb{G},\mathbb{H},\mathbb{K}$).

\begin{remark}
Our notion of a monoidal functor is equal to that of a \lq$\gamma$-monoidal functor' in \cite{RMV}.
\end{remark}

\begin{definition}\label{TRANSF}
For monoidal functors $F,G:$ $\mathbb{G}\rightarrow\mathbb{H}$, a natural transformation $\varphi$ from $F$ to $G$ is said to be a monoidal transformation if it satisfies
\begin{eqnarray}
\varphi_{A\otimes B}\circ G_{A,B}&=&F_{A,B}\circ(\varphi_A\otimes\varphi_B)\nonumber\\
F_I&=&\varphi_0\circ G_I.
\label{1-2}
\end{eqnarray}
\end{definition}

The following remark is from \cite{VPicBr}.
\begin{remark}\label{DualityRem}
By condition {\rm (b2)}, it is shown that there exists a 2-cell $\varepsilon_A\in\mathbb{G}(A^{\ast}\otimes A,0)$ for each object $A$, such that the following compositions are identities:
\begin{eqnarray*}
& A\underset{\lambda_A^{-1}}{\longrightarrow}0\otimes A\underset{\eta
_A\otimes 1}{\longrightarrow}(A\otimes A^{\ast})\otimes A\underset{\alpha^{-1}}{\longrightarrow}A\otimes(A^{\ast}\otimes A)\underset{1\otimes\varepsilon_A}{\longrightarrow}A\otimes0\underset{\rho_A}{\longrightarrow}A\\
& A^{\ast}\underset{\rho_{A^{\ast}}^{-1}}{\longrightarrow}A^{\ast}\otimes 0\underset{1\otimes\eta_A}{\longrightarrow}A^{\ast}\otimes(A\otimes A^{\ast})\underset{\alpha}{\longrightarrow}(A^{\ast}\otimes A)\otimes A^{\ast}\underset{\varepsilon_A\otimes1}{\longrightarrow}0\otimes A^{\ast}\underset{\lambda_{A^{\ast}}}{\longrightarrow}A^{\ast}
\end{eqnarray*}

For each monoidal functor $F:\mathbb{G}\rightarrow\mathbb{H}$, there exists a natural morphism $\iota_{F,A}:F(A^{\ast})\rightarrow F(A)^{\ast}$.
\end{remark}

\begin{definition}
$\mathrm{SCG}$ is defined to be the 2-category whose 0-cells are symmetric categorical groups, 1-cells are monoidal functors, and 2-cells are monoidal transformations.
\end{definition}

The following two propositions are satisfied in $\mathrm{SCG}$ (see for example
\cite{BV}).

\begin{proposition}
For any symmetric categorical groups $\mathbb{G}$ and $\mathbb{H}$, if we define a monoidal functor $F\otimes_{\mathbb{G},\mathbb{H}}G:\mathbb{G}\rightarrow\mathbb{H}$ by
\begin{eqnarray*}
F\otimes_{\mathbb{G},\mathbb{H}}G(A)&:=&F(A)\otimes_{\mathbb{H}}G(A)\\
(F\otimes_{\mathbb{G},\mathbb{H}}G)_{A,B}&:=&(F(A\otimes B)\otimes G(A\otimes B)\\
&\overset{F_{A,B}\otimes G_{A,B}}{\longrightarrow}&F(A)\otimes F(B)\otimes G(A)\otimes G(B)\\
&\overset{\simeq}{\longrightarrow}&F(A)\otimes G(A)\otimes F(B)\otimes G(B))\\
(F\otimes_{\mathbb{G},\mathbb{H}}G)_I&:=&(F(I)\otimes G(I)\overset{F_I\otimes G_I}{\longrightarrow}I\otimes I\overset{\simeq}{\longrightarrow}I),
\end{eqnarray*}
then $(\mathrm{SCG}(\mathbb{G},\mathbb{H)},\otimes_{\mathbb{G},\mathbb{H}},0_{\mathbb{G},\mathbb{H}})$ becomes again a symmetric categorical group with appropriately defined $\alpha,\lambda,\rho,\gamma$, and
\[ \mathrm{Hom}=\mathrm{SCG}(-,-):\mathrm{SCG}\times\mathrm{SCG}\rightarrow\mathrm{SCG} \]
becomes a 2-functor $($cf. section 6 in {\rm \cite{BV}}$)$.
\end{proposition}

In $\mathrm{SCG}$, by definition of the zero categorical group we have $\mathbf{S}^1(\mathbb{G},0)=\{0_{\mathbb{G},0}\}$, while $\mathbf{S}^1(0,\mathbb{G})$ may have more than one objects. In this point $\mathrm{SCG}$ might be said to have \lq non self-dual' structure, but $\mathbf{S}^1(\mathbb{G},0)$ and $\mathbf{S}^1(0,\mathbb{G})$ have the following \lq self-dual' property.
\begin{remark}
\label{RRRRem}
{\rm (1)} For any symmetric categorical group $\mathbb{G}$ and any monoidal functor $F:\mathbb{G}\rightarrow 0$, there exists a unique 2-cell $\varphi:F\Longrightarrow0_{\mathbb{G},0}$.

{\rm (2)} For any symmetric categorical group $\mathbb{G}$ and any monoidal functor $F:0\rightarrow\mathbb{G}$, there exists a unique 2-cell $\varphi:F\Longrightarrow0_{0,\mathbb{G}}$.
\end{remark}

\begin{proof}
{\rm (1)} follows from the fact that the zero categorical group has only one morphism $\mathrm{id}_0$. {\rm (2)} follows from condition $(\ref{1-2})$ in Definition \ref{TRANSF}.
\end{proof}

The usual compatibility arguments show the following Lemma.

\begin{lemma}
\label{CriticalLemma}
Let $F:\mathbb{G}\rightarrow\mathbb{H}$ be a monoidal functor. For any $A,B\in \mathrm{Ob}(\mathbb{G}),$
\[
\xy
(-34,4)*+{\Phi_{A,B} :}="12";
(-20,4)*+{\mathbb{G}(A,B)}="0";
(20,4)*+{\mathbb{G}(A\otimes B^{\ast},0)}="2";
(-20,0)*+{\rotatebox{90}{$\in$}}="4";
(20,0)*+{\rotatebox{90}{$\in$}}="6";
(-20,-4)*+{f}="8";
(20,-4)*+{(f\otimes 1_{B^{\ast}})\circ\eta^{-1}_B}="10";
{\ar_{} "0";"2"};
{\ar@{|->}_{} "8";"10"};
\endxy
\]
and
\[
\xy
(-42,4)*+{\Psi_{A,B} :}="12";
(-24,4)*+{\mathbb{G}(A \otimes B^{\ast},0)}="0";
(28,4)*+{\mathbb{G}(A,B)}="2";
(-24,0)*+{\rotatebox{90}{$\in$}}="4";
(28,0)*+{\rotatebox{90}{$\in$}}="6";
(-24,-4)*+{g}="8";
(28,-4)*+{\rho_A^{-1}\circ(1_A\otimes\varepsilon^{-1}_B)\circ\alpha_{_{A,B^{\ast},B}}\circ(g\otimes 1_B)\circ\lambda_B}="10";
{\ar_{} "0";"2"};
{\ar@{|->}_{} "8";"10"};
\endxy
\]
are mutually inverse, and the following diagram is commutative$;$
\[
\xy
(-16,14)*+{\mathbb{G}(A,B)}="0";
(16,14)*+{\mathbb{G}(A\otimes B^{\ast},0)}="2";
(-28,0)*+{\mathbb{H}(F(A),F(B))}="4";
(0,-16)*+{\mathbb{H}(F(A)\otimes F(B)^{\ast},0)}="6";
(28,0)*+{\mathbb{H}(F(A\otimes B^{\ast}),F(0)),}="8";
(0,16)*+{}="10";
{\ar^{\Phi_{A,B}} "0";"2"};
{\ar_{F} "0";"4"};
{\ar^{F} "2";"8"};
{\ar@{}|\circlearrowright"6";"10"};
{\ar_{\Phi_{F(A),F(B)}} "4";"6"};
{\ar^{\Theta_{A,B }^F} "8";"6"};
\endxy
\]
where $\Theta_{A,B}^F$ is defined by
\[
\xy
(-48,4)*+{\Theta^F_{A,B}:}="12";
(-24,4)*+{\mathbb{H}(F(A\otimes B^{\ast}),F(0))}="0";
(28,4)*+{\mathbb{H}(F(A)\otimes F(B)^{\ast},0)}="2";
(-24,0)*+{\rotatebox{90}{$\in$}}="4";
(28,0)*+{\rotatebox{90}{$\in$}}="6";
(-24,-4)*+{h}="8";
(28,-4)*+{(1_{F(A)}\otimes(\iota^F_B)^{-1})\circ(F_{A,B^{\ast}})^{-1}\circ h\circ F_I.}="10";
{\ar_{} "0";"2"};
{\ar@{|->}_{} "8";"10"};
\endxy
\]
\end{lemma}

\section{Definition of a relatively exact 2-category}\label{Definition of a RelEx2Cat}

\paragraph{Locally $\mathrm{SCG}$ 2-category}\quad\\

We define a locally $\mathrm{SCG}$ 2-category not only as a 2-category whose $\mathrm{Hom}$-categories are $\mathrm{SCG}$, but with some more conditions, in order to let it be a 2-dimensional analogue of that of an additive category.

\begin{definition}
\label{LocSCG}
A locally small 2-category $\mathbf{S}$ is said to be locally $\mathrm{SCG}$ if the following conditions are satisfied$:$\\
{\rm (A1)} For every $A,B\in\mathbf{S}^0$, there is a given functor $\otimes_{A,B}:\mathbf{S}(A,B)\times\mathbf{S}(A,B)\rightarrow\mathbf{S}(A,B)$, and a given object $0_{A,B}\in \mathrm{Ob}(\mathbf{S}(A,B))=\mathbf{S}^1(A,B)$ such that $(\mathbf{S}(A,B),\otimes_{A,B},0_{A,B})$ becomes a symmetric categorical group, and the following naturality conditions are satisfied$:$
\[ 0_{A,B}\circ0_{B,C}=0_{A,C}\quad\quad(\forall A,B,C\in\mathbf{S}^0) \]
\noindent{\rm (A2)} $\mathrm{Hom}=\mathbf{S}(-,-):\mathbf{S}\times\mathbf{S}\rightarrow\mathrm{SCG}$ is a 2-functor which satisfies for any $A,B,C\in\mathbf{S}^0$,
\begin{eqnarray}
(0_{A,B})_I^{\sharp}&=\mathrm{id}_{0_{A,C}}\in\mathbf{S}^2(0_{A,C},0_{A,C})\label{1-6}\\
(0_{A,B})_I^{\flat}&=\mathrm{id}_{0_{C,B}}\in\mathbf{S}^2(0_{C,B},0_{C,B}).
\label{1-7}
\end{eqnarray}
\noindent{\rm (A3)} There is a 0-cell $0\in\mathbf{S}^0$ called a zero object, which satisfy the following conditions:\\
{\rm (a3-1)} $\mathbf{S}(0,0)$ is the zero categorical group.\\
{\rm (a3-2)} For any $A\in\mathbf{S}^0$ and $f\in\mathbf{S}^1(0,A)$, there exists a unique 2-cell $\theta_f\in\mathbf{S}^2(f,0_{0,A})$.\\
{\rm (a3-3)} For any $A\in\mathbf{S}^0$ and $f\in\mathbf{S}^1(A,0)$, there exists a unique 2-cell $\tau_f\in\mathbf{S}^2(f,0_{A,0})$.\\
{\rm (A4)} For any $A,B\in\mathbf{S}^0$, their product and coproduct exist.
\end{definition}

Let us explain about these conditions.
\begin{remark}
By condition {\rm (A1)} of Definition \ref{LocSCG}, every 2-cell in a locally $\mathrm{SCG}$ 2-category becomes invertible, as in the case of $\mathrm{SCG}$ (cf. \cite{VPicBr}). This helps us to avoid being fussy about the directions of 2-cells in many propositions and lemmas, and we use the word \lq dual' simply to reverse 1-cells.
\end{remark}

\begin{remark}
By condition {\rm (A2)} in Definition \ref{LocSCG},
\begin{eqnarray*}
f^{\sharp}&:=&f\circ-:\mathbf{S}(B,C)\rightarrow\mathbf{S}(A,C)\\
f^{\flat}&:=&-\circ f:\mathbf{S}(C,A)\rightarrow\mathbf{S}(C,B)
\end{eqnarray*}
are monoidal functors $(\forall C\in\mathbf{S}^0)$ for any $f\in\mathbf{S}^1(A,B)$, and the following naturality conditions are satisfied:\\
{\rm (a2-1)} For any $f\in\mathbf{S}^1(A,B),g\in\mathbf{S}^1(B,C)$ and $D\in\mathbf{S}^0$, we have $(f\circ g)^{\sharp}=g^{\sharp}\circ f^{\sharp}$ as monoidal functors.
\[
\xy
(-24,0)*+{A}="0";
(-8,0)*+{B}="2";
(8,0)*+{C}="4";
(24,0)*+{D}="6";
{\ar^{f} "0";"2"};
{\ar^{g} "2";"4"};
{\ar^{} "4";"6"};
\endxy
\quad
\xy
(-12,7)*+{\mathbf{S}(C,D)}="0";
(12,7)*+{\mathbf{S}(B,D)}="2";
(0,-8)*+{\mathbf{S}(A,D)}="4";
(0,8)*+{}="6";
{\ar^{g^{\sharp}} "0";"2"};
{\ar_{(f\circ g)^{\sharp}} "0";"4"};
{\ar^{f^{\sharp} } "2";"4"};
{\ar@{}|\circlearrowright"6";"4"};
\endxy
\]
\noindent{\rm (a2-2)} The dual of {\rm (a2-1)} for $-^{\flat}$.
\\
{\rm (a2-3)} For any $f\in\mathbf{S}^1(A,B),g\in\mathbf{S}^1(C,D)$, we have $f^{\sharp}\circ g^{\flat}=g^{\flat}\circ f^{\sharp}$ as monoidal functors.
\[
\xy
(-24,0)*+{A}="0";
(-8,0)*+{B}="2";
(8,0)*+{C}="4";
(24,0)*+{D}="6";
{\ar^{f} "0";"2"};
{\ar^{} "2";"4"};
{\ar^{g} "4";"6"};
\endxy
\quad
\xy
(-12,6)*+{\mathbf{S}(B,C)}="0";
(12,6)*+{\mathbf{S}(A,C)}="2";
(-12,-6)*+{\mathbf{S}(B,D)}="4";
(12,-6)*+{\mathbf{S}(A,D)}="6";
{\ar^{f^{\sharp}} "0";"2"};
{\ar_{g^{\flat}} "0";"4"};
{\ar^{g^{\flat}} "2";"6"};
{\ar_{f^{\sharp}} "4";"6"};
{\ar@{}|\circlearrowright"0";"6"};
\endxy
\]

Since already $(f\circ g)^{\sharp}=g^{\sharp}\circ f^{\sharp}$ as functors, {\rm (a2-1)} means $\ (f\circ g)_I^{\sharp}=(g^{\sharp}\circ f^{\sharp})_I$, and by $(\ref{1-4})$ in Remark \ref{AddedRemark}, this is equivalent to
\[ (f\circ g)_I^{\sharp}=f^{\sharp}(g_I^{\sharp})\cdot f_I^{\sharp}=(f\circ g_I^{\sharp})\cdot f_I^{\sharp}. \]

Similarly, we obtain
\begin{eqnarray}
(f\circ g)_I^{\flat}&=(f_I^{\flat}\circ g)\cdot g_I^{\flat} \label{1-9},\\
(f_I^{\sharp}\circ g)\cdot g_I^{\flat}&=(f\circ g_I^{\flat})\cdot f_I^{\sharp}. \label{1-10}
\end{eqnarray}
\end{remark}

\begin{remark}
\label{Contra}
By condition {\rm (A2)}, for any $f,g\in\mathbf{S}^1(A,B)$ and any $\alpha\in\mathbf{S}^2(f,g)$, $\alpha\circ-:f^{\sharp}\Rightarrow g^{\sharp}$ becomes a monoidal transformation. So, the diagrams
\[
\xy
(-24,6)*+{f\circ(k\otimes h)}="0";
(-24,-6)*+{(f\circ k)\otimes(f\circ h)}="2";
(24,6)*+{g\circ(k\otimes h)}="4";
(24,-6)*+{(g\circ k)\otimes(g\circ h)}="6";
{\ar@{=>}^{\alpha\circ(k\otimes h)} "0";"4"};
{\ar@{=>}_{f^{\sharp}_{k,h}} "0";"2"};
{\ar@{=>}_{(\alpha\circ k)\otimes(\alpha\circ h)} "2";"6"};
{\ar@{=>}^{g^{\sharp}_{k,h}} "4";"6"};
{\ar@{}|\circlearrowright"0";"6"};
\endxy
\ \text{and}\ 
\xy
(-12,6)*+{f\circ 0_{B,C}}="0";
(12,6)*+{g\circ 0_{B,C}}="2";
(0,-8)*+{0_{A,C}}="4";
(0,9)*+{}="6";
{\ar@{=>}^{\alpha\circ0_{B,C}} "0";"2"};
{\ar@{=>}_{f^{\sharp}_I} "0";"4"};
{\ar@{=>}^{g^{\sharp}_I} "2";"4"};
{\ar@{}|\circlearrowright"4";"6"};
\endxy
\]
are commutative for any $C\in\mathbf{S}^0$ and $k,h\in\mathbf{S}^1(B,C)$. Similar statement also holds for \ $-\circ\alpha:f^{\flat}\Rightarrow g^{\flat}$.
\end{remark}

\begin{corollary}
\label{ContraCor}
In a locally $\mathrm{SCG}$ 2-category $\mathbf{S}$, the following are satisfied$:$

{\rm (1)} For any diagram in $\mathbf{S}$
\[
\xy
(-16,0)*+{C}="0";
(0,0)*+{A}="2";
(16,0)*+{B}="4";
{\ar@/^0.75pc/^{h} "0";"2"};
{\ar@/_0.75pc/_{0_{C,A}} "0";"2"};
{\ar@/^0.75pc/^{f} "2";"4"};
{\ar@/_0.75pc/_{g} "2";"4"};
{\ar@{=>}^{\varepsilon} (-8,2)*{};(-8,-2)*{}} ;
{\ar@{=>}^{\alpha} (8,2)*{};(8,-2)*{}} ;
\endxy
\,
\]
we have
\begin{equation}
h\circ\alpha=(\varepsilon\circ f)\cdot f_I^{\flat}\cdot g_I^{\flat-1}\cdot(\varepsilon^{-1}\circ g).
\label{1-13}
\end{equation}

{\rm (2)} For any diagram in $\mathbf{S}$
\[
\xy
(-16,0)*+{A}="0";
(0,0)*+{B}="2";
(16,0)*+{C}="4";
{\ar@/^0.75pc/^{f} "0";"2"};
{\ar@/_0.75pc/_{g} "0";"2"};
{\ar@/^0.75pc/^{h} "2";"4"};
{\ar@/_0.75pc/_{0_{B,C}} "2";"4"};
{\ar@{=>}^{\alpha} (-8,2)*{};(-8,-2)*{}};
{\ar@{=>}^{\varepsilon} (8,2)*{};(8,-2)*{}};
\endxy
,
\]
we have
\begin{equation}
\alpha\circ h=(f\circ\varepsilon)\cdot f_I^{\sharp}\cdot g_I^{\sharp-1}\cdot(g\circ\varepsilon^{-1}).
\label{1-14}
\end{equation}

{\rm (3)} For any diagram in $\mathbf{S}$
\[
\xy
(-16,0)*+{A}="0";
(0,0)*+{B}="2";
(16,0)*+{C}="4";
{\ar@/^0.75pc/^{f} "0";"2"};
{\ar@/_0.75pc/_{0_{A,B}} "0";"2"};
{\ar@/^0.75pc/^{g} "2";"4"};
{\ar@/_0.75pc/_{0_{B,C}} "2";"4"};
{\ar@{=>}^{\alpha} (-8,2)*{};(-8,-2)*{}};
{\ar@{=>}^{\beta} (8,2)*{};(8,-2)*{}};
\endxy
,
\]
we have
\begin{equation}
(f\circ\beta)\cdot f_I^{\sharp}=(\alpha\circ g)\cdot g_I^{\flat}.
\label{eqtn}
\end{equation}
\end{corollary}
\begin{proof}
{\rm (1)} $(h\circ\alpha)\underset{\ref{2-cell}}{=}(\varepsilon\circ f)\cdot(0_{C,A}\circ\alpha)\cdot(\varepsilon^{-1}\circ g)=(\varepsilon\circ f)\cdot f_I^{\flat}\cdot g_I^{\flat-1}\cdot(\varepsilon^{-1}\circ g)$. {\rm (2)} is the dual of {\rm (1)}. And {\rm (3)} follows from $(\ref{1-6})$, $(\ref{1-7})$, $(\ref{1-13})$, $(\ref{1-14})$.
\end{proof}

\begin{remark}
We don't require a locally $\mathrm{SCG}$ 2-category to satisfy $\mathbf{S}^1(A,0)=\{0_{A,0}\}$, for the sake of duality (see the comments before Remark \ref{RRRRem} ).
\end{remark}

\paragraph{Relatively exact 2-category}\quad\\

\begin{definition}\label{RelExCat}
Let $\mathbf{S}$ be a locally $\mathrm{SCG}$ 2-category. $\mathbf{S}$ is said to be relatively exact if the following conditions are satisfied$:$

{\rm (B1)} For any 1-cell $f\in\mathbf{S}^1(A,B)$, its kernel and cokernel exist.

{\rm (B2)} For any 1-cell $f\in\mathbf{S}^1(A,B)$, $f$ is faithful if and only if $f=\mathrm{ker}(\mathrm{cok}(f))$.

{\rm (B3)} For any 1-cell $f\in\mathbf{S}^1(A,B)$, $f$ is cofaithful if and only if $f=\mathrm{cok}(\mathrm{ker}(g))$.


\noindent It is shown in {\rm \cite{VPicBr}} that $\mathrm{SCG}$ satisfies these conditions.
\end{definition}


Let us explain about these conditions.

\begin{definition}\label{Kernel}
Let $\mathbf{S}$ be a locally $\mathrm{SCG}$ 2-category. For any $f\in\mathbf{S}^1(A,B)$, its kernel $(\mathrm{Ker}(f),\mathrm{ker}(f),\varepsilon_f)$ is defined by universality as follows $($we abbreviate $\mathrm{ker}(f)$ to $k(f)):$\\
{\rm (a)} $\mathrm{Ker}(f)\in\mathbf{S}^0$, $k(f)\in\mathbf{S}^1(\mathrm{Ker}(f),A)$, $\varepsilon_f\in\mathbf{S}^2(k(f)\circ f,0)$.\\
{\rm (b1)} $($existence of a factorization$)$

For any $K\in\mathbf{S}^0$, $k\in\mathbf{S}^1(K,A)$ and $\varepsilon\in\mathbf{S}^2(k\circ f,0)$, there exist $\underline{k}\in\mathbf{S}^1(K,\mathrm{Ker}(f))$ and $\underline{\varepsilon}\in\mathbf{S}^2(\underline{k}\circ k(f),k)$ such that $(\underline{\varepsilon}\circ f)\cdot\varepsilon=(\underline{k}\circ\varepsilon_f)\cdot(\underline{k})_I^{\sharp}$.
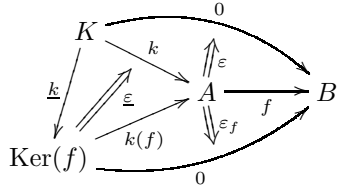
\[
\xy
(-29,-9)*+{\mathrm{Ker}(f)}="0";
(-24,8)*+{K}="2";
(-8,0)*+{A}="4";
(8,0)*+{B}="6";
{\ar_<<<<<<{k(f)} "0";"4"};
{\ar^{k} "2";"4"};
{\ar_{\underline{k}} "2";"0"};
{\ar_{f} "4";"6"};
{\ar@/^1.25pc/^{0} "2";"6"};
{\ar@/_1.25pc/_{0} "0";"6"};
{\ar@{=>}_<<<{\varepsilon} (-8,2)*{};(-7,7)*{}};
{\ar@{=>}^<<<<{\varepsilon_f} (-8,-2)*{};(-7,-7)*{}};
{\ar@{=>}_<<<<<<<<{\underline{\varepsilon}} (-25,-5)*{};(-18,3)*{}};
\endxy
\]
\noindent{\rm (b2)} $($uniqueness of the factorization$)$

For any factorizations $(\underline{k},\underline{\varepsilon})$ and $(\underline{k}^{\prime},\underline{\varepsilon}^{\prime})$ which satisfy {\rm (b1)}, there exists a unique 2-cell $\xi\in\mathbf{S}^2(\underline{k},\underline{k}^{\prime})$ such that $(\xi\circ k(f))\cdot\underline{\varepsilon}^{\prime}=\underline{\varepsilon}$.
\end{definition}

\begin{remark}\label{Cokernel}
{\rm (1)} By its universality, the kernel of $f$ is unique up to an equivalence. We write this equivalence class again $\mathrm{Ker}(f)=[\mathrm{Ker}(f),k(f),\varepsilon_f]$.

{\rm (2)} It is also easy to see that if $f$ and $f^{\prime}$ are equivalent, then
\[ [\mathrm{Ker}(f),k(f),\varepsilon_f]=[\mathrm{Ker}(f^{\prime}),k(f^{\prime}),\varepsilon_{f^{\prime}}]. \]

For any $f$, its cokernel $\mathrm{Cok}(f)=[\mathrm{Cok}(f),c(f),\pi_{f}]$ is defined dually, and the dual statements also hold for the cokernel.
\end{remark}

\begin{remark}
Let $\mathbf{S}$ be a locally $\mathrm{SCG}$ 2-category, and let $f$ be in $\mathbf{S}^1(A,B)$.

For any pair $(k,\varepsilon)$ with $k\in\mathbf{S}^1(0,A),\varepsilon\in\mathbf{S}^2(k\circ f,0)$
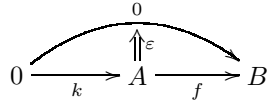
\[
\xy
(-8,0)*+{0}="4";
(8,0)*+{A}="6";
(24,0)*+{B}="8";
{\ar_{k} "4";"6"};
{\ar_{f} "6";"8"};
{\ar@/^1.65pc/^{0} "4";"8"};
{\ar@{=>}_<<<{\varepsilon} (8,2)*{};(8,6)*{}} ;
\endxy
\]
and for any pair $(k^{\prime},\varepsilon^{\prime})$ with $k^{\prime}\in\mathbf{S}^1(0,A),\varepsilon^{\prime}\in\mathbf{S}^2(k\circ f,0)$,
there exists a unique 2-cell $\xi\in\mathbf{S}^2(k,k^{\prime})$ such that $(\xi\circ f)\cdot\varepsilon^{\prime}=\varepsilon.$
\end{remark}
\begin{proof}
By condition {\rm (a3-2)} of Definition \ref{LocSCG}, $\varepsilon\in\mathbf{S}^2(k\circ f,0)$ must be equal to the unique 2-cell $(\theta_k\circ f)\cdot f_I^{\flat}$. Similarly we have $\varepsilon^{\prime}=(\theta_{k^{\prime}}\circ f)\cdot f_I^{\flat}$, and, $\xi$ should be the unique 2-cell $\theta_k\cdot\theta_{k^{\prime}}^{-1}\in\mathbf{S}^2(k,k^{\prime})$, which satisfies $(\xi\circ f)\cdot\varepsilon^{\prime}=\varepsilon$.
\end{proof}

From this, it makes no ambiguity if we abbreviate $\mathrm{Ker}(f)=[0,0_{0,A},f_I^{\flat}]$ to $\mathrm{Ker}(f)=0$, because $[0,k,\varepsilon]=[0,k^{\prime},\varepsilon^{\prime}]$ for any $(k,\varepsilon)$ and $(k^{\prime},\varepsilon^{\prime})$. Dually, we abbreviate $\mathrm{Cok}(f)=[0,0_{A,0},f_I^{\sharp}]$ to $\mathrm{Cok}(f)=0$.

By using condition {\rm (A3)} of Definition \ref{LocSCG}, we can show the following easily:

\begin{example}
{\rm (1)} For any $A\in\mathbf{S}^0$, $\mathrm{Ker}(0_{A,0}:A\rightarrow0)=[A,\mathrm{id}_A,\mathrm{id}_0]$.

{\rm (2)} For any $A\in\mathbf{S}^0$, $\mathrm{Cok}(0_{0,A}:0\rightarrow A)=[A,\mathrm{id}_A,\mathrm{id}_0]$.
\end{example}

\begin{caution}
{\rm (1)} $\mathrm{Ker}(0_{0,A}:0\rightarrow A)$ need not be equivalent to $0$. Indeed, in the case of $\mathrm{SCG}$, for any symmetric categorical group $\mathbb{G}$, $\mathrm{Ker}(0_{0,\mathbb{G}}:0\rightarrow\mathbb{G})$ is equivalent to an important invariant $\pi_1(\mathbb{G)}[0]$.

{\rm (2)} $\mathrm{Cok}(0_{A,0}:A\rightarrow0)$ need not be equivalent to $0$ either. In the case of $\mathrm{SCG}$, $\mathrm{Cok}(0_{\mathbb{G},0}:\mathbb{G}\rightarrow0)$ is equivalent to $\pi_0(\mathbb{G)}[1]$.
\end{caution}

\begin{remark}\label{B2B3Remark}
The precise meaning of condition {\rm (B2)} in Definition \ref{RelExCat} is that, for any 1-cell $f\in\mathbf{S}^1(A,B)$ and its cokernel $[\mathrm{Cok}(f),\mathrm{cok}(f),\pi_f]$,
$f$ is faithful if and only if $\mathrm{Ker}(\mathrm{cok}(f))=[A,f,\pi_f]$. Similarly for condition {\rm (B3)}.
\end{remark}

\paragraph{Relative $($co-$)$kernel and first properties of a relatively exact 2-category}\quad\\

Throughout this subsection, $\mathbf{S}$ is a relatively exact 2-category.

\begin{definition}
For any diagram in $\mathbf{S}$
\begin{equation}
\xy
(-8,0)*+{A}="4";
(8,0)*+{B}="6";
(24,0)*+{C}="8";
{\ar_{f} "4";"6"};
{\ar_{g} "6";"8"};
{\ar@/^1.65pc/^{0} "4";"8"};
{\ar@{=>}_<<<{\varphi} (8,2)*{};(8,6)*{}};
\endxy
,
\label{ElementaryDiagram}
\end{equation}
its relative kernel $(\mathrm{Ker}(f,\varphi),\mathrm{ker}(f,\varphi),\varepsilon_{(f,\varphi)})$ is defined as follows $($we abbreviate $\mathrm{ker}(f,\varphi)$ to $k(f,\varphi)):$\\
{\rm (a)} $\mathrm{Ker}(f,\varphi)\in\mathbf{S}^0$, $k(f,\varphi)\in\mathbf{S}^1(\mathrm{Ker}(f,\varphi),A)$, $\varepsilon_{(f,\varphi)}\in\mathbf{S}^2(k(f,\varphi)\circ f,0)$.\\
{\rm (b0)} $($compatibility of the 2-cells$)$

$\varepsilon_{(f,\varphi)}$ is compatible with $\varphi$ i.e. $(k(f,\varphi)\circ\varphi)\cdot k(f,\varphi)_I^{\sharp}=(\varepsilon_{(f,\varphi)}\circ g)\cdot g_I^{\flat}$.
\\
\noindent{\rm (b1)} $($existence of a factorization$)$

For any $K\in\mathbf{S}^0$, $k\in\mathbf{S}^1(K,A)$ and $\varepsilon\in\mathbf{S}^2(k\circ f,0)$ which are compatible with $\varphi$, i.e. $(k\circ\varphi)\cdot k_I^{\sharp}=(\varepsilon\circ g)\cdot g_I^{\flat}$, there exist $\underline{k}\in\mathbf{S}^1(K,\mathrm{Ker}(f,\varphi))$ and $\underline{\varepsilon}\in\mathbf{S}^2(\underline{k}\circ k(f,\varphi),k)$ such that $(\underline{\varepsilon}\circ f)\cdot\varepsilon=(\underline{k}\circ\varepsilon_{(f,\varphi)})\cdot(\underline{k})_I^{\sharp}$.
\[
\xy
(-29,-9)*+{\mathrm{Ker}(f,\varphi)}="0";
(-24,8)*+{K}="2";
(-8,0)*+{A}="4";
(8,0)*+{B}="6";
(24,0)*+{C}="8";
{\ar_<<<<<<{k(f,\varphi)} "0";"4"};
{\ar^{k} "2";"4"};
{\ar_{\underline{k}} "2";"0"};
{\ar_{f} "4";"6"};
{\ar_{g} "6";"8"};
{\ar@/^1.25pc/^{0} "2";"6"};
{\ar@/_1.25pc/_{0} "0";"6"};
{\ar@/^1.65pc/^{0} "4";"8"};
{\ar@{=>}_<<<{\varphi} (8,2)*{};(8,6)*{}};
{\ar@{=>}_<<<{\varepsilon} (-8,2)*{};(-7,7)*{}};
{\ar@{=>}^<<<<{\varepsilon_{(f,\varphi)}} (-8,-2)*{};(-7,-7)*{}};
{\ar@{=>}_<<<<<<<<{\underline{\varepsilon}} (-25,-5)*{};(-18,3)*{}};
\endxy
\]
\noindent{\rm (b2)} $($uniqueness of the factorization$)$

For any factorizations $(\underline{k},\underline{\varepsilon})$ and $(\underline{k}^{\prime},\underline{\varepsilon}^{\prime})$ which satisfy {\rm (b1)}, there exists a unique 2-cell $\xi\in\mathbf{S}^2(\underline{k},\underline{k}^{\prime})$ such that $(\xi\circ k(f,\varphi))\cdot\underline{\varepsilon}^{\prime}=\underline{\varepsilon}.$
\end{definition}

\begin{remark}
By its universality, the relative kernel of $(f,\varphi)$ is unique up to\ an equivalence. We write this equivalence class $[\mathrm{Ker}(f,\varphi),k(f,\varphi),\varepsilon_{(f,\varphi)}]$.
\end{remark}

\begin{definition}
Let $\mathbf{S}$ be a relatively exact 2-category. For any diagram $(\ref{ElementaryDiagram})$ in $\mathbf{S}$, its relative cokernel $(\mathrm{Cok}(g,\varphi),\mathrm{cok}(g,\varphi),\pi_{(g,\varphi)})$ is defined dually by universality. We abbreviate $\mathrm{cok}(g,\varphi)$ to $c(g,\varphi)$, and write the equivalence class of the relative cokernel $[\mathrm{Cok}(g,\varphi),c(g,\varphi),\pi_{(g,\varphi)}]$.
\end{definition}

\begin{caution}
In the rest of this paper, $\mathbf{S}$ denotes a relatively exact 2-category, unless otherwise specified. In the following propositions and lemmas, we often omit the statement and the proof of their duals. Each term should be replaced by its dual. For example, kernel and cokernel, faithfulness and cofaithfulness are mutually dual.

\end{caution}

\begin{remark}\label{3Rem}
By using condition {\rm (A3)} of Definition \ref{LocSCG}, we can show the following easily. (These are also corollaries of Proposition \ref{KerAndRelKer}.)\\
{\rm (1)} $\mathrm{Ker}(f,f_I^{\sharp})=\mathrm{Ker}(f)$ (and thus the ordinary kernel can be regarded as a relative kernel).
\[
\xy
(-8,0)*+{A}="4";
(8,0)*+{B}="6";
(24,0)*+{0}="8";
{\ar_{f} "4";"6"};
{\ar_{0} "6";"8"};
{\ar@/^1.65pc/^{0} "4";"8"};
{\ar@{=>}_{f^{\sharp}_I} (8,2)*{};(8,6)*{}};
\endxy
\]
\noindent{\rm (2)} $\mathrm{ker}(f,\varphi)$ is faithful.
\end{remark}

\begin{lemma}\label{KernelLemma}
Let $f\in\mathbf{S}^1(A,B)$ and take its kernel $[\mathrm{Ker}(f),k(f),\varepsilon_f]$. If $K\in\mathbf{S}^0$, $k\in\mathbf{S}^1(K,\mathrm{Ker}(f))$ and $\sigma\in\mathbf{S}^2(k\circ k(f),0)$
\[
\xy
(-30,0)*+{K}="2";
(-10,0)*+{\mathrm{Ker}(f)}="4";
(10,0)*+{A}="6";
(30,0)*+{B}="8";
{\ar_{k} "2";"4"};
{\ar_{k(f)} "4";"6"};
{\ar_{f} "6";"8"};
{\ar@/^1.65pc/^{0} "4";"8"};
{\ar@/_1.65pc/_{0} "2";"6"};
{\ar@{=>}_{\varepsilon_f} (10,2)*{};(10,6)*{}};
{\ar@{=>}_{\sigma} (-10,-2)*{};(-10,-6)*{}};
\endxy
\]
is compatible with $\varepsilon_{f}$, i.e. if $\sigma$ satisfies
\begin{equation}
(\sigma\circ f)\cdot f_I^{\flat}=(k\circ\varepsilon_f)\cdot k_I^{\sharp},
\label{firstequality2}
\end{equation}
then there exists a unique 2-cell $\zeta\in\mathbf{S}^2(k,0)$ such that $\sigma=(\zeta\circ k(f))\cdot k(f)_I^{\flat}$.
\end{lemma}
\begin{proof}
By $(\ref{firstequality2})$, $\sigma:k\circ k(f)\Longrightarrow0$ is a factorization compatible with $\varepsilon_f$ and $f_I^{\flat}$.
On the other hand, by Remark \ref{Contra}, $k(f)_I^{\flat}:0\circ k(f)\Rightarrow0$ is also a factorization compatible with $\varepsilon_f$, $f_I^{\flat}$.
So, by the universality of the kernel, there exists a unique 2-cell $\zeta\in\mathbf{S}^2(k,0)$ such that $\sigma=(\zeta\circ k(f))\cdot k(f)_I^{\flat}$.
\end{proof}

It is easy to see that the same statement also holds for relative (co-)kernels. In any relatively exact 2-category, the relative (co-)kernel always exist. More precisely, the following proposition holds.
\begin{proposition}\label{KerAndRelKer}
Consider diagram $(\ref{ElementaryDiagram})$ in $\mathbf{S}$. By the universality of $\mathrm{Ker}(g)=[\mathrm{Ker}(g),\ell,\varepsilon]$, $f$ factors through $\ell$ uniquely up to an equivalence as $\underline{\varphi}:\underline{f}\circ\ell\Longrightarrow f$, where $\underline{f}\in\mathbf{S}^1(A,\mathrm{Ker}(g))$ and $\underline{\varphi}\in\mathbf{S}^2(\underline{f}\circ\ell,f):$
\[
(\underline{f}\circ\varepsilon)\cdot(\underline{f})_I^{\sharp}=(\underline{\varphi}\circ g)\cdot\varphi
\quad
\xy
(-16,-12)*+{\mathrm{Ker}(\underline{f})}="0";
(0,12)*+{\mathrm{Ker}(g)}="2";
(-8,0)*+{A}="4";
(8,0)*+{B}="6";
(24,0)*+{C}="8";
{\ar^{\ell} "2";"6"};
{\ar^{\underline{f}} "4";"2"};
{\ar_{k(\underline{f})} "0";"4"};
{\ar_{f} "4";"6"};
{\ar_{f} "4";"6"};
{\ar_{g} "6";"8"};
{\ar@/^0.50pc/^{0} "2";"8"};
{\ar@/_1.65pc/_{0} "4";"8"};
{\ar@/^1.80pc/^{0} "0";"2"};
{\ar@{=>}_<<<{\varphi} (8,-2)*{};(8,-6)*{}};
{\ar@{=>}_<<<{\underline{\varphi}} (0,6)*{};(0,2)*{}};
{\ar@{=>}^<<<{\varepsilon_{\underline{f}}} (-10,2)*{};(-13,4)*{}};
{\ar@{=>}_<<<{\varepsilon} (10,2)*{};(13,7)*{}};
\endxy
\]
Then we have $\mathrm{Ker}(f,\varphi)=[\mathrm{Ker}(\underline{f}),k(\underline{f}),\eta],$ where $\eta:=(k(\underline{f})\circ\underline{\varphi}^{-1})\cdot(\varepsilon_{\underline{f}}\circ\ell)\cdot\ell_I^{\flat}\in\mathbf{S}^2(k(\underline{f})\circ f,0)$. We abbreviate this to $\mathrm{Ker}(f,\varphi)=\mathrm{Ker}(\underline{f})$.
\end{proposition}
\begin{proof}
For any $K\in\mathbf{S}^0$, $k\in\mathbf{S}^1(K,A)$ and $\sigma\in\mathbf{S}^2(k\circ f,0)$ which are compatible with $\varphi$, i.e. $(\sigma\circ g)\cdot g_I^{\flat}=(k\circ\varphi)\cdot k_I^{\sharp}$, if we put
\[ \rho:=(k\circ\underline{\varphi})\cdot\sigma\in\mathbf{S}^2(k\circ\underline{f}\circ\ell,0), \]
then
$\rho$ is compatible with $\varepsilon$. By Lemma \ref{KernelLemma}, there exists a 2-cell $\zeta:k\circ\underline{f}\Rightarrow0$ such that $\rho=(\zeta\circ\ell)\cdot\ell_I^{\flat}$. So, by the universality of $\mathrm{Ker}(\underline{f})$, there exist $\underline{k}\in\mathbf{S}^1(K,\mathrm{Ker}(\underline{f}))$ and $\underline{\sigma}\in\mathbf{S}^2(\underline{k}\circ k(\underline{f}),k)$
such that $(\underline{\sigma}\circ\underline{f})\cdot\zeta=(\underline{k}\circ\varepsilon_{\underline{f}})\cdot(\underline{k})_I^{\sharp}$. Then, $\underline{\sigma}$ is compatible with $\sigma$ and $\eta$,
\[
\xy
(-29,-9)*+{\mathrm{Ker}(\underline{f})}="0";
(-24,8)*+{K}="2";
(-8,0)*+{A}="4";
(8,0)*+{B}="6";
{\ar_<<<<<<{k(\underline{f})} "0";"4"};
{\ar^{k} "2";"4"};
{\ar_{\underline{k}} "2";"0"};
{\ar_{f} "4";"6"};
{\ar@/^1.25pc/^{0} "2";"6"};
{\ar@/_1.25pc/_{0} "0";"6"};
{\ar@{=>}_<<<{\sigma} (-8,2)*{};(-7,7)*{}};
{\ar@{=>}^<<<<{\eta} (-8,-2)*{};(-7,-7)*{}};
{\ar@{=>}_<<<<<<<<{\underline{\sigma}} (-25,-5)*{};(-18,3)*{}};
\endxy
\]
and the existence of a factorization is shown. To show the uniqueness of the factorization, let $(\underline{k}^{\prime},\underline{\sigma}^{\prime})$ be another factorization which is compatible with $\sigma$, $\eta$, i.e. $(\underline{\sigma}^{\prime}\circ f)\cdot\sigma=(\underline{k}^{\prime}\circ\eta)\cdot(\underline{k}^{\prime})_I^{\sharp}$. Then, by using $\ \eta=(k(\underline{f})\circ\underline{\varepsilon}^{-1})\cdot(\varepsilon_{\underline{f}}\circ\ell)\cdot\ell_I^{\flat}$ and $\zeta\circ\ell=\rho\cdot\ell_I^{\flat-1}=(k\circ\underline{\varepsilon})\cdot\sigma\cdot\ell_I^{\flat-1}$, we can show $((\underline{\sigma}^{\prime}\circ\underline{f})\cdot\zeta)\circ\ell=((\underline{k}^{\prime}\circ\varepsilon_{\underline{f}})\cdot(\underline{k}^{\prime})_I^{\sharp})\circ\ell$. Since $\ell$ is faithful, we obtain $((\underline{\sigma}^{\prime}\circ\underline{f})\cdot\zeta)=(\underline{k}^{\prime}\circ\varepsilon_{\underline{f}})\cdot(\underline{k}^{\prime})_I^{\sharp}$. Thus, $\underline{\sigma}^{\prime}$ is compatible with $\zeta$ and $\varepsilon_{\underline{f}}$.
By the universality of $\mathrm{Ker}(\underline{f})$, there exists a 2-cell $\xi\in\mathbf{S}^2(\underline{k},\underline{k}^{\prime})$ such that $(\xi\circ k(\underline{f}))\cdot\underline{\sigma}^{\prime}=\underline{\sigma}$. Uniqueness of such $\xi\in\mathbf{S}^2(\underline{k},\underline{k}^{\prime})$ follows from the faithfulness of $k(\underline{f})$.
\end{proof}

\begin{proposition}\label{KerFullFaith}
Let $f\in\mathbf{S}^1(A,B)$, $g\in\mathbf{S}^1(B,C)$
and suppose $g$ is fully faithful. Then, $\mathrm{Ker}(f\circ g)=[\mathrm{Ker}(f),k(f),(\varepsilon_f\circ g)\cdot g_I^{\flat}]$. We abbreviate this to $\mathrm{Ker}(f\circ g)=\mathrm{Ker}(f)$.
\end{proposition}
\begin{proof}
Since $g$ is fully faithful, for any $K\in\mathbf{S}^0,k\in\mathbf{S}^1(K,A)$ and $\sigma\in\mathbf{S}^2(k\circ f\circ g,0)$,
there exists $\rho\in\mathbf{S}^2(k\circ f,0)$ such that $\sigma=(\rho\circ g)\cdot g_I^{\flat}$. And by the universality of $\mathrm{Ker}(f)$, there are $\underline{k}\in\mathbf{S}^1(K,\mathrm{Ker}(f))$ and $\underline{\sigma}\in\mathbf{S}^2(\underline{k}\circ k(f),k)$ such that $(\underline{\sigma}\circ f)\cdot\rho=(\underline{k}\circ\varepsilon_f)\cdot(\underline{k})_I^{\sharp}$.
Then, it can be easily seen that $\underline{\sigma}$ is compatible with $\sigma$ and $(\varepsilon_f\circ g)\cdot g_I^{\flat}$:
\[ (\underline{\sigma}\circ f\circ g)\cdot\sigma=(\underline{k}\circ((\varepsilon_f\circ g)\cdot g_I^{\flat}))\cdot(\underline{k})_I^{\sharp}
\quad
\xy
(-29,-9)*+{\mathrm{Ker}(f)}="0";
(-24,8)*+{K}="2";
(-8,0)*+{A}="4";
(8,0)*+{C}="6";
{\ar_<<<<<<{k(f)} "0";"4"};
{\ar^{k} "2";"4"};
{\ar_{\underline{k}} "2";"0"};
{\ar_{f \circ g} "4";"6"};
{\ar@/^1.25pc/^{0} "2";"6"};
{\ar@/_1.25pc/_{0} "0";"6"};
{\ar@{=>}_<<<{\sigma} (-8,2)*{};(-7,7)*{}};
{\ar@{=>}^<<<<{} (-8,-2)*{};(-7,-7)*{}};
{\ar@{=>}_<<<<<<<<{\underline{\sigma}} (-25,-5)*{};(-18,3)*{}};
{\ar@{..}^<<<<{} (-7.5,-5)*{};(-7,-16)*{_{(\varepsilon_f\circ g)\cdot g^{\flat}_I}}};
\endxy
\]
Thus we obtain a desired factorization. To show the uniqueness of the factorization, let $(\underline{k}^{\prime},\underline{\sigma}^{\prime})$ be another factorization of $k$ which satisfies
\[ (\underline{\sigma}^{\prime}\circ f\circ g)\cdot\sigma=(\underline{k}^{\prime}\circ((\varepsilon_f\circ g)\cdot g_I^{\flat}))\cdot(\underline{k}^{\prime})_I^{\sharp} \]
Then, we can show $\underline{\sigma}^{\prime}$ is compatible with $\rho$ and $\varepsilon_f$. By the universality of $\mathrm{Ker}(f)$, there exists a 2-cell $\xi\in\mathbf{S}^2(\underline{k},\underline{k}^{\prime})$ such that $(\xi\circ k(f))\cdot\underline{\sigma}^{\prime}=\underline{\sigma}$. Uniqueness of such $\xi$ follows from the faithfulness of $k(f)$.
\end{proof}

By definition, $f\in\mathbf{S}^1(A,B)$ is faithful (resp. fully faithful) if and only if \ $-\circ f:\mathbf{S}^2(g,h)\rightarrow\mathbf{S}^2(g\circ f,h\circ f)$ is injective (resp. bijective) for any $K\in\mathbf{S}^0$ and any $g,h\in\mathbf{S}^1(K,A)$. Concerning this, we have the following lemma.

\begin{lemma}\label{LemLem}
Let $f\in\mathbf{S}^1(A,B)$.\\
{\rm (1)} $f$ is faithful if and only if for any $K\in\mathbf{S}^0$ and $k\in\mathbf{S}^1(K,A)$,
\[ -\circ f:\mathbf{S}^2(k,0)\rightarrow\mathbf{S}^2(k\circ f,0\circ f)\ \text{is injective}. \]
\noindent{\rm (2)} $f$ is fully faithful if and only if for any $K\in\mathbf{S}^0$ and $k\in\mathbf{S}^1(K,A)$,
\[ -\circ f:\mathbf{S}^2(k,0)\rightarrow\mathbf{S}^2(k\circ f,0\circ f)\ \text{is bijective}. \]
\end{lemma}
\begin{proof}
By Lemma \ref{CriticalLemma}, we have the following commutative diagram for any $g,h\in\mathbf{S}^1(K,A)$:
\[
\xy
(-16,10)*+{\mathbf{S}^2(g,h)}="0";
(16,10)*+{\mathbf{S}^2(g\otimes h^{\ast},0)}="2";
(-28,-4)*+{\mathbf{S}^2(g\circ f,h\circ f)}="4";
(0,-28)*+{\mathbf{S}^2((g\circ f)\otimes(h\circ f)^{\ast},0)}="6";
(28,-4)*+{\mathbf{S}^2((g\otimes h^{\ast})\circ f,0\circ f)}="8";
(0,12)*+{}="10";
{\ar_{\text{bij.}}^{\Phi_{g,h}} "0";"2"};
{\ar_{-\circ f} "0";"4"};
{\ar^{-\circ f} "2";"8"};
{\ar@{}|\circlearrowright"6";"10"};
{\ar_{\text{bij.}}^{\Phi_{g \circ f,h\circ f}} "4";"6"};
{\ar^{\text{bij.}}_{\Theta_{g,h}^{f^{\flat}}} "8";"6"};
{\ar^{-\circ f} "2";"8"};
\endxy
\]
So we have
\begin{eqnarray*}
&-\circ f:\mathbf{S}^2(g,h)\rightarrow\mathbf{S}^2(g\circ f,h\circ f)\ \text{is injective (resp.bijective)}\\
\Leftrightarrow&-\circ f:\mathbf{S}^2(g\otimes h^{\ast},0)\rightarrow\mathbf{S}^2((g\otimes h^{\ast})\circ f,0\circ f)\ \text{is injective (resp.bijective)}.
\end{eqnarray*}
\end{proof}

\begin{corollary}\label{CorOfLem1}
For any $f\in\mathbf{S}^1(A,B)$, $f$ is faithful if and only if the following condition is satisfied$:$
\begin{equation}
\alpha\circ f=\mathrm{id}_{0\circ f}\Rightarrow\alpha=\mathrm{id}_0\ (\forall K\in\mathbf{S}^0,\forall\alpha\in\mathbf{S}^2(0_{K,A},0_{K,A}))
\label{CorCor}
\end{equation}
\end{corollary}
\begin{proof}
If $f$ is faithful, $(\ref{CorCor})$ is trivially satisfied, since we have $\mathrm{id}_{0\circ f}=\mathrm{id}_0\circ f$. To show the other implication, by Lemma \ref{LemLem}, it suffices to show that $-\circ f:\mathbf{S}^2(k,0)\rightarrow\mathbf{S}^2(k\circ f,0\circ f)\ \text{is injective}$.
For any $\alpha_1,\alpha_2\in\mathbf{S}^2(k,0)$ which satisfy $\alpha_1$ $\circ f=\alpha_2\circ f$, we have $(\alpha_1^{-1}\cdot\alpha_2)\circ f=(\alpha_1\circ f)^{-1}\cdot(\alpha_2\circ f)=\mathrm{id}_{0\circ f}$.
From the assumption we obtain $\alpha_1^{-1}\cdot\alpha_2=\mathrm{id}_0$, i.e. $\alpha_1=\alpha_2.$
\end{proof}

The next corollary immediately follows from Lemma \ref{LemLem}.
\begin{corollary}\label{CorOfLemLem}
For any $f\in\mathbf{S}^1(A,B)$, $f$ is fully faithful if and only if for any $K\in\mathbf{S}^0$, $k\in\mathbf{S}^1(K,A)$, and any $\sigma\in\mathbf{S}^2(k\circ f,0)$, there exists unique $\tau\in\mathbf{S}^2(k,0)$ such that $\sigma=(\tau\circ f)\cdot f_I^{\flat}$.
\end{corollary}

\begin{corollary}\label{RevisionLabel1}
For any $f\in\mathbf{S}^1(A,B)$, the following are equivalent$:$\\
{\rm (1)} $f$ is fully faithful.\\
{\rm (2)} $\mathrm{Ker}(f)=0$.
\end{corollary}
\begin{proof}
\underline{{\rm (1)}$\Rightarrow${\rm (2)}}

Since $f$ is fully faithful, for any $k\in\mathbf{S}^1(K,A)$ and $\varepsilon\in\mathbf{S}^2(k\circ f,0)$, there exists a 2-cell $\underline{\varepsilon}\in\mathbf{S}^2(0_{K,A},k)$ such that $(\underline{\varepsilon}\circ f)=(0\circ f_I^{\flat})\cdot0_I^{\sharp}\cdot\varepsilon^{-1}=(0\circ f_I^{\flat})\cdot\varepsilon^{-1}$, and the existence of a factorization is shown. To show the uniqueness of the factorization, it suffices to show that for any other factorization $(\underline{k}^{\prime},\underline{\varepsilon}^{\prime})$ with $(\underline{\varepsilon}^{\prime}\circ f)\cdot\varepsilon=(\underline{k}^{\prime}\circ f_I^{\flat})\cdot(\underline{k}^{\prime})_I^{\sharp}$, the unique 2-cell $\tau\in\mathbf{S}^2(\underline{k}^{\prime},0)$ (see condition {\rm (a3-2)} in Definition \ref{LocSCG}) satisfies $(\tau\circ0)\cdot\underline{\varepsilon}=\underline{\varepsilon}^{\prime}$. Since $f$ is faithful, this is equivalent to $(\tau\circ0\circ f)\cdot(\underline{\varepsilon}\circ f)\cdot\varepsilon=(\underline{\varepsilon}^{\prime}\circ f)\cdot\varepsilon$, and this follows easily from \ $\tau\circ0=(\tau\circ0)\cdot0_I^{\sharp}=(\underline{k}^{\prime})_I^{\sharp}$ and $(\tau\circ0\circ f)\cdot(0\circ f_I^{\flat})=(\underline{k}^{\prime}\circ f_I^{\flat})\cdot(\tau\circ0)$. (see Corollary \ref{ContraCor}.)\\
\underline{{\rm (2)}$\Rightarrow${\rm (1)}}
Since $\mathrm{Ker}(f)=[0,0,f^{\flat}_I]$, for any $K\in \mathbf{S}^0$, $k\in\mathbf{S}^1(K,A)$ and any $\sigma\in\mathbf{S}^2(k\circ f,0)$, there exist $\underline{k}\in\mathbf{S}^1(K,0)$ and $\underline{\sigma}\in\mathbf{S}^2(\underline{k}\circ 0,k)$ such that $(\underline{\sigma}\circ f)\cdot \sigma=(\underline{k}\circ f^{\flat}_I)\cdot (\underline{k})^{\sharp}_I$. Thus $\tau:=\underline{\sigma}^{-1}\cdot\underline{k}^{\sharp}_I$ satisfies $\sigma=(\tau\circ f)\cdot f^{\flat}_I$.
If there exists another $\tau^{\prime}\in \mathbf{S}^2(k,0)$ satisfying $\sigma=(\tau^{\prime}\circ f)\cdot f^{\flat}_I$, then by the universality of the kernel, there exists $\upsilon\in\mathbf{S}^2(\underline{k},0)$ such that $(\upsilon\circ 0)\cdot\tau^{\prime -1}=\underline{\tau}$. Since $\upsilon\circ 0=\underline{k}^{\sharp}_I$ by $(\ref{eqtn})$, we obtain $\tau=\tau^{\prime}$. Thus $\tau$ is uniquely determined.

\end{proof}

\begin{proposition}\label{final1}
For any $f\in\mathbf{S}^1(A,B)$, the following are equivalent.\\
{\rm (1)} $f$ is an equivalence.\\
{\rm (2)} $f$ is cofaithful and fully faithful.\\
{\rm (3)} $f$ is faithful and fully cofaithful.
\end{proposition}
\begin{proof}
Since {\rm (1)}$\Leftrightarrow${\rm (3)} is the dual of {\rm (1)}$\Leftrightarrow${\rm (2)}, we show only {\rm (1)}$\Leftrightarrow${\rm (2)}.\\
{\rm (1)}$\Rightarrow${\rm (2)} : trivial.\\
{\rm (2)}$\Rightarrow${\rm (1)} : Since $f$ is cofaithful, we have $f=\mathrm{cok}(\mathrm{ker}(f))$, $\mathrm{Cok}(k(f))=[B,f,\varepsilon_f]$. On the other hand, since $f$ is fully faithful, we have $\mathrm{Ker}(f)=[0,0,f_I^{\flat}]$, and so we have $\mathrm{Cok}(k(f))=[A,\mathrm{id}_A,\mathrm{id}_0]$. And by the uniqueness (up to an equivalence) of the cokernel, there is an equivalence from $A$ to $B$, which is equivalent to $f$. Thus, $f$ becomes an equivalence.
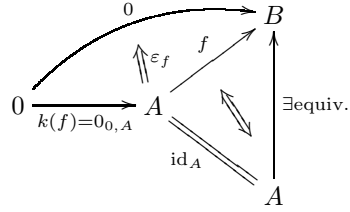
\[
\xy
(-18,0)*+{0}="2";
(16,-12)*+{A}="4";
(0,0)*+{A}="6";
(16,12)*+{B}="8";
{\ar@{=}_>>>>>>>>>{\mathrm{id}_A} "6";"4"};
{\ar_{k(f)=0_{0,A}} "2";"6"};
{\ar^{f} "6";"8"};
{\ar_{\exists\text{equiv.} } "4";"8"};
{\ar@/^1.25pc/^{0} "2";"8"};
{\ar@{=>}_{\varepsilon_f} (-1,3)*{};(-2,8)*{}};
{\ar@{<=>} (9,2)*{};(13,-4)*{}};
\endxy
\]
\end{proof}

\begin{lemma}\label{final3}
Let $f:A\rightarrow B$ be a faithful 1-cell in $\mathbf{S}$. Then, for any $K\in\mathbf{S}^0$ and $k\in\mathbf{S}^1(K,0)$, we have $\mathbf{S}^2(k\circ0_{0,\mathrm{Ker}(f)},0_{K,\mathrm{Ker}(f)})=\{k_I^{\sharp}\}$.
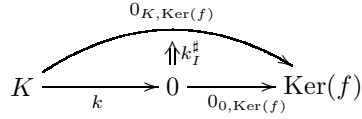
\[
\xy
(-20,0)*+{K}="10";
(0,0)*+{0}="12";
(20,0)*+{\mathrm{Ker}(f)}="22";
{\ar_{k} "10";"12"};
{\ar_{0_{0,\mathrm{Ker}(f)}} "12";"22"};
{\ar@/^1.80pc/^{0_{K,\mathrm{Ker}(f)}} "10";"22"};
{\ar@{=>}_{k^{\sharp}_I} (0,3)*{};(0,6)*{}};
\endxy
\]
\end{lemma}

\begin{proof}
For any $\sigma\in\mathbf{S}^2(k\circ0_{0,\mathrm{Ker}(f)},0_{K,\mathrm{Ker}(f)})$, we can show $((\sigma\circ k(f))\cdot k(f)_I^{\flat})\circ f=((k\circ k(f)_I^{\flat})\cdot k_I^{\sharp})\circ f$. By the faithfulness of $f$, we have $(\sigma\circ k(f))\cdot k(f)_I^{\flat}=(k\circ k(f)_I^{\flat})\cdot k_I^{\sharp}$.
Thus, we have $\sigma\circ k(f)=k_I^{\sharp}\circ k(f)$. By the faithfulness of $k(f)$, we obtain $\sigma=k_I^{\sharp}$.
\end{proof}

\begin{corollary}\label{final4}
$f:A\rightarrow B$ is faithful if and only if $\mathrm{Ker}(0_{0,A},f_I^{\flat})=0$.
\end{corollary}
\begin{proof}
Since there is a factorization diagram with $(0_{0,\mathrm{Ker}(f)}\circ\varepsilon_f)\cdot(0_{0,\mathrm{Ker}(f)})_I^{\sharp}=(k(f)_I^{\flat}\circ f)\cdot f_I^{\flat}$
\[
\xy
(19,0)*+{B}="0";
(22,-1)*+{,}="1";
(0,0)*+{A}="2";
(-20,12)*+{\mathrm{Ker}(f)}="4";
(-20,-12)*+{0}="6";
{\ar_{f} "2";"0"};
{\ar^{k(f)} "4";"2"};
{\ar|*+{_{0_{0,A}}} "6";"2"};
{\ar@/^0.60pc/^{0_{0,\mathrm{Ker}(f)}} "6";"4"};
{\ar@/^1.20pc/^<<<<<<<<<<{0} "4";"0"};
{\ar@/_1.20pc/_<<<<<<<<<<{0} "6";"0"};
{\ar@{=>}_{\varepsilon_f} (2,3)*{};(4,7)*{}};
{\ar@{=>}_{f^{\flat}_I}(2,-3)*{};(4,-7)*{}};
{\ar@{=>}_{k(f)^{\flat}_I} (-14,4)*{};(-10,-2)*{}};
\endxy
\]
(see {\rm (a3-2)} in Definition \ref{LocSCG}) we have $\mathrm{Ker}(0_{0,A},f_I^{\flat})=\mathrm{Ker}(0_{0,\mathrm{Ker}(f)})$ by Proposition \ref{KerAndRelKer}. So, it suffices to show $\mathrm{Ker}(0_{0,\mathrm{Ker}(f)})=0$. For any $K\in\mathbf{S}^0$ and $k\in\mathbf{S}^1(K,0)$, we have $\mathbf{S}^2(k\circ0_{0,\mathrm{Ker}(f)},0_{K,\mathrm{Ker}(f)})=\{k_I^{\sharp}\}$ by the Lemma \ref{final3}. So
$0_{0,\mathrm{Ker}(f)}$ becomes fully faithful, and thus $\mathrm{Ker}(0_{0,\mathrm{Ker}(f)})=0$.

Conversely, assume $\mathrm{Ker}(0_{0,A},f_I^{\flat})=0$. For any $K\in \mathbf{S}^0$ and $\alpha\in\mathbf{S}^2(0_{K,A},0_{K,A})$ satisfying $\alpha\circ f=\mathrm{id}_{0\circ f}$, we show $\alpha=\mathrm{id}_0$ (Corollary \ref{CorOfLem1}).

By $\alpha\circ f=\mathrm{id}_{0\circ f}$, $\alpha$ is compatible with $f^{\flat}_I$ :
\[
\xy
(36,0)*+{B}="8";
(18,0)*+{A}="6";
(0,0)*+{0}="2";
(-36,12)*+{\mathrm{Ker}(0_{0,A},f_I^{\flat})}="0";
(-20,12)*+{=0}="1";
(-20,-12)*+{K}="4";
{\ar_{f} "6";"8"};
{\ar_{0_{0,A}} "2";"6"};
{\ar^{0_{K,0}} "4";"2"};
{\ar_{\mathrm{id}_0} "1";"2"};
{\ar@/^1.20pc/^{0} "1";"6"};
{\ar@/_1.20pc/_{0} "4";"6"};
{\ar@/^1.80pc/^{0} "2";"8"};
{\ar@{=>}^{\mathrm{id}_0} (1,3)*{};(2,8)*{}};
{\ar@{=>}_{\alpha}(1,-3)*{};(2,-8)*{}};
{\ar@{=>}_{f^{\flat}_I} (18,2)*{};(18,6)*{}};
\endxy
\]
So there exist $k\in \mathbf{S}^1(K,0)$ and $\varepsilon\in\mathbf{S}^2(k\circ\mathrm{id}_0,0_{K,0})$ satisfying
\[ (\varepsilon\circ0_{0,A})\cdot\alpha=(k\circ \mathrm{id_0})\cdot k^{\sharp}_I. \]
Since $\varepsilon\circ0_{0,A}=k^{\sharp}_I$ by $(\ref{1-6})$ and $(\ref{1-14})$, we obtain $\alpha=\mathrm{id}_0$.
\end{proof}

In any relatively exact 2-category $\mathbf{S}$, the difference kernel of any pair of 1-cells $g,h:A\rightarrow B$ always exists. More precisely, we have the following proposition:

\begin{proposition}\label{ExistDiffKer}
For any $g,h\in\mathbf{S}^1(A,B)$, if we take the kernel $\mathrm{Ker}(g\otimes h^{\ast})=[\mathrm{Ker}(g\otimes h^{\ast}),k,\varepsilon]$ of $g\otimes h^{\ast}$
and put $\widetilde{\varepsilon}:=\Psi_{k\circ g,k\circ h}(\Theta_{g,h}^{k^{\sharp}}(\varepsilon\cdot k_I^{\sharp-1}))\in\mathbf{S}^2(k\circ g,k\circ h)$,
then $(\mathrm{Ker}(g\otimes h^{\ast}),k,\widetilde{\varepsilon})$ is the difference kernel of $g$ and $h$.
\end{proposition}
\begin{proof}
For any $K\in\mathbf{S}^0$ and $\ell\in\mathbf{S}^1(K,A)$,
there exists a natural isomorphism (Lemma \ref{CriticalLemma})
\[
\xy
(-28,4)*+{\mathbf{S}^2(\ell\circ(g\otimes h^{\ast}),0)}="0";
(28,4)*+{\mathbf{S}^2(\ell\circ g,\ell\circ h)}="2";
(-28,0)*+{\rotatebox{90}{$\in$}}="4";
(28,0)*+{\rotatebox{90}{$\in$}}="6";
(-28,-4)*+{\sigma}="8";
(28,-4)*+{\widetilde{\sigma}:=\Psi_{\ell\circ g,\ell\circ h}(\Theta^{\ell^{\sharp}}_{g,h} (\sigma\cdot\ell^{\sharp}_I)).}="10";
{\ar "0";"2"};
{\ar@{|->} "8";"10"};
\endxy
\]
So, to give a 2-cell $\sigma\in\mathbf{S}^2(\ell\circ(g\otimes h^{\ast}),0)$
is equivalent to give a 2-cell $\widetilde{\sigma}\in\mathbf{S}^2(\ell\circ g,\ell\circ h)$.
And, by using Remark \ref{Contra} and Corollary \ref{ContraCor}, the usual compatibility argument shows the proposition.
\end{proof}

In any relatively exact 2-category $\mathbf{S}$, the pullback of any pair of morphisms $f_i:A_i\rightarrow B$ $(i=1,2)$ always exists. More precisely, we have the following proposition:

\begin{proposition}\label{ExistPullback}
For any $f_i\in\mathbf{S}^1(A_i,B)\ (1=1,2)$, if we take the product of $A_1$ and $A_2$ $(A_1\times A_2,p_1,p_2)$, and take the difference kernel $(D,d,\varphi)$ of $p_1\circ f_1$ and $p_2\circ f_2$
\[
\xy
(-24,0)*+{D}="0";
(0,0)*+{A_1\times A_2}="2";
(20,0)*+{B}="4";
{\ar^{d} "0";"2"};
{\ar@/^1.25pc/^{p_1\circ f_1} "2";"4"};
{\ar@/_1.25pc/_{p_2\circ f_2} "2";"4"};
\endxy
\quad
\xy
(-12,0)*+{D}="0";
(0,8)*+{A_1}="2";
(0,-8)*+{A_2}="4";
(12,0)*+{B}="6";
{\ar^{d\circ p_1} "0";"2"};
{\ar_{d\circ p_2} "0";"4"};
{\ar^{f_1} "2";"6"};
{\ar_{f_2} "4";"6"};
{\ar@{=>}^{\varphi} (0,4)*{};(0,-4)*{}};
\endxy
,
\]
then, $(D,d\circ p_1,d\circ p_2,\varphi)$ is the pullback of $f_1$ and $f_2$.
\end{proposition}
\begin{proof}[of condition {\rm (b1)} (in Definition \ref{DefPullback})]
For any $X\in\mathbf{S}^0$, $g_i\in\mathbf{S}^1(X,A_i)$ $(i=1,2)$ and $\eta\in\mathbf{S}^2(g_1\circ f_1,g_2\circ f_2)$, by the universality of $A_1\times A_2$, there exist $g\in\mathbf{S}^1(X,A_1\times A_2)$ and $\xi_i\in\mathbf{S}^2(d\circ p_i,g_i)\ (i=1,2)$.
Applying the universality of the difference kernel to the 2-cell
\begin{equation}
\zeta:=(\xi_1\circ f_1)\cdot\eta\cdot(\xi_2^{-1}\circ f_2)\in\mathbf{S}^2(g\circ p_1\circ f_1,g\circ p_2\circ f_2),
\label{firstequality6}
\end{equation}
we see there exist $\underline{g}\in\mathbf{S}^1(X,D)$ and $\underline{\zeta}\in\mathbf{S}^2(\underline{g}\circ d,g)$
\begin{equation}
\xy
(-28,0)*+{D}="0";
(0,0)*+{A_1\times A_2}="2";
(32,0)*+{B}="4";
(-20,10)*+{X}="6";
{\ar_{d} "0";"2"};
{\ar@/^0.50pc/^{p_1\circ f_1} "2";"4"};
{\ar@/_0.50pc/_{p_2\circ f_2} "2";"4"};
{\ar^{g} "6";"2"};
{\ar_{\underline{g}} "6";"0"};
{\ar@{=>}^{\underline{\zeta}} (-20,2)*{};(-14,4)*{}};
\endxy
\label{firstequality9}
\end{equation}
such that
\begin{equation}
(\underline{g}\circ\varphi)\cdot(\underline{\zeta}\circ p_2\circ f_2)=(\underline{\zeta}\circ p_1\circ f_1)\cdot\zeta.
\label{firstequality7}
\end{equation}

By $(\ref{firstequality6})$ and $(\ref{firstequality7})$, we have $(\underline{g}\circ\varphi)\cdot(((\underline{\zeta}\circ p_2)\cdot\xi_2)\circ f_2)=(((\underline{\zeta}\circ p_1)\cdot\xi_1)\circ f_1)\cdot\eta$, and thus condition {\rm (b1)} is satisfied.
\[
\xy
(-10,0)*+{D}="0";
(0,10)*+{A_1}="2";
(0,-10)*+{A_2}="4";
(10,)*+{B}="6";
(-32,0)*+{X}="8";
{\ar^{d\circ p_1} "0";"2"};
{\ar_{d \circ p_2} "0";"4"};
{\ar_{f_2} "4";"6"};
{\ar^{f_1} "2";"6"};
{\ar_{\underline{g}} "8";"0"};
{\ar@/^1.20pc/^{g_1} "8";"2"};
{\ar@/_1.20pc/_{g_2} "8";"4"};
{\ar@{=>}^{\varphi} (0,5)*{};(0,-5)*{}};
{\ar@{=>} (-12,2)*{};(-16,8)*{}};
{\ar@{=>} (-12,-2)*{};(-16,-8)*{}};
{\ar@{..}@/_1.60pc/^{} (-14,6)*{};(-15,20)*{_{(\underline{\zeta}\circ p_1)\cdot\xi_1}}};
{\ar@{..}@/^1.60pc/^{} (-14,-6)*{};(-15,-20)*{_{(\underline{\zeta}\circ p_2)\cdot\xi_2}}};
\endxy
\]
\noindent\underline{proof of condition (b2)}

If we take $h\in\mathbf{S}^1(X,D)$ and $\eta_i\in\mathbf{S}^2(h\circ d\circ p_i,g_i)\ (i=1,2)$
which satisfy $(h\circ\varphi)\cdot(\eta_2\circ f_2)=(\eta_1\circ f_1)\cdot\eta$, then by the universality of $A_1\times A_2$, there exists a unique 2-cell $\kappa\in\mathbf{S}^2(h\circ d,g)$ such that
\begin{equation}
(\kappa\circ p_i)\cdot\xi_i=\eta_i\ (i=1,2).\label{firstequality8}
\end{equation}
Then, $\kappa$ becomes compatible with $\varphi$ and $\zeta$, i.e. $(h\circ\varphi)\cdot(\kappa\circ p_2\circ f_2)=(\kappa\circ p_1\circ f_1)\cdot\zeta$.
So, comparing this with factorization $(\ref{firstequality9})$, by the universality of the difference kernel, we see there exists a unique 2-cell $\chi\in\mathbf{S}^2(h,\underline{g})$ which satisfies
\begin{equation}
(\chi\circ d)\cdot\underline{\zeta}=\kappa\label{firstequality10}
\end{equation}
Then we have $(\chi\circ d\circ p_i)\cdot(\underline{\zeta}\circ p_i)\cdot\xi_i=(\kappa\circ p_i)\cdot\xi_i=\eta_i\ (i=1,2)$. Thus $\chi$ is the desired 2-cell in condition {\rm (b2)}, and the uniqueness of such a $\chi$ follows from the uniqueness of $\kappa$ and $\chi$ which satisfy $(\ref{firstequality8})$ and $(\ref{firstequality10})$.
\end{proof}

By the universality of the pullback, we have the following remark:

\begin{remark}
\label{PreviousRemark}
Let
\begin{eqnarray}
\xy
(-14,0)*+{A_1\times_B A_2}="0";
(0,8)*+{A_2}="2";
(0,-8)*+{A_1}="4";
(12,0)*+{B}="6";
{\ar^{f^{\prime}_1} "0";"2"};
{\ar_{f^{\prime}_2} "0";"4"};
{\ar^{f_2} "2";"6"};
{\ar_{f_1} "4";"6"};
{\ar@{=>}^{\xi} (0,3)*{};(0,-3)*{}};
\endxy
\label{DocPullBack}
\end{eqnarray}
be a pull-back diagram. Then, for any $K\in\mathbf{S}^0$, $g,h\in\mathbf{S}^1(K,A_1\times_B A_2)$ and $\alpha,\beta\in\mathbf{S}^2(g,h)$, we have 
\[ \alpha\circ f_i^{\prime}=\beta\circ f_i^{\prime}\ (i=1,2)\Longrightarrow\alpha=\beta. \]
\end{remark}
\begin{proof}
To the diagram
\[
\xy
(-14,0)*+{K}="0";
(0,-8)*+{A_1}="2";
(0,8)*+{A_2}="4";
(12,0)*+{B}="6";
{\ar^{g\circ f^{\prime}_1} "0";"4"};
{\ar_{g\circ f^{\prime}_2} "0";"2"};
{\ar_{f_1} "2";"6"};
{\ar^{f_2} "4";"6"};
{\ar@{=>}_{g\circ\xi} (0,3)*{};(0,-3)*{}};
\endxy
,
\]
the following diagram gives a factorization which satisfies condition {\rm (b1)} in Definition \ref{DefPullback}.
\[
\xy
(-12,0)*+{A_1\times_B A_2}="0";
(0,10)*+{A_2}="2";
(0,-10)*+{A_1}="4";
(10,0)*+{B}="6";
(-32,0)*+{K}="8";
(-18,12)*+{}="10";
(-18,-12)*+{}="12";
{\ar^{f^{\prime}_1} "0";"2"};
{\ar_{f^{\prime}_2} "0";"4"};
{\ar_{f_1} "4";"6"};
{\ar^{f_2} "2";"6"};
{\ar_{g} "8";"0"};
{\ar@/^1.20pc/^{g\circ f^{\prime}_1} "8";"2"};
{\ar@/_1.20pc/_{g \circ f^{\prime}_2} "8";"4"};
{\ar@{=>}^{\xi} (0,5)*{};(0,-5)*{}};
{\ar@{}|\circlearrowright"0";"10"};
{\ar@{}|\circlearrowright"0";"12"};
\endxy
\]
Since each of $\mathrm{id}_g:g\Longrightarrow g$ and $\alpha\circ\beta^{-1}:g\Longrightarrow g$ gives a 2-cell which satisfies condition {\rm (b2)}, we have $\alpha\circ\beta^{-1}=\mathrm{id}$ by the uniqueness. Thus $\alpha=\beta$.
\end{proof}

\begin{proposition}\label{FaithLemma}\label{FaithProp}
$($See also Proposition \ref{final13}.$)$
Let $(\ref{DocPullBack})$ be a pull-back diagram. We have

{\rm (1)} $f_1$: faithful $\Rightarrow f_1^{\prime}$: faithful.

{\rm (2)} $f_1$: fully faithful $\Rightarrow f_1^{\prime}$: fully faithful.

{\rm (3)} $f_1$: cofaithful $\Rightarrow f_1^{\prime}$: cofaithful.
\end{proposition}
\begin{proof}
\noindent\underline{proof of {\rm (1)}} By Corollary \ref{CorOfLem1}, it suffices to show $\alpha\circ f_1^{\prime}=\mathrm{id}_{0\circ f_1^{\prime}}\Rightarrow\alpha=\mathrm{id}_0$ for any $K\in\mathbf{S}^0$ and $\alpha\in\mathbf{S}^2(0_{K,A_1\times_B A_2},0_{K,A_1\times_B A_2})$. Since $(0\circ\xi)\cdot(\alpha\circ f_2^{\prime}\circ f_1)=(\alpha\circ f_1^{\prime}\circ f_2)\cdot(0\circ\xi)=\mathrm{id}_{0\circ f_1^{\prime}\circ f_2}\cdot(0\circ\xi)=0\circ\xi$, we have $\alpha\circ f_2^{\prime}\circ f_1=\mathrm{id}_{0\circ f_2^{\prime}\circ f_1}=\mathrm{id}_{0\circ f_2^{\prime}}\circ f_1$. Since $f_1$ is faithful, we obtain $\alpha\circ f_2^{\prime}=\mathrm{id}_{0\circ f_2^{\prime}}$. Thus, we have $\alpha\circ f_i^{\prime}=\mathrm{id}_{0\circ f_i^{\prime}}=\mathrm{id}_0\circ f_i^{\prime}\ (i=1,2)$.
By Remark \ref{PreviousRemark}, this implies $\alpha=\mathrm{id}_0$.

\noindent\underline{proof of {\rm (2)}} By {\rm (1)}, $f_1^{\prime}$ is already faithful.
By Corollary \ref{CorOfLem1}, it suffices to show that for any $K\in\mathbf{S}^0$, $k\in\mathbf{S}^1(K,A_1\times_B A_2)$ and any $\sigma\in\mathbf{S}^2(k\circ f_1^{\prime},0)$, there exists a unique 2-cell $\kappa\in\mathbf{S}^2(k,0)$ such that $\sigma=(\kappa\circ f_1^{\prime})\cdot(f_1^{\prime})_I^{\flat}$. Since $f_1$ is fully faithful, for any $K\in\mathbf{S}^0$, $k\in\mathbf{S}^1(K,A_1\times_B A_2)$ and any $\sigma\in\mathbf{S}^2(k\circ f_1^{\prime},0)$, there exists $\tau\in\mathbf{S}^2(k\circ f_2^{\prime},0)$ such that $(\tau\circ f_1)\cdot(f_1)_I^{\flat}=(k\circ\xi^{-1})\cdot(\sigma\circ f_2)\cdot(f_2)_I^{\flat}$.
Then, for the diagram
\[
\xy
(-14,0)*+{K}="0";
(0,10)*+{A_1}="2";
(0,-10)*+{A_2}="4";
(14,0)*+{B}="6";
{\ar^{0} "0";"2"};
{\ar_{0} "0";"4"};
{\ar_{f_2} "4";"6"};
{\ar^{f_1} "2";"6"};
{\ar@{=>}|*+{_{(f_1)^{\flat}_I\cdot(f_2)^{\flat-1}_I}} (0,5)*{};(0,-5)*{}};
\endxy,
\]
each of the factorizations 
\[
\xy
(-14,0)*+{A_1\times_B A_2}="0";
(0,10)*+{A_2}="2";
(0,-10)*+{A_1}="4";
(14,0)*+{B}="6";
(-36,0)*+{K}="8";
{\ar^{f^{\prime}_1} "0";"2"};
{\ar_{f^{\prime}_2} "0";"4"};
{\ar_{f_1} "4";"6"};
{\ar^{f_2} "2";"6"};
{\ar_{k} "8";"0"};
{\ar@/^1.20pc/^{0} "8";"2"};
{\ar@/_1.20pc/_{0} "8";"4"};
{\ar@{=>}^{\xi} (0,4)*{};(0,-4)*{}};
{\ar@{=>}^{\sigma} (-16,3)*{};(-20,8)*{}};
{\ar@{=>}_{\tau} (-16,-3)*{};(-20,-8)*{}};
\endxy
\quad
\xy
(-14,0)*+{A_1\times_B A_2}="0";
(0,10)*+{A_2}="2";
(0,-10)*+{A_1}="4";
(14,0)*+{B}="6";
(-36,0)*+{K}="8";
{\ar^{f^{\prime}_1} "0";"2"};
{\ar_{f^{\prime}_2} "0";"4"};
{\ar_{f_1} "4";"6"};
{\ar^{f_2} "2";"6"};
{\ar_{0} "8";"0"};
{\ar@/^1.20pc/^{0} "8";"2"};
{\ar@/_1.20pc/_{0} "8";"4"};
{\ar@{=>}^{ \xi}  (0,4)*{};(0,-4)*{}};
{\ar@{=>}^{_{(f^{\prime}_1)^{\flat}_I}} (-16,3)*{};(-20,8)*{}};
{\ar@{=>}_{_{(f^{\prime}_2)^{\flat}_I}} (-16,-3)*{};(-20,-8)*{}};
\endxy
\]
satisfies condition {\rm (b1)} in Definition \ref{DefPullback}. So there exists a 2-cell $\kappa\in\mathbf{S}^2(k,0)$ such that $\sigma=(\kappa\circ f_1^{\prime})\cdot(f_1^{\prime})_I^{\flat}$. Uniqueness of such $\kappa$ follows from the faithfulness of $f_1^{\prime}$.\\
\underline{proof of {\rm (3)}} Let $(A_1\times A_2,p_1,p_2)$ be the product of $A_1$ and $A_2$. For $\mathrm{id}_{A_1}\in\mathbf{S}^1(A_1,A_1)$ and $0\in\mathbf{S}^1(A_1,A_2)$, by the universality of $A_1\times A_2$, there exist $i_1\in\mathbf{S}^1(A_1,A_1\times A_2)$, $\xi_1\in\mathbf{S}^2(i_1\circ p_1,\mathrm{id}_{A_1})$ and $\xi_2\in\mathbf{S}^2(i_2\circ p_2,0)$.
Similarly, there is a 1-cell $i_2\in\mathbf{S}^1(A_2,A_1\times A_2)$ such that there are equivalences $i_2\circ p_2\simeq \mathrm{id}_{A_2},\ i_2\circ p_1\simeq0$. If we put $t:=(p_1\circ f_1)\otimes(p_2\circ f_2)^{\ast}$,
then by Proposition \ref{ExistDiffKer} and \ref{ExistPullback}, we have $A_1\times_B A_2=\mathrm{Ker}(t)$. So we may assume $\mathrm{Ker}(t)=[A_1\times_BA_2,d,\varepsilon_t ]$ and $f_1^{\prime}=d\circ p_2$.
\[
\xy(-12,0)*+{A_1\times_B A_2}="4";
(12,0)*+{A_1\times A_2}="6";
(36,0)*+{B}="8";
{\ar_{d} "4";"6"};
{\ar_{t} "6";"8"};
{\ar@/^1.65pc/^{0} "4";"8"};
{\ar@{=>}_{\varepsilon_t} (12,2)*{};(12,6)*{}};
\endxy
\]
Since $i_1\circ t$ and $f_1$ are equivalent;
\[ i_1\circ t\simeq(i_1\circ p_1\circ f_1)\otimes(i_1\circ p_2\circ f_2^{\ast})\simeq(\mathrm{id}_{A_1}\circ f_1)\otimes(0\circ f_2^{\ast})\simeq f_1, \]
by the cofaithfulness of $f_1$, it follows that $t$ is cofaithful. Thus, we have $B=\mathrm{Cok}(\mathrm{ker}(t))$, i.e. $\mathrm{Cok}(d)=[B,t,\varepsilon_t]$.
By (the dual of) Corollary \ref{CorOfLem1}, it suffices to show $f_1^{\prime}\circ\alpha=\mathrm{id}_{f_1^{\prime}\circ 0}\Rightarrow\alpha=\mathrm{id}_0$ for any $C\in\mathbf{S}^0$ and any $\alpha\in\mathbf{S}^2(0_{A_2,C},0_{A_2,C})$. For the 2-cell $(d\circ p_2)_I^{\sharp}\in\mathbf{S}^2(d\circ p_2\circ 0_{A_2,C},0)$ (see the following diagram), by the universality of $\mathrm{Cok}(d)$, there exist $u\in\mathbf{S}^1(B,C)$ and $\gamma\in\mathbf{S}^2(t\circ u,p_2\circ0)$ such that $(d\circ\gamma)\cdot(d\circ p_2)_I^{\sharp}=(\varepsilon_t\circ u)\cdot u_I^{\flat}$. Thus, if we put $\gamma^{\prime}:=\gamma\cdot(p_2\circ\alpha)$, we have
\begin{eqnarray*}
(d\circ\gamma^{\prime})\cdot(d\circ p_2)_I^{\sharp}=(d\circ\gamma)\cdot(d\circ p_2\circ\alpha)\cdot(d\circ p_2)_I^{\sharp}\\
=(d\circ\gamma)\cdot(f_1^{\prime}\circ\alpha)\cdot(d\circ p_2)_I^{\sharp}
=(\varepsilon_t\circ u)\cdot u_I^{\flat}.
\end{eqnarray*}
So, $\gamma$ and $\gamma^{\prime}\in\mathbf{S}^2(t\circ u,p_2\circ0)$ give two factorization of $p_2\circ0$ compatible with $\varepsilon_t$ and $(d\circ p_2)_I^{\sharp}$. By the universality of $\mathrm{Cok}(d)=[B,t,\varepsilon_t]$, there exists a unique 2-cell $\beta\in\mathbf{S}^2(u,u)$ such that
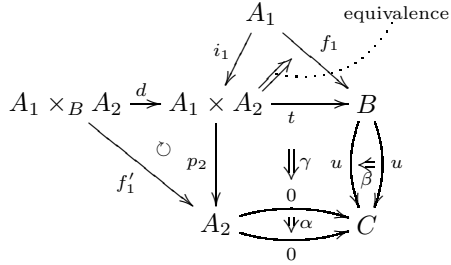
\begin{equation}
(t\circ\beta)\cdot\gamma=\gamma^{\prime}.\label{proppullback}
\end{equation}
\[
\xy
(-20,8)*+{A_1\times_B A_2}="0";
(0,8)*+{A_1\times A_2}="2";
(20,8)*+{B}="4";
(0,-8)*+{A_2}="6";
(20,-8)*+{C}="8";
(-14,-4)*+{}="10";
(6,20)*+{A_1}="12";
{\ar^{d} "0";"2"};
{\ar_{i_1} "12";"2"};
{\ar^{f_1} "12";"4"};
{\ar_{p_2} "2";"6"};
{\ar_{f^{\prime}_1} "0";"6"};
{\ar_{t} "2";"4"};
{\ar@/^0.50pc/^{u} "4";"8"};
{\ar@/_0.50pc/_{u} "4";"8"};
{\ar@/^0.50pc/^{0} "6";"8"};
{\ar@/_0.50pc/_{0} "6";"8"};
{\ar@{=>}^{\alpha} (10,-7)*{};(10,-9)*{}};
{\ar@{=>}^{\gamma} (10,2)*{};(10,-2)*{}};
{\ar@{=>}^{\beta} (21,0)*{};(19,0)*{}};
{\ar@{=>} (6,10)*{};(10,14)*{}};
{\ar@{..}@/_0.80pc/^{} (8,12)*{};(24,20)*{_{\text{equivalence}}}};
{\ar@{}|\circlearrowright"2";"10"};
\endxy
\]
Then we have $(i_1\circ t\circ\beta)\cdot(i_1\circ\gamma)=i_1\circ\gamma^{\prime}=(i_1\circ\gamma)\cdot(i_1\circ p_2\circ\alpha)=(i_1\circ\gamma)\cdot(\xi_2\circ0)\cdot(0\circ\alpha)\cdot(\xi_2^{-1}\circ0)=(i_1\circ\gamma)$, and thus, $(i_1\circ t)\circ\beta=\mathrm{id}_{i_1\circ t\circ u}$. Since $i_1\circ t\simeq f_1$ is cofaithful, we obtain $\beta=\mathrm{id}_u$. So, by $(\ref{proppullback})$, we have $\gamma=\gamma^{\prime}=\gamma\cdot(p_2\circ\alpha)$, and consequently $p_2\circ\alpha=\mathrm{id}_{p_2\circ0}$.
Since $p_2$ is cofaithful (because $i_2\circ p_2\simeq \mathrm{id}_{A_2}$ is cofaithful), we obtain $\alpha=\mathrm{id}_0$.
\end{proof}

\begin{proposition}\label{KerFactRelKer}
Consider diagram $(\ref{ElementaryDiagram})$ in $\mathbf{S}$. If we take $\mathrm{Ker}(f,\varphi)=[\mathrm{Ker}(f,\varphi),\ell,\varepsilon]$, then by the universality of $\mathrm{Ker}(f)=[\mathrm{Ker}(f),k(f),\varepsilon_f],$ $\ell$ factors uniquely up to an equivalence as
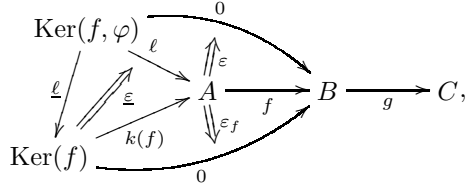
\[
\xy
(-29,-9)*+{\mathrm{Ker}(f)}="0";
(-24,8)*+{\mathrm{Ker}(f,\varphi)}="2";
(-8,0)*+{A}="4";
(8,0)*+{B}="6";
(24,0)*+{C}="8";
(26,-1)*+{,}="10";
{\ar_<<<<<<{k(f)} "0";"4"};
{\ar^{\ell} "2";"4"};
{\ar_{\underline{\ell}} "2";"0"};
{\ar_{f} "4";"6"};
{\ar_{g} "6";"8"};
{\ar@/^1.25pc/^{0} "2";"6"};
{\ar@/_1.25pc/_{0} "0";"6"};
{\ar@{=>}_<<<{\varepsilon} (-8,2)*{};(-7,7)*{}};
{\ar@{=>}^<<<<{\varepsilon_f} (-8,-2)*{};(-7,-7)*{}};
{\ar@{=>}_<<<<<<<<{\underline{\varepsilon}} (-25,-5)*{};(-18,3)*{}};
\endxy
\]
where $(\underline{\varepsilon}\circ f)\cdot\varepsilon=(\underline{\ell}\circ\varepsilon_f)\cdot(\underline{\ell})_I^{\sharp}$. Then, $\underline{\ell}$ becomes fully faithful.
\end{proposition}
\begin{proof}
Since $\underline{\ell}\circ k(f)$ is equivalent to a faithful 1-cell $\ell$, so $\underline{\ell}$ becomes faithful. For any $K\in\mathbf{S}^0$, $k\in\mathbf{S}^1(K,\mathrm{Ker}(f,\varphi))$ and $\sigma\in\mathbf{S}^2(k\circ\underline{\ell},0)$, if we put $\sigma^{\prime}:=(k\circ\underline{\varepsilon}^{-1})\cdot(\sigma\circ k(f))\cdot k(f)_I^{\flat}\in\mathbf{S}^2(k\circ\ell,0)$, then $\sigma^{\prime}$ becomes compatible with $\varepsilon$.
So, by Lemma \ref{KernelLemma}, there exists $\tau\in\mathbf{S}^2(k,0)$ such that $\sigma^{\prime}=(\tau\circ\ell)\cdot\ell_I^{\flat}$, i.e.
\begin{equation}
(k\circ\underline{\varepsilon}^{-1})\cdot(\sigma\circ k(f))\cdot (k(f))_I^{\flat}=(\tau\circ\ell)\cdot\ell_I^{\flat}.
\label{firstequality5}
\end{equation}
Now, since $(k\circ\underline{\varepsilon})\cdot(\tau\circ\ell)\cdot\ell_I^{\flat}=(\tau\circ\underline{\ell}\circ k(f))\cdot(\underline{\ell}\circ k(f))_I^{\flat}$ by Corollary \ref{ContraCor}, $(\ref{firstequality5})$ is equivalent to $(\sigma\circ k(f))\cdot(k(f))_I^{\flat}=(\tau\circ\underline{\ell}\circ k(f))\cdot(\underline{\ell}_I^{\flat}\circ k(f))\cdot(k(f))_I^{\flat}$.

Thus, we obtain $\sigma\circ k(f)=((\tau\circ\underline{\ell})\cdot\underline{\ell}_I^{\flat})\circ k(f)$. Since $k(f)$ is faithful, it follows that $\sigma=(\tau\circ\underline{\ell})\cdot\underline{\ell}_I^{\flat}$. Uniqueness of such $\tau$ follows from the faithfulness of\ $\underline{\ell}$. Thus $\underline{\ell}$ becomes fully faithful by Corollary \ref{CorOfLemLem}.
\end{proof}

\section{Existence of proper factorization systems}\label{Existence of PFS}

\begin{definition}
For any $A,B\in\mathbf{S}^0$ and $f\in\mathbf{S}^1(A,B)$, we define its image as $\mathrm{Ker}(\mathrm{cok}(f))$.
\end{definition}

\begin{remark}
By the universality of the kernel, there exist $i(f)\in\mathbf{S}^1(A,\mathrm{Im}(f))$ and $\iota\in\mathbf{S}^2(i(f)\circ k(c(f)),f)$
such that $(\iota\circ c(f))\cdot\pi_f=(i(f)\circ\varepsilon_{c(f)})\cdot i(f)_I^{\sharp}$.
Coimage of $f$ is defined dually, and we obtain a factorization through $\mathrm{Coim}(f)$.
\end{remark}

\begin{proposition}[1st factorization]\label{FirstFactorization}
For any $f\in\mathbf{S}^1(A,B)$, the factorization $\iota:i(f)\circ k(c(f))\Longrightarrow f$ through $\mathrm{Im}(f)$
\[
\xy
(-10,0)*+{A}="2";
(10,0)*+{B}="4";
(0,-12)*+{\mathrm{Im}(f)}="8";
{\ar^{f} "2";"4"};
{\ar_{k(c(f))} "8";"4"};
{\ar_{i(f)} "2";"8"};
{\ar@{=>}_{\iota} (0,-8)*{};(0,-2)*{}};
\endxy
\]
satisfies the following properties$:$\\
{\rm (A)} $i(f)$ is fully cofaithful and $k(c(f))$ is faithful.\\
{\rm (B)} For any factorization $\eta:i\circ m\Longrightarrow f$
where $m$ is faithful, following {\rm (b1)} and {\rm (b2)} hold$:$\\
{\rm (b1)} There exist $t\in\mathbf{S}^1(\mathrm{Im}(f),C)$, $\zeta_m\in\mathbf{S}^2(t\circ m,k(c(f)))$, $\zeta_i\in\mathbf{S}^2(i(f)\circ t,i)$
\[
\xy
(-12,0)*+{A}="2";
(12,0)*+{B}="4";
(0,12)*+{C}="6";
(0,-12)*+{\mathrm{Im}(f)}="8";
{\ar^{t} "8";"6"};
{\ar_{k(c(f))} "8";"4"};
{\ar_{i(f)} "2";"8"};
{\ar^{m} "6";"4"};
{\ar^{i} "2";"6"};
{\ar@{=>}^{\zeta_m} (3,4)*{};(5,-4)*{}};
{\ar@{=>}^{\zeta_i} (-3,-4)*{};(-5,4)*{}};
\endxy
\]
such that $(i(f)\circ\zeta_m)\cdot\iota=(\zeta_i\circ m)\cdot\eta$.
\\{\rm (b2)} If both $(t,\zeta_m,\zeta_i)$ and $(t^{\prime},\zeta_m^{\prime},\zeta_i^{\prime})$ satisfy {\rm (b1)}, then there is a unique 2-cell $\kappa\in\mathbf{S}^2(t,t^{\prime})$ such that $(i(f)\circ\kappa)\cdot\zeta_i^{\prime}=\zeta_i$ and $(\kappa\circ m)\cdot\zeta_m^{\prime}=\zeta_m$.
\end{proposition}

Dually, we obtain the following proposition for the coimage factorization.
\begin{proposition}[2nd factorization]\label{SecondFactorization}
For any $f\in\mathbf{S}^1(A,B)$, the factorization $\mu: c(k(f))\circ j(f)\Longrightarrow f$ through $\mathrm{Coim}(f)$
\[
\xy(-10,0)*+{A}="2";
(10,0)*+{B}="4";
(0,12)*+{\mathrm{Coim}(f)}="6";
{\ar_{f} "2";"4"};
{\ar^{j(f)} "6";"4"};
{\ar^{c(k(f))} "2";"6"};
{\ar@{=>}_{\mu} (0,8)*{};(0,2)*{}} ;
\endxy
\]
satisfies the following properties$:$\\
{\rm (A)} $j(f)$ is fully faithful and $c(k(f))$ is cofaithful.\\
{\rm (B)} For any factorization $\nu:e\circ j\Longrightarrow f$
where $e$ is cofaithful, following {\rm (b1)} and {\rm (b2)} $($the dual of the conditions in Proposition \ref{FirstFactorization}$)$ hold$:$\\
{\rm (b1)} There exists $s\in\mathbf{S}^1(C,\mathrm{Coim}(f))$, $\zeta_e\in\mathbf{S}^2(e\circ s,c(k(f)))$, and $\zeta_j\in\mathbf{S}^2(s\circ j(f),j)$
\[
\xy
(-12,0)*+{A}="2";
(12,0)*+{B}="4";
(0,12)*+{\mathrm{Coim}(f)}="6";
(0,-12)*+{C}="8";
{\ar^{s} "8";"6"};
{\ar_{j} "8";"4"};
{\ar_{e} "2";"8"};
{\ar^{j(f)} "6";"4"};
{\ar^{c(k(f))} "2";"6"};
{\ar@{=>}^{\zeta_j} (3,4)*{};(5,-4)*{}};
{\ar@{=>}^{\zeta_e} (-3,-4)*{};(-5,4)*{}};
\endxy
\]
such that $(e\circ\zeta_j)\cdot\nu=(\zeta_e\circ j(f))\cdot\mu$.\\
{\rm (b2)} If both $(s,\zeta_e,\zeta_j)$ and $(s^{\prime},\zeta_e^{\prime},\zeta_j^{\prime})$ satisfy {\rm (b1)}, then there is a unique 2-cell $\lambda\in\mathbf{S}^2(t,t^{\prime})$ such that \ $(\lambda\circ j(f))\cdot\zeta_j^{\prime}=\zeta_j$ and $(e\circ\lambda)\cdot\zeta_e^{\prime}=\zeta_e$.
\end{proposition}

In the rest of this section, we show Proposition \ref{FirstFactorization}.
\begin{lemma}
For any $f\in\mathbf{S}^1(A,B)$, $i(f)$ is cofaithful.
\end{lemma}
\begin{proof}
It suffices to show that for any $C\in\mathbf{S}^0$ and $\alpha\in\mathbf{S}^2(0_{\mathrm{Im}(f),C},0_{\mathrm{Im}(f),C})$
\[
\xy
(-18,0)*+{A}="0";
(0,0)*+{\mathrm{Im}(f)}="2";
(18,0)*+{C}="4";
{\ar^{i(f)} "0";"2"};
{\ar@/^1.25pc/^{0} "2";"4"};
{\ar@/_1.25pc/_{0} "2";"4"};
{\ar@{=>}^{\alpha} (9,2)*{};(9,-2)*{}};
\endxy
,
\]
we have $i(f)\circ\alpha=\mathrm{id}_{i(f)\circ0}\Longrightarrow\alpha=\mathrm{id}_0$. Take the pushout of $k(c(f))$ and $0_{\mathrm{Im}(f),C}$
\[
\xy
(-14,0)*+{\mathrm{Im}(f)}="2";
(14,0)*+{C\displaystyle\coprod_{\mathrm{Im}(f)}B}="4";
(0,8)*+{B}="6";
(0,-8)*+{C}="8";
{\ar_{i_C} "8";"4"};
{\ar_{0} "2";"8"};
{\ar^{i_B} "6";"4"};
{\ar^{k(c(f))} "2";"6"};
{\ar@{=>}^{\xi} (0,3)*{};(0,-3)*{}};
\endxy
\]
and put
\begin{eqnarray*}
\xi_1&:=&\xi\cdot(\xi_1\circ f_2)\cdot\eta=(g\circ\xi)\cdot(\xi_2\circ f_1)
(i_C)_I^{\flat}=(k(c(f))\circ i_B\overset{\xi}{\Longrightarrow}0\circ i_C\overset{(i_C)_I^{\flat}}{\Longrightarrow}0)\\
\xi_2&:=&\xi\cdot(\alpha\circ i_C)\cdot(i_C)_I^{\flat}=(k(c(f))\circ i_B\overset{\xi}{\Longrightarrow}0\circ i_C\overset{\alpha\circ i_C}{\Longrightarrow}0\circ i_C\overset{(i_C)_I^{\flat}}{\Longrightarrow}0).
\end{eqnarray*}
Then, since $i_C$ is faithful by (the dual of) Lemma \ref{FaithLemma}, we have
\[ \alpha=\mathrm{id}_0\Longleftrightarrow\alpha\circ i_C=\mathrm{id}_{0\circ i_C}\Longleftrightarrow\xi\cdot(\alpha\circ i_C)\cdot(i_C)_0^{\flat}=\xi\cdot \mathrm{id}_{0\circ i_C}\cdot(i_C)_I^{\flat}\Longleftrightarrow\xi_1=\xi_2. \]
So, it suffices to show $\xi_1=\xi_2$. For each $i=1,2$, since $\mathrm{Cok}(k(c(f))=[\mathrm{Cok}(f),c(f),\varepsilon_{c(f)}]$, there exist $e_i\in\mathbf{S}^1(\mathrm{Cok}(f),C\underset{\mathrm{Im}(f)}{\coprod}B)$ and $\varepsilon_i\in\mathbf{S}^2(c(f)\circ e_i,i_B)$ such that
\begin{equation}\label{lambda5-1}
(k(c(f))\circ\varepsilon_i)\cdot\xi_i=(\varepsilon_{c(f)}\circ e_i)\cdot(e_i)_{I}^{\flat}.
\end{equation}
\[
\xy
(-18,-14)*+{\mathrm{Im}(f)}="2";
(18,-14)*+{C\displaystyle\coprod_{\mathrm{Im}(f)}B}="4";
(2,0)*+{B}="6";
(20,12)*+{\mathrm{Cok}(f)}="8";
{\ar_{0} "2";"4"};
{\ar^{i_B} "6";"4"};
{\ar|*+{_{k(c(f)}} "2";"6"};
{\ar|*+{_{c(f)}} "6";"8"};
{\ar^{e_i} "8";"4"};
{\ar@/^2.20pc/^{0} "2";"8"};
{\ar@{=>}^{\xi_i} (0,-5)*{};(0,-11)*{}};
{\ar@{=>}^{\varepsilon_i} (14,4)*{};(11,-2)*{}};
{\ar@{=>}_{\varepsilon_{c(f)}} (-2,2)*{};(-4,4)*{}};
\endxy
\]
Since by assumption $i(f)\circ\alpha=\mathrm{id}_{i(f)\circ0},$ we have
\begin{eqnarray*}
i(f)\circ\xi_2&=&(i(f)\circ\xi)\cdot(i(f)\circ\alpha\circ i_C)\cdot(i(f)\circ(i_C)_I^{\flat})\\
&=&(i(f)\circ\xi)\cdot(\mathrm{id}_{i(f)\circ0\circ i_C})\cdot(i(f)\circ(i_C)_I^{\flat})=i(f)\circ\xi_1.
\end{eqnarray*}
So, if we put $\varpi:=(\iota^{-1}\circ i_B)\cdot(i(f)\circ\xi_i)\cdot(i(f))_I^{\sharp}\in\mathbf{S}^2(f\circ i_B,0)$, this doesn't depend on $i=1,2$. We can show easily
$(f\circ\varepsilon_i)\cdot\varpi=(\pi_f\circ e_i)\cdot(e_i)_I^{\flat}\ (i=1,2)$. Thus $(e_1,\varepsilon_1)$ and $(e_2,\varepsilon_2)$ are two factorizations of $i_B$ compatible with $\varpi$ and $\pi_f$.
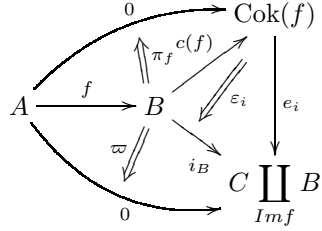
\[
\xy(-18,0)*+{A}="2";
(16,-12)*+{C\displaystyle\coprod_{Imf}B}="4";
(0,0)*+{B}="6";
(16,12)*+{\mathrm{Cok}(f)}="8";
{\ar@/_1.65pc/_{0} "2";"4"};
{\ar_{i_B} "6";"4"};
{\ar^{f} "2";"6"};
{\ar^{c(f)} "6";"8"};
{\ar^{e_i} "8";"4"};
{\ar@/^1.25pc/^{0} "2";"8"};
{\ar@{=>}_{\pi_f} (-1,3)*{};(-2,9)*{}};
{\ar@{=>}_{\varpi} (-1,-3)*{};(-4,-10)*{}};
{\ar@{=>}^{\varepsilon_i} (12,6)*{};(6,-2)*{}};
\endxy
\]
By the universality of $\mathrm{Cok}(f)$, there exists a 2-cell $\beta\in\mathbf{S}^2(e_1,e_2)$ such that $(c(f)\circ\beta)\cdot\varepsilon_2=\varepsilon_1$,
and thus we have $\varepsilon_1^{-1}=\varepsilon_2^{-1}\cdot(c(f)\circ\beta^{-1})$.
So, by $(\ref{lambda5-1})$, we have
\begin{eqnarray*}
\xi_1&=&(k(c(f))\circ\varepsilon_1^{-1})\cdot(\varepsilon_{c(f)}\circ e_1)\cdot(e_1)_I^{\flat}\\
&=&(k(c(f))\circ\varepsilon_2^{-1})\cdot(k(c(f))\circ c(f)\circ\beta^{-1})\cdot(\varepsilon_{c(f)}\circ e_1)\cdot(e_1)_I^{\flat}\\
&\underset{\ref{1-13}}{=}&(k(c(f))\circ\varepsilon_2^{-1})\cdot(\varepsilon_{c(f)}\circ e_2)\cdot(e_2)_I^{\flat}=\xi_2.
\end{eqnarray*}
\end{proof}

\begin{lemma}\label{Lemma4-2}
Let $f\in\mathbf{S}^{1}(A,B).$ Let $\iota:i(f)\circ k(c(f))\Longrightarrow f$ be the factorization of $f$ through $\mathrm{Im}(f)$ as before. If we are given a factorization $\eta:i\circ m\Longrightarrow f$ of $f$
where $i\in\mathbf{S}^1(A,C)$, $m\in\mathbf{S}^1(C,B)$ and $m$ is faithful, then there exist $t\in\mathbf{S}^1(\mathrm{Im}(f),C)$, $\zeta_i\in\mathbf{S}^2(i(f)\circ t,i)$ and $\zeta_m\in\mathbf{S}^2(t\circ m,k(c(f)))$ such that $(\zeta_i\circ m)\cdot\eta=(i(f)\circ\zeta_m)\cdot\iota.$
\end{lemma}
\begin{proof}
By the universality of $\mathrm{Cok}(f)$, for $\pi:=(\eta^{-1}\circ c(m))\cdot(i\circ\pi_m)\cdot i_I^{\sharp}\in\mathbf{S}^2(f\circ c(m),0)$, there exist $\overline{m}\in\mathbf{S}^1(\mathrm{Cok}(f),\mathrm{Cok}(m))$ and $\overline{\eta}\in\mathbf{S}^2(c(f)\circ\overline{m},c(m))$ such that
\begin{equation} (f\circ\overline{\eta})\cdot\pi=(\pi_f\circ\overline{m})\cdot(\overline{m})_I^{\flat}. \label{etaeta}
\end{equation}
\[
\xy
(-18,0)*+{A}="2";
(14,-12)*+{\mathrm{Cok}(m)}="4";
(0,0)*+{B}="6";
(16,12)*+{\mathrm{Cok}(f)}="8";
{\ar@/_1.65pc/_{0} "2";"4"};
{\ar^>>>>>>>{c(m)} "6";"4"};
{\ar^{f} "2";"6"};
{\ar^{c(f)} "6";"8"};
{\ar^{\overline{m} } "8";"4"};
{\ar@/^1.25pc/^{0} "2";"8"};
{\ar@{=>}^{\pi_f} (-1,3)*{};(-2,9)*{}};
{\ar@{=>}^{\pi} (-1,-3)*{};(-4,-10)*{}};
{\ar@{=>}^{\overline{\eta}} (12,6)*{};(9,0)*{}};
\endxy
\]
Since $m$ is faithful by assumption, it follows $\mathrm{Ker}(c(m))=[C,m,\pi_m]$. By the universality of $\mathrm{Ker}(c(m))$, for the 2-cell
\begin{equation}
\zeta:=(k(c(f))\circ\overline{\eta}^{-1})\cdot(\varepsilon_{c(f)}\circ\overline{m})\cdot(\overline{m})_I^{\flat}\in\mathbf{S}^2(k(c(f))\circ c(m),0),
\label{firstequality11}
\end{equation}
there exist $t\in\mathbf{S}^1(\mathrm{Im}(f),C)$ and $\zeta_m\in S^2(t\circ m,k(c(f)))$
such that $(\zeta_m\circ c(m))\cdot\zeta=(t\circ\pi_m)\cdot t_I^{\sharp}.$

If we put $\overline{\zeta}:=(i(f)\circ\zeta_{m})\cdot\iota$,
then the following claim holds:

\begin{claim}\label{Claim4-2}
Each of the two factorizations of $f$ through $\mathrm{Ker}(c(m))$
\[ \eta:i\circ m\Longrightarrow f\quad\text{and}\quad\overline{\zeta}:i(f)\circ t\circ m\Longrightarrow f \]
is compatible with $\pi_m$ and $\pi$.
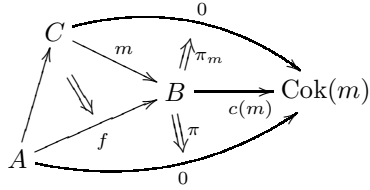
\[
\xy
(-29,-9)*+{A}="0";
(-24,8)*+{C}="2";
(-8,0)*+{B}="4";
(12,0)*+{\mathrm{Cok}(m)}="6";
{\ar_{f} "0";"4"};
{\ar^{m} "2";"4"};
{\ar "0";"2"};
{\ar_{c(m)} "4";"6"};
{\ar@/^1.25pc/^{0} "2";"6"};
{\ar@/_1.25pc/_{0} "0";"6"};
{\ar@{=>}_{\pi_m} (-7,3)*{};(-6,7)*{}};
{\ar@{=>}^{\pi} (-8,-3)*{};(-7,-8)*{}};
{\ar@{=>}_{} (-22,2)*{};(-19,-3)*{}};
\endxy
\]
\end{claim}
If the above claim is proven, then by the universality of $\mathrm{Ker}(c(m))=[C,m,\pi_m]$, there exists a unique 2-cell $\zeta_i\in\mathbf{S}^2(i(f)\circ t,i)$ such that $(\zeta_i\circ m)\cdot\eta=\overline{\zeta}$.
Thus we obtain $(t,\zeta_m,\zeta_i)$ which satisfies $(\zeta_i\circ m)\cdot\eta=\overline{\zeta}=(i(f)\circ\zeta_m)\cdot\iota$, and the lemma is proven. So, we show Claim \ref{Claim4-2}.

\noindent\underline{{\rm (a)} compatibility of $\eta$ with $\pi_m$, $\pi$}

This follows immediately from the definition of $\pi$.

\noindent\underline{{\rm (b)} compatibility of $\overline{\zeta}$ with $\pi_m$, $\pi$}

We have
\begin{eqnarray*}
i(f)\circ\zeta&\underset{\ref{firstequality11}}{=}&(\iota\circ c(m))\cdot(f\circ\overline{\eta}^{-1})\cdot(\iota^{-1}\circ c(f)\circ\overline{m})\\
& &\cdot\ (i(f)\circ\varepsilon_{c(f)}\circ\overline{m}))\cdot(i(f)\circ(\overline{m})_I^{\flat})\\
&\underset{\ref{etaeta}}{=}&(\iota\circ c(m))\cdot\pi\cdot i(f)_I^{\sharp-1}.
\end{eqnarray*}
From this, we obtain $(i(f)\circ t\circ\pi_m)\cdot(i(f)\circ t_I^{\sharp})=(\overline{\zeta}\circ c(m))\cdot\pi\cdot i(f)_I^{\sharp-1}$.
So we have
\[ (\overline{\zeta}\circ c(m))\cdot\pi=(i(f)\circ t\circ\pi_m)\cdot(i(f)\circ
t_I^{\sharp})\cdot i(f)_I^{\sharp-1}=(i(f)\circ t\circ\pi_m)\cdot(i(f)\circ t)_I^{\sharp}. \]
\end{proof}

\begin{lemma}\label{Lemma4-3}
Let $A,B,C\in\mathbf{S}^0$, $f,m,i\in\mathbf{S}^1$, $\eta\in\mathbf{S}^2$ be as in Lemma \ref{Lemma4-2}. If a triplet $(t^{\prime},\zeta_m^{\prime},\zeta_i^{\prime})$ $($where $t^{\prime}\in\mathbf{S}^1(\mathrm{Im}(f),C),\zeta_m^{\prime}\in\mathbf{S}^2(t^{\prime}\circ m,k(c(f))),\zeta_i^{\prime}\in\mathbf{S}^2(i(f)\circ t^{\prime},i)$ satisfies
\begin{equation}
(i(f)\circ\zeta_m^{\prime})\cdot\iota=(\zeta_i^{\prime}\circ m)\cdot\eta,
\label{lambda8-1}
\end{equation}
then $\zeta_m^{\prime}$ becomes compatible with $\zeta$ and $\pi_m$ $($in the notation of the proof of Lemma \ref{Lemma4-2}$)$, i.e. we have $(\zeta_m^{\prime}\circ c(m))\cdot\zeta=(t^{\prime}\circ\pi_m)\cdot(t^{\prime})_I^{\sharp}$.
\end{lemma}

\begin{remark}
Since $m$ is faithful, $\zeta_m^{\prime}$ which satisfies $(\ref{lambda8-1})$ is uniquely determined by $t^{\prime}$ and $\zeta_i^{\prime}$ if it exists.
\end{remark}
\begin{proof}[of Lemma \ref{Lemma4-3}]
Since we have
\begin{eqnarray*}
\lefteqn{i(f)\circ((\zeta_m^{\prime}\circ c(m)\cdot\zeta)}\hspace{1cm}\\
&\underset{\ref{lambda8-1},\ref{2-cell}}{=}&(\zeta_i^{\prime}\circ m\circ c(m))\cdot(\eta\circ c(m))\cdot(f\circ\overline{\eta}^{-1})\cdot(\iota^{-1}\circ c(f)\circ\overline{m})\\
& &\cdot\ (i(f)\circ\varepsilon_{c(f)}\circ\overline{m})\cdot(i(f)\circ(\overline{m})_I^{\flat})\\
&\underset{\ref{etaeta}}{=}&((i(f)\circ t^{\prime}\circ\pi_m)\cdot(i(f)\circ(t^{\prime})_I^{\sharp}),
\end{eqnarray*}
we obtain $(\zeta_m^{\prime}\circ c(m))\cdot\zeta=(t^{\prime}\circ\pi_m)\cdot(t^{\prime})_I^{\sharp}$ by the cofaithfulness of $i(f)$.
\end{proof}

\begin{corollary}
Let $A$, $B$, $C$, $f$, $m$, $i$, $\eta$ as in Proposition \ref{FirstFactorization}. If both $(t,\zeta_m,\zeta_i)$ and $(t^{\prime},\zeta_m^{\prime},\zeta_i^{\prime})$ satisfy {\rm (b1)}, then there exists a unique 2-cell $\kappa\in\mathbf{S}^2(t,t^{\prime})$ such that $(i(f)\circ\kappa)\cdot\zeta_i^{\prime}=\zeta_i$ and $(\kappa\circ m)\cdot\zeta_m^{\prime}=\zeta_m$.
\end{corollary}
\begin{proof}
By Lemma \ref{Lemma4-3}, there exists a 2-cell $\kappa\in\mathbf{S}^2(t,t^{\prime})$ such that $(\kappa\circ m)\cdot\zeta_m^{\prime}=\zeta_m$ by the universality of $\mathrm{Ker}(c(m))=[C,m,\pi_m]$. This $\kappa$ also satisfies $\zeta_i=(i(f)\circ\kappa)\cdot\zeta_i^{\prime}$, and unique by the cofaithfulness of $i(f)$.
\end{proof}
Considering the case of $C=\mathrm{Im}(f)$, we obtain the following corollary.
\begin{corollary}\label{Corollary4-2}
For any $t\in\mathbf{S}^1(\mathrm{Im}(f),\mathrm{Im}(f))$, $\zeta_m\in\mathbf{S}^2(t\circ k(c(f)),k(c(f)))$ and $\zeta_i\in\mathbf{S}^2(i(f)\circ t,i(f))$ satisfying $(\zeta_i\circ k(c(f)))\cdot\iota=(i(f)\circ\zeta_m)\cdot\iota,$ there exists a unique 2-cell $\kappa\in\mathbf{S}^2(t,\mathrm{id}_{\mathrm{Im}(f)})$ such that $i(f)\circ\kappa=\zeta_i$ and $\kappa\circ k(c(f))=\zeta_m$.
\end{corollary}

Now, we can prove Proposition \ref{FirstFactorization}.
\begin{proof}[of Proposition \ref{FirstFactorization}]

Since all the other is already shown, it suffices to show
the following:

\begin{claim}\label{Claim4-3}
For any $C\in\mathbf{S}^0$ and any $g,h\in\mathbf{S}^1(\mathrm{Im}(f),C)$,
\[ i(f)\circ-:\mathbf{S}^2(g,h)\longrightarrow\mathbf{S}^2(i(f)\circ g,i(f)\circ h) \]
is surjective.
\end{claim}
So, we show Claim \ref{Claim4-3}. If we take the difference kernel of $g$ and $h$;
\[ d_{(g,h)}:\mathrm{DK}(g,h)\longrightarrow \mathrm{Im}(f),\ \varphi_{(g,h)}:d_{(g,h)}\circ g\Longrightarrow d_{(g,h)}\circ h, \]
then by the universality of the difference kernel, for any $\beta\in\mathbf{S}^2(i(f)\circ g,i(f)\circ h)$ there exist $i\in\mathbf{S}^1(A,\mathrm{DK}(g,h))$ and $\lambda\in\mathbf{S}^2(i\circ d_{(g,h)},i(f))$
\[
\xy
(-28,0)*+{\mathrm{DK}(g,h)}="0";
(0,0)*+{\mathrm{Im}(f)}="2";
(24,0)*+{C}="4";
(-20,14)*+{A}="6";
{\ar_{d_{(g,h)}} "0";"2"};
{\ar@/^0.50pc/^{g} "2";"4"};
{\ar@/_0.50pc/_{h} "2";"4"};
{\ar^{i(f)} "6";"2"};
{\ar_{i} "6";"0"};
{\ar@{=>}_{\lambda} (-20,4)*{};(-14,6)*{}};
\endxy
\]
such that $(i\circ\varphi_{(g,h)})\cdot(\lambda\circ h)=(\lambda\circ g)\cdot\beta$.

If we put $m:=d_{(g,h)}\circ k(c(f))$, then $m$ becomes faithful since $d_{(g,h)}$ and $k(c(f))$ are faithful. Applying Lemma \ref{Lemma4-2} to the factorization $\eta:=(\lambda\circ k(c(f)))\cdot\iota:i\circ m\Longrightarrow f$, we obtain $t\in\mathbf{S}^1(\mathrm{Im}(f),\mathrm{DK}(g,h))$, $\zeta_m\in\mathbf{S}^2(t\circ m,k(c(f)))$ and $\zeta_i\in\mathbf{S}^2(i(f)\circ t,i)$ such that $(\zeta_i\circ m)\cdot\eta=(i(f)\circ\zeta_m)\cdot\iota$. Thus we have 
\[ (\zeta_i\circ d_{(g,h)}\circ k(c(f)))\cdot(\lambda\circ k(c(f)))\cdot\iota=(i(f)\circ\zeta_m)\cdot\iota. \]
So, if we put $\overline{\zeta}_i:=(\zeta_i\circ d_{(g,h)})\cdot\lambda\in\mathbf{S}^2(i(f)\circ t\circ d_{(g,h)},i(f))$, then we have
\[ (\overline{\zeta}_i\circ k(c(f)))\cdot\iota=(i(f)\circ\zeta_m)\cdot\iota. \]
By Corollary \ref{Corollary4-2}, there exists a 2-cell $\kappa\in\mathbf{S}^2(t\circ d_{(g,h)},\mathrm{id}_{\mathrm{Im}(f)})$ such that \ $\kappa\circ k(c(f))=\zeta_m$ and \ $i(f)\circ\kappa=\overline{\zeta}_{i}$. If we put $\alpha:=(\kappa^{-1}\circ g)\cdot(t\circ\varphi_{(g,h)})\cdot(\kappa\circ h)\in\mathbf{S}^2(g,h)$, we can show that $\alpha$ satisfies $i(f)\circ\alpha=\beta$.
Thus $i(f)\circ-:\mathbf{S}^2(g,h)\longrightarrow\mathbf{S}^2(i(f)\circ g,i(f)\circ h)$ \ is surjective.
\end{proof}

\begin{remark}
In condition {\rm (B)} of Proposition \ref{FirstFactorization}, if moreover $i$ is fully cofaithful, then $t$ becomes fully cofaithful since $i$ and $i(f)$ are fully cofaithful. On the other hand, $t$ is faithful since $k(c(f))$ is faithful. So, in this case $t$ becomes an equivalence by Proposition \ref{final1}.
\end{remark}
Together with Corollary \ref{Corollary4-2}, we can show easily the following corollary:

\begin{corollary}
For any $f\in\mathbf{S}^1(A,B)$, the following {\rm (b1)} and {\rm (b2)} hold:\\{\rm (b1)} If in the factorizations
\[
\xy
(-10,0)*+{A}="2";
(10,0)*+{B}="4";
(0,12)*+{C}="6";
{\ar_{f} "2";"4"};
{\ar^{m} "6";"4"};
{\ar^{i} "2";"6"};
{\ar@{=>}_{\eta} (0,8)*{};(0,2)*{}};
\endxy
\quad
\xy
(-10,0)*+{A}="2";
(10,0)*+{B}="4";
(12,-1)*+{,}="5";
(0,12)*+{C^{\prime}}="6";
{\ar_{f} "2";"4"};
{\ar^{m^{\prime}} "6";"4"};
{\ar^{i^{\prime}} "2";"6"};
{\ar@{=>}_{\eta'} (0,8)*{};(0,2)*{}};
\endxy
\]
$m,m^{\prime}$ are faithful and $i,i^{\prime}$ are fully cofaithful, then there exist $t\in\mathbf{S}^1(C,C^{\prime})$, $\zeta_m\in\mathbf{S}^2(t\circ m^{\prime},m),$ and $\zeta_i\in\mathbf{S}^2(i\circ t,i^{\prime})$ such that $(i\circ\zeta_m)\cdot\eta=(\zeta_i\circ m^{\prime})\cdot\eta^{\prime}$.\\
{\rm (b2)} If both $(t,\zeta_m,\zeta_i)$ and $(t^{\prime},\zeta_m^{\prime},\zeta_i^{\prime})$ satisfy {\rm (b1)}, then there is a unique 2-cell $\kappa\in\mathbf{S}^2(t,t^{\prime})$ such that $(i\circ\kappa)\cdot\zeta_i^{\prime}=\zeta_i$ and $(\kappa\circ m^{\prime})\cdot\zeta_m^{\prime}=\zeta_m$.
\end{corollary}

\begin{remark}
Proposition \ref{FirstFactorization} and Proposition \ref{SecondFactorization} implies respectively the existence of (2,1)-proper factorization system and (1,2)-proper factorization system in any relatively exact 2-category, in the sense of \cite{DV}.
\end{remark}

In the notation of this section, condition {\rm (B2)} and {\rm (B3)} in Definition \ref{RelExCat} can be written as follows:

\begin{corollary}\label{final2}
For any $f\in\mathbf{S}^1(A,B)$, we have$;$\\
{\rm (1)} $f$ is faithful iff $i(f):A\longrightarrow \mathrm{Im}(f)$ is an equivalence.\\
{\rm (2)} $f$ is cofaithful iff $j(f):\mathrm{Coim}(f)\longrightarrow B$ is an equivalence.\end{corollary}

\begin{proof}
Since {\rm (1)} is the dual of {\rm (2)}, we show only {\rm (2)}.

In the coimage factorization diagram
\[
\xy
(-10,0)*+{A}="2";
(10,0)*+{B}="4";
(0,12)*+{\mathrm{Coim}(f)}="6";
{\ar_{f} "2";"4"};
{\ar^{j(f)} "6";"4"};
{\ar^{c(k(f))} "2";"6"};
{\ar@{=>}_{\mu_f} (0,8)*{};(0,2)*{}};
\endxy
,
\]
since $c(k(f))$ is cofaithful and $j(f)$ is fully faithful, we have\\
$f$ is cofaithful$\Longleftrightarrow j(f)$ is cofaithful$\underset{\text{Prop. \ref{final1}}}{\Longleftrightarrow}j(f)$ is an equivalence.
\end{proof}

\section{Definition of relative 2-exactness}\label{Definition of Rel2Ex}

\paragraph{Diagram lemmas {\rm (1)}}\quad\\

\begin{definition}\label{DefComplex}
A complex $A_{\bullet}=(A_n,d_n^A,\delta_n^A)$ is a diagram
\[
\xy
(-48,0)*+{}="0";
(-36,0)*+{\cdots A_{n-2}}="1";
(-16,0)*+{A_{n-1}}="2";
(0,0)*+{A_n}="3";
(16,0)*+{A_{n+1}}="4";
(36,0)*+{A_{n+2}\cdots}="5";
(48,0)*+{}="6";
{\ar_{d^A_{n-2}} "1";"2"};
{\ar_{d^A_{n-1}} "2";"3"};
{\ar_{d^A_n} "3";"4"};
{\ar_{d^A_{n+1}} "4";"5"};
{\ar@/^1.65pc/^{0} "1";"3"};
{\ar@/_1.65pc/_{0} "2";"4"};
{\ar@/^1.65pc/^{0} "3";"5"};
{\ar@/_1.65pc/_{0} "0";"2"};
{\ar@/_1.65pc/_{0} "4";"6"};
{\ar@{=>}_<<<{\delta^A_{n-1}} (-18,2)*{};(-18,6)*{}};
{\ar@{=>}^<<<{\delta^A_n} (0,-2)*{};(0,-6)*{}};
{\ar@{=>}_<<<{\delta^A_{n+1}} (18,2)*{};(18,6)*{}};
{\ar@{=>}_<<<{} (-32,-2)*{};(-32,-6)*{}};
{\ar@{=>}_<<<{} (32,-2)*{};(32,-6)*{}};
\endxy
\]
where $A_n\in\mathbf{S}^0$, $d_n^A\in\mathbf{S}^1(A_n,A_{n+1})$, $\delta_n^A\in\mathbf{S}^2(d_{n-1}^A\circ d_n^A,0)$, and satisfies the following compatibility condition for each $n\in \mathbb{Z}:$
\[ (d_{n-1}^A\circ\delta_{n+1}^A)\cdot(d_{n-1}^A)_I^{\sharp}=(\delta_n^A\circ d_{n+1}^A)\circ(d_{n+1}^A)_I^{\flat} \]
\end{definition}

\begin{remark}\label{RemRemComplex}
We consider a bounded complex as a particular case of a complex, by adding zeroes.
\[
\xy
(-52,0)*+{\cdots}="-1";
(-48,0)*+{0}="0";
(-32,0)*+{0}="1";
(-16,0)*+{A_0}="2";
(0,0)*+{A_1}="3";
(16,0)*+{A_2}="4";
(36,0)*+{A_3\cdots}="5";
(48,0)*+{}="6";
{\ar_{0} "0";"1"};
{\ar_{0} "1";"2"};
{\ar_{d^A_0} "2";"3"};
{\ar_{d^A_1} "3";"4"};
{\ar_{d^A_2} "4";"5"};
{\ar@/^1.65pc/^{0} "1";"3"};
{\ar@/_1.65pc/_{0} "2";"4"};
{\ar@/^1.65pc/^{0} "3";"5"};
{\ar@/_1.65pc/_{0} "0";"2"};
{\ar@/_1.65pc/_{0} "4";"6"};
{\ar@{=>}_<<<{(d^A_0)^{\flat}_I} (-16,2)*{};(-16,6)*{}};
{\ar@{=>}^<<<{\delta^A_1} (0,-2)*{};(0,-6)*{}};
{\ar@{=>}_<<<{\delta^A_2} (18,2)*{};(18,6)*{}};
{\ar@{=>}_<<<{\mathrm{id}} (-32,-2)*{};(-32,-6)*{}};
{\ar@{=>}_<<<{} (32,-2)*{};(32,-6)*{}};
\endxy
\]
\end{remark}

\begin{definition}
For any complexes $A_{\bullet}=(A_n,d_n^A,\delta_n^A)$ and $B_{\bullet}=(B_n,d_n^B,\delta_n^B)$, a complex morphism $f_{\bullet}=(f_n,\lambda_n):A_{\bullet}\longrightarrow B_{\bullet}$ consists of $f_n\in\mathbf{S}^1(A_n,B_n)$ and $\lambda_n\in\mathbf{S}^2(d_n^A\circ f_{n+1},f_n\circ d_n^B)$ for each $n$, satisfying
\[ (\delta_n^A\circ f_{n+1})\cdot(f_{n+1})_I^{\flat}=(d_{n-1}^A\circ\lambda_n)\cdot(\lambda_{n-1}\circ d_n^B)\cdot(f_{n-1}\circ\delta_n^B)\cdot(f_{n-1})_I^{\sharp}.
\]
\[
\xy
(-32,4)*+{}="0";
(-36,6)*+{\cdots A_{n-2}}="1";
(-16,6)*+{A_{n-1}}="2";
(0,6)*+{A_n}="3";
(16,6)*+{A_{n+1}}="4";
(36,6)*+{A_{n+2}\cdots}="5";
(32,4)*+{}="6";
(-32,-4)*+{}="10";
(-36,-6)*+{\cdots B_{n-2}}="11";
(-16,-6)*+{B_{n-1}}="12";
(0,-6)*+{B_n}="13";
(16,-6)*+{B_{n+1}}="14";
(36,-6)*+{B_{n+2} \cdots}="15";
(32,-4)*+{}="16";
{\ar^{d^A_{n-2}} "1";"2"};
{\ar^{d^A_{n-1}} "2";"3"};
{\ar^{d^A_n} "3";"4"};
{\ar^{d^A_{n+1}} "4";"5"};
{\ar_{d^B_{n-2}} "11";"12"};
{\ar_{d^B_{n-1}} "12";"13"};
{\ar_{d^B_n} "13";"14"};
{\ar_{d^B_{n+1}} "14";"15"};
{\ar|{f_{n-2}} "0";"10"};
{\ar|{f_{n-1}} "2";"12"};
{\ar|{f_n} "3";"13"};
{\ar|{f_{n+1}} "4";"14"};
{\ar|{f_{n+2}} "6";"16"};
{\ar@{=>}|{\lambda_{n-2}} (-24,3)*{};(-24,-3)*{}};
{\ar@{=>}|{\lambda_{n-1}} (-8,3)*{};(-8,-3)*{}};
{\ar@{=>}|{\lambda_n} (8,3)*{};(8,-3)*{}};
{\ar@{=>}|{\lambda_{n+1}} (24,3)*{};(24,-3)*{}};
\endxy
\]
\end{definition}

\begin{proposition}\label{final5}
Consider the following diagram in $\mathbf{S}$.
\begin{eqnarray}\label{TACDiag}
\xy
(-16,6)*+{A_1}="0";
(-16,-6)*+{A_2}="2";
(0,6)*+{B_1}="4";
(0,-6)*+{B_2}="6";
{\ar^{f_1} "0";"4"};
{\ar_{a} "0";"2"};
{\ar_{f_2} "2";"6"};
{\ar^{b} "4";"6"};
{\ar@{=>}_{\lambda} (-8,2)*{};(-8,-2)*{}};
\endxy
\end{eqnarray}
If we take the cokernels of $f_1$ and $f_2$, then there exist $\overline{b}\in\mathbf{S}^1(\mathrm{Cok}(f_1),\mathrm{Cok}(f_2))$ and $\overline{\lambda}\in\mathbf{S}^2(c(f_1)\circ\overline{b},b\circ c(f_2))$ such that
\[ (\pi_{f_1}\circ\overline{b})\cdot(\overline{b})_I^{\flat}=(f_1\circ\overline{\lambda})\cdot(\lambda\circ c(f_2))\cdot(a\circ\pi_{f_2})\cdot a_I^{\sharp}. \]
\[
\xy
(-16,6)*+{A_1}="0";
(-16,-6)*+{A_2}="2";
(0,6)*+{B_1}="4";
(0,-6)*+{B_2}="6";
(18,6)*+{\mathrm{Cok}(f_1)}="8";
(18,-6)*+{\mathrm{Cok}(f_2)}="10";
{\ar^{f_1} "0";"4"};
{\ar_{a} "0";"2"};
{\ar_{f_2} "2";"6"};
{\ar|*+{_b} "4";"6"};
{\ar^{\overline{b}} "8";"10"};
{\ar^{c(f_1)} "4";"8"};
{\ar_{c(f_2)} "6";"10"};
{\ar@/^1.65pc/^{0} "0";"8"};
{\ar@/_1.65pc/_{0} "2";"10"};
{\ar@{=>}_{\lambda} (-8,2)*{};(-8,-2)*{}};
{\ar@{=>}^{\overline{\lambda}} (8,2)*{};(8,-2)*{}};
{\ar@{=>}_{\pi_{f_1}} (0,8)*{};(0,12)*{}};
{\ar@{=>}_{\pi_{f_2}} (0,-8)*{};(0,-12)*{}};
\endxy
\]
If $(\overline{b}^{\prime},\overline{\lambda}^{\prime})$ also satisfies this condition, there exists a unique 2-cell $\xi\in\mathbf{S}^2(\overline{b},\overline{b}^{\prime})$ such that $(c(f_1)\circ\xi)\cdot\overline{\lambda}^{\prime}=\overline{\lambda}$.
\end{proposition}
\begin{proof}
This follows immediately if we apply the universality of $\mathrm{Cok}(f_1)$ to $(\lambda\circ c(f_2))\cdot(a\circ\pi_{f_2})\cdot a_I^{\sharp}\in\mathbf{S}^2(f_1\circ b\circ c(f_2),0)$.
\end{proof}

\begin{proposition}\label{final6}
Consider the following diagrams in $\mathbf{S}$,
\[
\xy
(-16,12)*+{A_1}="0";
(-16,0)*+{A_2}="2";
(-16,-12)*+{A_3}="12";
(0,12)*+{B_1}="4";
(0,0)*+{B_2}="6";
(0,-12)*+{B_3}="16";
{\ar^{f_1} "0";"4"};
{\ar_{a_1} "0";"2"};
{\ar|*+{_{f_2}} "2";"6"};
{\ar^{b_1} "4";"6"};
{\ar_{f_3} "12";"16"};
{\ar_{a_2} "2";"12"};
{\ar^{b_2} "6";"16"};
{\ar@/_1.65pc/_{a} "0";"12"};
{\ar@/^1.65pc/^{b} "4";"16"};
{\ar@{=>}_{\lambda_1} (-8,8)*{};(-8,4)*{}};
{\ar@{=>}_{\lambda_2} (-8,-4)*{};(-8,-8)*{}};
{\ar@{=>}_{\alpha} (-19,0)*{};(-22,0)*{}};
{\ar@{=>}_{\beta} (3,0)*{};(6,0)*{}};
\endxy
\quad
\xy
(-16,6)*+{A_1}="0";
(-16,-6)*+{A_3}="2";
(0,6)*+{B_1}="4";
(0,-6)*+{B_3}="6";
{\ar^{f_1} "0";"4"};
{\ar_{a} "0";"2"};
{\ar_{f_3} "2";"6"};
{\ar^{b} "4";"6"};
{\ar@{=>}_{\lambda} (-8,2)*{};(-8,-2)*{}};
\endxy
\]
which satisfy $(f_1\circ\beta)\cdot\lambda=(\lambda_1\circ b_2)\cdot(a_1\circ\lambda_2)\cdot(\alpha\circ f_3)$. Applying Proposition \ref{final5}, we obtain diagrams
\[
\xy
(-16,6)*+{B_1}="0";
(-16,-6)*+{B_3}="2";
(2,6)*+{\mathrm{Cok}(f_1)}="4";
(2,-6)*+{\mathrm{Cok}(f_3)}="6";
{\ar^{c(f_1)} "0";"4"};
{\ar_{b} "0";"2"};
{\ar_{c(f_3)} "2";"6"};
{\ar^{\overline{b}} "4";"6"};
{\ar@{=>}^{\overline{\lambda}} (-8,2)*{};(-8,-2)*{}};
\endxy
\quad
\xy
(-16,6)*+{B_1}="0";
(-16,-6)*+{B_2}="2";
(2,6)*+{\mathrm{Cok}(f_1)}="4";
(2,-6)*+{\mathrm{Cok}(f_2)}="6";
{\ar^{c(f_1)} "0";"4"};
{\ar_{b_1} "0";"2"};
{\ar_{c(f_2)} "2";"6"};
{\ar^{\overline{b}_1} "4";"6"};
{\ar@{=>}^{\overline{\lambda}_1} (-8,2)*{};(-8,-2)*{}};
\endxy
\quad
\xy
(-16,6)*+{B_2}="0";
(-16,-6)*+{B_3}="2";
(2,6)*+{\mathrm{Cok}(f_2)}="4";
(2,-6)*+{\mathrm{Cok}(f_3)}="6";
{\ar^{c(f_2)} "0";"4"};
{\ar_{b_2} "0";"2"};
{\ar_{c(f_3)} "2";"6"};
{\ar^{\overline{b}_2} "4";"6"};
{\ar@{=>}^{\overline{\lambda}_2} (-8,2)*{};(-8,-2)*{}};
\endxy
\]
with
\begin{eqnarray}
(\pi_{f_1}\circ\overline{b})\cdot(\overline{b})_I^{\flat}&=&(f_1\circ\overline{\lambda})\cdot(\lambda\circ c(f_3))\cdot(a\circ\pi_{f_3})\cdot a_I^{\sharp}\label{final6-1}\\
(\pi_{f_1}\circ\overline{b}_1)\cdot(\overline{b}_1)_I^{\flat}&=&(f_1\circ\overline{\lambda}_1)\cdot(\lambda_1\circ c(f_2))\cdot(a_1\circ\pi_{f_2})\cdot(a_1)_I^{\sharp}\nonumber\\
(\pi_{f_2}\circ\overline{b}_2)\cdot(\overline{b}_2)_I^{\flat}&=&(f_2\circ\overline{\lambda}_2)\cdot(\lambda_2\circ c(f_3))\cdot(a_2\circ\pi_{f_3})\cdot(a_2)_I^{\sharp}.\nonumber
\end{eqnarray}
Then, there exists a unique 2-cell $\overline{\beta}\in\mathbf{S}^2(\overline{b}_1\circ\overline{b}_2,\overline{b})$ such that
\[ (c(f_1)\circ\overline{\beta})\cdot\overline{\lambda}=(\overline{\lambda}_1\circ\overline{b}_2)\cdot(b_1\circ\overline{\lambda}_2)\cdot(\beta\circ c(f_3)). \]\end{proposition}
\begin{proof}
By $(\ref{final6-1})$, $\overline{\lambda}$ is compatible with $\pi_{f_1}$ and $(\lambda\circ c(f_3))\cdot(a\circ\pi_{f_3})\cdot a_I^{\sharp}$.
\[
\xy
(-19,0)*+{A_1}="0";
(0,0)*+{B_1}="2";
(16,14)*+{\mathrm{Cok}(f_1)}="4";
(16,-14)*+{\mathrm{Cok}(f_3)}="6";
{\ar^{f_1 } "0";"2"};
{\ar|*+{_{c(f_1)}} "2";"4"};
{\ar|*+{_{b\circ c(f_3)}} "2";"6"};
{\ar^{\overline{b} } "4";"6"};
{\ar@/^1.20pc/^{0} "0";"4"};
{\ar@/_1.20pc/_{0} "0";"6"};
{\ar@{=>}_{\pi_{f_1}} (-2,4)*{};(-5,9)*{}};
{\ar@{=>}_{} (-2,-4)*{};(-5,-9)*{}};
{\ar@{=>}^{\overline{\lambda}} (12,4)*{};(8,-2)*{}};
{\ar@{.}@/^1.80pc/^{} (-4,-7)*{};(-8,-22)*{_{(\lambda\circ c(f_3))\cdot(a \circ\pi_{f_3})\cdot a^{\sharp} _I}}};
\endxy
\]
On the other hand, $\overline{\lambda}^{\prime}:=(\overline{\lambda}_1\circ\overline{b}_2)\cdot(b_1\circ\overline{\lambda}_2)\cdot(\beta\circ c(f_3))$ is also compatible with $\pi_{f_1}$ and $(\lambda\circ c(f_3))\cdot(a\circ\pi_{f_3})\cdot a_I^{\sharp}$. So, by the universality of the $\mathrm{Cok}(f_1)$, there exists a unique 2-cell $\overline{\beta}\in\mathbf{S}^2(\overline{b}_1\circ\overline{b}_2,\overline{b})$ such that $(c(f_1)\circ\overline{\beta})\cdot\overline{\lambda}=\overline{\lambda}^{\prime}$.
\end{proof}

\begin{corollary}\label{final7}
Let $(f_n,\lambda_n):(A_n,d_n^A,\delta_n^A)\longrightarrow (B_n,d_n^B,\delta_n^B)$ be a complex morphism.
Then, by taking the cokernels, we obtain a complex morphism $(c(f_n),\overline{\lambda}_n):(B_n,d_n^B,\delta_n^B)\longrightarrow(\mathrm{Cok}(f_n),\overline{d}_n^B,\overline{\delta}_n^B)$
which satisfies
\begin{equation}
(d_n^A\circ\pi_{f_{n+1}})\cdot(d_n^A)_I^{\sharp}=(\lambda_n\circ c(f_{n+1}))\cdot(f_n\circ\overline{\lambda}_n)\cdot(\pi_{f_n}\circ\overline{d}_n^B)\cdot(\overline{d}_n^B)_I^{\flat}
\label{kome5}
\end{equation}
for each $n$.
\end{corollary}
\begin{proof}
By Proposition \ref{final5}, we obtain $\overline{d}_n^B$ and $\overline{\lambda}_n$ which satisfy $(\ref{kome5})$. And by Proposition \ref{final6}, for each $n$, there exists a unique 2-cell $\overline{\delta}_n^B\in\mathbf{S}^2(\overline{d}_{n-1}^B\circ\overline{d}_n^B,0)$ such that $((\delta_n^B\circ c(f_{n+1}))\cdot c(f_{n+1})_I^{\flat}=(d_{n-1}^B\circ\overline{\lambda}_n)\cdot(\overline{\lambda}_{n-1}\circ\overline{d}_n^B)\cdot(c(f_{n-1})\circ\overline{\delta}_n^B)\cdot c(f_{n-1})_I^{\sharp}$. By the uniqueness of $\overline{\beta}$ in Proposition \ref{final6}, it is easy to see that
\[ (\overline{\delta}_n^B\circ\overline{d}_{n+1}^B)\cdot(\overline{d}_{n+1}^B)_I^{\flat}=(\overline{d}_{n-1}^B\circ\overline{\delta}_{n+1}^B)\cdot(\overline{d}_{n-1}^B)_I^{\sharp}. \]
These are saying that $(\mathrm{Cok}(f_n),\overline{d}_n^B,\overline{\delta}_n^B)$ is a complex and $(c(f_n),\overline{\lambda}_n)$ is a complex morphism.
\end{proof}

\begin{proposition}\label{final8}
Consider the following diagram in $\mathbf{S}$.
\[
\xy
(-16,6)*+{A_1}="0";
(-16,-6)*+{A_3}="2";
(0,6)*+{B_1}="4";
(0,-6)*+{B_3}="6";
{\ar^{f_1} "0";"4"};
{\ar_{a} "0";"2"};
{\ar_{f_3} "2";"6"};
{\ar^{b} "4";"6"};
{\ar@{=>}_{\lambda} (-8,2)*{};(-8,-2)*{}};
\endxy
\]
By taking the cokernels of $f_1$ and $f_2$, we obtain
\[
\xy
(-16,6)*+{A_1}="0";
(-16,-6)*+{A_2}="2";
(0,6)*+{B_1}="4";
(0,-6)*+{B_2}="6";
(20,6)*+{\mathrm{Cok}(f_1)}="8";
(20,-6)*+{\mathrm{Cok}(f_2)}="10";
(27,-7)*+{,}="12";
{\ar^{f_1} "0";"4"};
{\ar_{a} "0";"2"};
{\ar_{f_2} "2";"6"};
{\ar|*+{_b} "4";"6"};
{\ar^{\overline{b}} "8";"10"};
{\ar^{c(f_1)} "4";"8"};
{\ar_{c(f_2)} "6";"10"};
{\ar@{=>}_{\lambda} (-8,2)*{};(-8,-2)*{}} ;
{\ar@{=>}^{\lambda^{\prime}} (9,2)*{};(9,-2)*{}} ;
\endxy
\]
and from this diagram, by taking the cokernels of $a,b,\overline{b}$, we obtain
\[
\xy
(-20,8)*+{A_2}="0";
(0,8)*+{B_2}="2";
(-20,-8)*+{\mathrm{Cok}(a)}="10";
(0,-8)*+{\mathrm{Cok}(b)}="12";
(20,8)*+{\mathrm{Cok}(f_2)}="4";
(20,-8)*+{\mathrm{Cok}(\overline{b})}="14";
(27,-9)*+{.}="16";
{\ar_{f_2} "0";"2"};
{\ar_{c(f_2)} "2";"4"};
{\ar_{c(a)} "0";"10"};
{\ar|*+{_{c(b)}} "2";"12"};
{\ar^{c(\overline{b})} "4";"14"};
{\ar^{_{\overline{c(f_2)}}} "12";"14"};
{\ar_{\overline{f}_2} "10";"12"};
{\ar@{=>}_{\overline{\lambda}} (-10,2)*{};(-10,-2)*{}};
{\ar@{=>}^{\overline{\lambda}^{\prime}} (9,2)*{};(9,-2)*{}};
{\ar@/^1.65pc/^{0} "0";"4"};
{\ar@/_1.65pc/_{0} "10";"14"};
{\ar@{=>}_{\pi_{f_2}} (0,10)*{};(0,14)*{}};
{\ar@{=>}^{ \overline{\pi}_{f_2}} (0,-10)*{};(0,-14)*{}};
\endxy
\]
Then we have $\mathrm{Cok}(\overline{f}_{2})=[\mathrm{Cok}(\overline{b}),\overline{c(f_{2})},\overline{\pi}_{f_{2}}]$. We abbreviate this to $\mathrm{Cok}(\overline{f}_{2})=\mathrm{Cok}(\overline{b})$.
\end{proposition}
\begin{proof}
Left to the reader.
\end{proof}

\begin{proposition}\label{final9}
In the following diagram, assume $f_{\bullet}:A_{\bullet}\longrightarrow B_{\bullet}$ is a complex morphism.
\begin{equation}
\xy
(-16,6)*+{A_1}="0";
(-16,-6)*+{B_1}="2";
(0,6)*+{A_2}="4";
(0,-6)*+{B_2}="6";
(16,6)*+{A_3}="8";
(16,-6)*+{B_3}="10";
{\ar^{d_1^A} "0";"4"};
{\ar_{f_1} "0";"2"};
{\ar_{d_1^B} "2";"6"};
{\ar|*+{_{f_2}} "4";"6"};
{\ar^{f_3} "8";"10"};
{\ar^{d_2^A} "4";"8"};
{\ar_{d_2^B} "6";"10"};
{\ar@/^1.65pc/^{0} "0";"8"};
{\ar@/_1.65pc/_{0} "2";"10"};
{\ar@{=>}_{\lambda_1} (-8,2)*{};(-8,-2)*{}};
{\ar@{=>}_{\lambda_2} (8,2)*{};(8,-2)*{}};
{\ar@{=>}_{\delta_2^A} (0,9)*{};(0,12)*{}};
{\ar@{=>}_{\delta_2^B} (0,-9)*{};(0,-12)*{}};
\endxy
\label{finaldia8-5}
\end{equation}

If we take the cokernels of $d_1^A$ and $d_1^B$,
\[
\xy
(-16,6)*+{A_1}="0";
(-16,-6)*+{B_1}="2";
(0,6)*+{A_2}="4";
(0,-6)*+{B_2}="6";
(20,6)*+{\mathrm{Cok}(d_1^A)}="8";
(20,-6)*+{\mathrm{Cok}(d_1^B)}="10";
{\ar^{d_1^A} "0";"4"};
{\ar_{f_1} "0";"2"};
{\ar_{d_1^B} "2";"6"};
{\ar|*+{_{f_2}} "4";"6"};
{\ar^{\overline{f}_2} "8";"10"};
{\ar^{c(d_1^A)} "4";"8"};
{\ar_{c(d_1^B)} "6";"10"};
{\ar@{=>}_{\lambda_1} (-8,2)*{};(-8,-2)*{}};
{\ar@{=>}_{\overline{\lambda}_1} (10,2)*{};(10,-2)*{}};
\endxy
\]
then by the universality of cokernel, we obtain $\overline{d}_2^A\in\mathbf{S}^1(\mathrm{Cok}(d_1^A),A_3)$ and $\overline{\delta}_2^A\in\mathbf{S}^2(c(d_1^A)\circ\overline{d}_2^A,d_2^A)$ such that $(d_1^A\circ\overline{\delta}_2^A)\cdot\delta_2^A=(\pi_{d_1^A}\circ\overline{d}_2^A)\cdot(\overline{d}_2^A)_I^{\flat}$. Similarly, we obtain $\overline{d}_2^B\in\mathbf{S}^1(\mathrm{Cok}(d_1^B),B_3)$, $\overline{\delta}_2^B\in\mathbf{S}^2(c(d_1^B)\circ\overline{d}_2^B,d_2^B)$ with $(d_1^B\circ\overline{\delta}_2^B)\cdot\delta_2^B=(\pi_{d_1^B}\circ\overline{d}_2^B)\cdot(\overline{d}_2^B)_I^{\flat}$. Then, there exists a unique 2-cell $\overline{\lambda}_2\in\mathbf{S}^2(\overline{d}_2^A\circ f_3,\overline{f}_2\circ\overline{d}_2^B)$
such that $(c(d_1^A)\circ\overline{\lambda}_2)\cdot(\overline{\lambda}_1\circ\overline{d}_2^B)\cdot(f_2\circ\overline{\delta}_2^B)=(\overline{\delta}_2^A\circ f_3)\cdot\lambda_2$.
\end{proposition}
\begin{proof}
If we put $\delta:=(d_1^A\circ\lambda_2^{-1})\cdot(\delta_2^A\circ f_3)\cdot(f_3)_I^{\flat}$, then both the factorizations
\begin{eqnarray*}
(\overline{\delta}_2^A\circ f_3)\cdot\lambda_2&:&c(d_1^A)\circ(\overline{d}_2^A\circ f_3)\Longrightarrow f_2\circ d_2^B\\
(\overline{\lambda}_1\circ\overline{d}_2^B)\cdot(f_2\circ\overline{\delta}_2^B)&:&c(d_1^A)\circ(\overline{f}_2\circ\overline{d}_2^B)\Longrightarrow f_2\circ d_2^B
\end{eqnarray*}
are compatible with $\pi_{d_1^A}$ and $\delta$. So the proposition follows from the universality of $\mathrm{Cok}(d_2^A)$.
\end{proof}

\begin{proposition}\label{final10}
In diagram $(\ref{TACDiag})$, if we take the coimage factorization $\mu_a:c(k(a))\circ j(a)\Longrightarrow a$ and $\mu_b:c(k(b))\circ j(b)\Longrightarrow b$, then there exist $f\in\mathbf{S}^1(\mathrm{Coim}(a),\mathrm{Coim}(b))$, $\lambda_1\in\mathbf{S}^2(f_1\circ c(k(b)),c(k(a))\circ f)$ and $\lambda_2\in\mathbf{S}^2(f\circ j(b),j(a)\circ f_2)$ such that $(f_1\circ\mu_b)\cdot\lambda=(\lambda_1\circ j(b))\cdot(c(k(a))\circ\lambda_2)\cdot(\mu_a\circ f_2)$.
\begin{equation}
\xy
(-12,14)*+{A_1}="0";
(-12,0)*+{\mathrm{Coim}(a)}="2";
(-12,-14)*+{A_2}="4";
(12,14)*+{B_1}="10";
(12,0)*+{\mathrm{Coim}(b)}="12";
(12,-14)*+{B_2}="14";
{\ar^{f_1} "0";"10"};
{\ar|*+{_f} "2";"12"};
{\ar_{f_2} "4";"14"};
{\ar|*+{_{c(k(a))}} "0";"2"};
{\ar|*+{_{j(a)}} "2";"4"};
{\ar|*+{_{c(k(b))}} "10";"12"};
{\ar|*+{_{j(b)}} "12";"14"};
{\ar@/_2.80pc/_{a} "0";"4"};
{\ar@/^2.80pc/^{b} "10";"14"};
{\ar@{=>}^{\lambda_1} (0,9)*{};(0,5)*{}};
{\ar@{=>}^{\lambda_2} (0,-5)*{};(0,-9)*{}};
{\ar@{=>}_{\mu_a} (-19,0)*{};(-22,0)*{}};
{\ar@{=>}_{\mu_b} (19,0)*{};(22,0)*{}};
\endxy
\label{coimfactorab}
\end{equation}
Moreover, for any other $(f^{\prime},\lambda_1^{\prime},\lambda_2^{\prime})$ with this property, there exists a unique 2-cell $\xi\in\mathbf{S}^2(f,f^{\prime})$ such that \ $\lambda_1\cdot(c(k(a))\circ\xi)=\lambda_1^{\prime}$ and $(\xi\circ j(b))\cdot\lambda_2^{\prime}=\lambda_2$.
\end{proposition}
\begin{proof}
Since the coimage factorization is unique up to an equivalence and is obtained
by the factorization which fills in the following diagram, we may assume $\mathrm{Ker}(a)=[\mathrm{Ker}(a),k(a),\varepsilon_a]$, $\mathrm{Cok}(k(a))=[\mathrm{Coim}(a),c(k(a)),\pi_{k(a)}]$, and $(k(a)\circ\mu_a)\cdot\varepsilon_a=(\pi_{k(a)}\circ j(a))\cdot j(a)_I^{\flat}$.
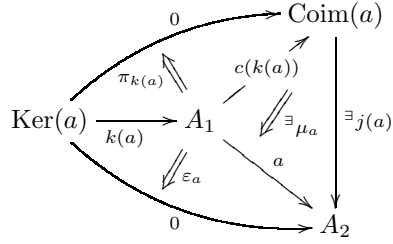
\[
\xy
(-20,0)*+{\mathrm{Ker}(a)}="0";
(0,0)*+{A_1}="2";
(18,14)*+{\mathrm{Coim}(a)}="4";
(18,-14)*+{A_2}="6";
{\ar_{k(a)} "0";"2"};
{\ar|*+{_{c(k(a))}} "2";"4"};
{\ar^{a} "2";"6"};
{\ar@/^1.20pc/^{0} "0";"4"};
{\ar@/_1.20pc/_{0} "0";"6"};
{\ar@{=>}^{\pi_{k(a) } } (-2,4)*{};(-5,9)*{}} ;
{\ar@{=>}^{\varepsilon_a } (-2,-4)*{};(-5,-9)*{}} ;
{\ar@{=>}^{^{\exists} \mu_a} (12,4)*{};(8,-2)*{}} ;
{\ar^{^{\exists} j(a) } "4";"6"};
\endxy
\]
Similarly, we may assume 
\begin{eqnarray*}
\mathrm{Ker}(b)&=&[\mathrm{Ker}(b),k(b),\varepsilon_b],\\
\mathrm{Cok}(k(b))&=&[\mathrm{Coim}(b),c(k(b)),\pi_{k(b)}]
\end{eqnarray*}
and $(k(b)\circ\mu_b)\cdot\varepsilon_b=(\pi_{k(b)}\circ j(b))\cdot j(b)_I^{\flat}$.
By (the dual of) Proposition \ref{final5}, there are $\underline{f}_1\in\mathbf{S}^1(\mathrm{Ker}(a),\mathrm{Ker}(b))$ and $\underline{\lambda}\in\mathbf{S}^2(\underline{f}_1\circ k(b),k(a)\circ f_1)$ such that $(\underline{\lambda}\circ b)\cdot(k(a)\circ\lambda)\cdot(\varepsilon_a\circ f_2)\cdot(f_2)_I^{\flat}=(\underline{f}_1\circ\varepsilon_b)\cdot(\underline{f}_1)_I^{\sharp}$. Applying Proposition \ref{final9}, we can show the existence of $(f,\lambda_1,\lambda_2)$. To show the uniqueness (up to an equivalence), let $(f^{\prime},\lambda_1^{\prime},\lambda_2^{\prime})$ satisfy
\[ (f_1\circ\mu_b)\cdot\lambda=(\lambda_1^{\prime}\circ j(b))\cdot (c(k(a))\circ\lambda_2^{\prime})\cdot(\mu_a\circ f_2). \]
From this, we can obtain
\[ (\underline{f}_1\circ\pi_{k(b)})\cdot(\underline{f}_1)_I^{\sharp}=(\underline{\lambda}\circ c(k(b)))\cdot(k(a)\circ\lambda_1^{\prime})\cdot(\pi_{k(a)}\circ f^{\prime})\cdot f_I^{\prime\flat}. \]
And the uniqueness follows from the uniqueness of 2-cells in Proposition \ref{final5} and Proposition \ref{final9}.
\end{proof}

\begin{proposition}\label{FINell}
Let $f_{\bullet}:A_{\bullet}\longrightarrow B_{\bullet}$ be a complex morphism as in diagram $(\ref{finaldia8-5})$. If we take the cokernels of $f_1,f_2,f_3$ and relative cokernels of the complex $A_{\bullet}$ and $B_{\bullet}$ as in the following diagram, then we have $\mathrm{Cok}(\overline{f}_3)=\mathrm{Cok}(\overline{d}_2^B,\overline{\delta}_2^B)$.
\[
\xy
(-20,4)*+{A_1}="11";
(0,4)*+{A_2}="12";
(20,4)*+{A_3}="13";
(48,4)*+{\mathrm{Cok}(d^A_2 ,\delta^A_2)}="14";
(-20,-10)*+{B_1}="21";
(0,-10)*+{B_2}="22";
(20,-10)*+{B_3}="23";
(48,-10)*+{\mathrm{Cok}(d^B_2 ,\delta^B_2)}="24";
(-20,-24)*+{\mathrm{Cok}(f_1)}="31";
(0,-24)*+{\mathrm{Cok}(f_2)}="32";
(20,-24)*+{\mathrm{Cok}(f_3)}="33";
{\ar^{d_1^A} "11";"12"};
{\ar^{d_2^A} "12";"13"};
{\ar^{c(d^A_2, \delta^A_2)} "13";"14"};
{\ar_{f_1} "11";"21"};
{\ar|*+{_{f_2}} "12";"22"};
{\ar|*+{_{f_3}} "13";"23"};
{\ar^{\overline{f}_3} "14";"24"};
{\ar_{d_1^B} "21";"22"};
{\ar_{d_2^B} "22";"23"};
{\ar_{c(d^B_2, \delta^B_2)} "23";"24"};
{\ar_{c(f_1)} "21";"31"};
{\ar|*+{_{c(f_2)}} "22";"32"};
{\ar^{c(f_3)} "23";"33"};
{\ar_{\overline{d}_1^B} "31";"32"};
{\ar_{\overline{d}_2^B} "32";"33"};
{\ar@/^1.65pc/^{0} "11";"13"};
{\ar@/_1.65pc/_{0} "31";"33"};
{\ar@{=>}_{\lambda_1} (-10,1)*{};(-10,-5)*{}};
{\ar@{=>}_{\lambda_2 } (10,1)*{};(10,-5)*{}};
{\ar@{=>}_{^{\exists}} (32,1)*{};(32,-5)*{}};
{\ar@{=>}_{\delta_2^A} (0,6)*{};(0,10)*{}};
{\ar@{=>}_{\overline{\delta} _2^B} (0,-26)*{};(0,-30)*{}};
{\ar@{=>}_{\overline{\lambda}_1} (-10,-15)*{};(-10,-21)*{}};
{\ar@{=>}_{\overline{\lambda}_2} (10,-15)*{};(10,-21)*{}};
\endxy
\]
\end{proposition}
\begin{proof}
Immediately follows from Proposition \ref{final8}, Proposition \ref{final9} and (the dual of) Proposition \ref{KerAndRelKer}.
\end{proof}
\begin{proposition}\label{final12}
In diagram $(\ref{TACDiag})$, if $a$ is fully cofaithful, then the following diagram obtained in Proposition \ref{final5} is a pushout diagram.
\[
\xy
(0,6)*+{B_1}="4";
(0,-6)*+{B_2}="6";
(20,6)*+{\mathrm{Cok}(f_1)}="8";
(20,-6)*+{\mathrm{Cok}(f_2)}="10";
{\ar_{b} "4";"6"};
{\ar^{\overline{b}} "8";"10"};
{\ar^{c(f_1)} "4";"8"};
{\ar_{c(f_2)} "6";"10"};
{\ar@{=>}^{\overline{\lambda}} (10,2)*{};(10,-2)*{}};
\endxy
\]
\end{proposition}
\begin{proof}
Left to the reader.
\end{proof}

Concerning Proposition \ref{FaithProp}, we have the following proposition.
\begin{proposition} \label{final13}
Let
\[
\xy
(-11,6)*+{A_1\times_BA_2}="0";
(-11,-6)*+{A_1}="2";
(11,6)*+{A_2}="4";
(11,-6)*+{B}="6";
{\ar^{f_1^{\prime}} "0";"4"};
{\ar_{f_2^{\prime}} "0";"2"};
{\ar_{f_1} "2";"6"};
{\ar^{f_2} "4";"6"};
{\ar@{=>}_{\xi} (0,2)*{};(0,-2)*{}};
\endxy
\]
be a pullback diagram in $\mathbf{S}$. If $f_1$ is fully cofaithful, then $f_1^{\prime}$ is fully cofaithful.
\end{proposition}
\begin{proof}
Since $f_1$ is cofaithful, in the notation of the proof of Proposition \ref{FaithProp}, $\mathrm{Cok}(i_1)=[A_2,p_2,\xi_2]$ and $\mathrm{Cok}(d)=[B,t,\varepsilon_t]$. Applying Proposition \ref{final8} to the diagram
\[
\xy
(-20,6)*+{0}="0";
(-20,-6)*+{A_1\times_BA_2}="2";
(0,6)*+{A_1}="4";
(0,-6)*+{A_1\times A_2}="6";
(7,-7)*+{,}="7";
{\ar^{0} "0";"4"};
{\ar_{0} "0";"2"};
{\ar_{d} "2";"6"};
{\ar^{i_1} "4";"6"};
{\ar@{=>} (-10,2)*{};(-10,-2)*{}};
\endxy
\]
we obtain
\[ \mathrm{Cok}(f_1)=0\Longleftrightarrow \mathrm{Cok}(f_1^{\prime})=0. \]
\[
\xy
(-24,12)*+{0}="0";
(-24,0)*+{A_1\times_BA_2}="2";
(-24,-12)*+{A_1\times_BA_2}="4";
(0,12)*+{A_1}="10";
(0,0)*+{A_1\times A_2}="12";
(0,-12)*+{A_2}="14";
(22,12)*+{A_1}="20";
(22,0)*+{B}="22";
(22,-12)*+{0}="24";
{\ar_{0} "0";"2"};
{\ar^{0} "0";"10"};
{\ar_{0} "14";"24"};
{\ar^{0} "22";"24"};
{\ar@{=}_{\mathrm{id}} "2";"4"};
{\ar@{=}^{\mathrm{id}} "10";"20"};
{\ar|*+{_d} "2";"12"};
{\ar|*+{_t} "12";"22"};
{\ar|*+{_{i_1}} "10";"12"};
{\ar|*+{_{p_2}} "12";"14"};
{\ar^{f_1} "20";"22"};
{\ar_{f_1^{\prime}} "4";"14"};
{\ar@{}|\circlearrowright"2";"14"};
{\ar@{=>} (-12,8)*{};(-12,4)*{}};
{\ar@{=>} (11,8)*{};(11,4)*{}};
{\ar@{=>} (11,-4)*{};(11,-8)*{}};
\endxy
\]
\end{proof}

\begin{proposition}\label{final15}
In diagram $(\ref{TACDiag})$, assume $a$ is cofaithful. By Proposition \ref{final10}, we obtain a coimage factorization diagram as $(\ref{coimfactorab})$. If we take the cokernel of this diagram as
\[
\xy
(-14,16)*+{B_1}="0";
(-14,0)*+{\mathrm{Coim}(b)}="2";
(-14,-16)*+{B_2}="4";
(14,16)*+{\mathrm{Cok}(f_1)}="10";
(14,0)*+{\mathrm{Cok}(f)}="12";
(14,-16)*+{\mathrm{Cok}(f_2)}="14";
(20,-17)*+{,}="16";
{\ar^{c(f_1)} "0";"10"};
{\ar|*+{_{c(f)}} "2";"12"};
{\ar_{c(f_2)} "4";"14"};
{\ar|*+{_{c(k(b))}} "0";"2"};
{\ar|*+{_{j(b)}} "2";"4"};
{\ar|*+{_{\overline{c(k(b))}}} "10";"12"};
{\ar|*+{_{\overline{j(b)}}} "12";"14"};
{\ar@/_2.80pc/_{b} "0";"4"};
{\ar@/^2.80pc/^{\overline{b}} "10";"14"};
{\ar@{=>}^{\overline{\lambda_1}} (0,10)*{};(0,6)*{}};
{\ar@{=>}^{\overline{\lambda_2}} (0,-6)*{};(0,-10)*{}};
{\ar@{=>}_{\mu_b} (-21,0)*{};(-24,0)*{}};
{\ar@{=>}^{\overline{\mu_b}} (21,0)*{};(24,0)*{}};
\endxy
\]
then the factorization
\[
\xy
(-14,-6)*+{\mathrm{Cok}(f_1)}="2";
(14,-6)*+{\mathrm{Cok}(f_2)}="4";
(0,8)*+{\mathrm{Cok}(f)}="6";
{\ar_{\overline{b}} "2";"4"};
{\ar^{\overline{j(b)}} "6";"4"};
{\ar^{\overline{c(k(b))}} "2";"6"};
{\ar@{=>}_{\overline{\mu_b}} (0,2)*{};(0,-3)*{}};
\endxy
\]
becomes again a coimage factorization.
\end{proposition}

\begin{proof}
It suffices to show that $\overline{c(k(b))}$ is cofaithful and $\overline
{j(b)}$ is fully faithful. Since $c(k(b))$ and $c(f)$ are cofaithful, it
follows that $\overline{c(k(b))}$ is cofaithful. Since $j(a)$ is an
equivalence,
\[
\xy
(-16,6)*+{\mathrm{Coim}(b)}="2";
(-16,-6)*+{B_2}="4";
(8,6)*+{\mathrm{Cok}(f)}="12";
(8,-6)*+{\mathrm{Cok}(f_2)}="14";
{\ar^{c(f)} "2";"12"};
{\ar_{c(f_2)} "4";"14"};
{\ar_{j(b)} "2";"4"};
{\ar^{\overline{j(b)}} "12";"14"};
{\ar@{=>}^{\overline{\lambda_2}} (-5,2)*{};(-5,-2)*{}};
\endxy
\]
is a pushout diagram by Proposition \ref{final12}. By (the dual of) Proposition \ref{final13}, $\overline{j(b)}$ becomes fully faithful.
\end{proof}

\paragraph{Definition of the relative 2-exactness}\quad\\

\begin{lemma}\label{final16}
Consider the following diagram in $\mathbf{S}$.
\begin{eqnarray}
\xy
(-8,0)*+{A}="4";
(8,0)*+{B}="6";
(24,0)*+{C}="8";
{\ar^{f} "4";"6"};
{\ar^{g} "6";"8"};
{\ar@/_1.65pc/_{0} "4";"8"};
{\ar@{=>}_{\varphi} (8,-2)*{};(8,-6)*{}};
\endxy
\label{CRUDiag}
\end{eqnarray}
If we factor it as
\begin{eqnarray}
\xy
(0,14)*+{\mathrm{Ker}(g)}="0";
(24,14)*+{\mathrm{Cok}(f)}="2";
(-12,0)*+{A}="4";
(12,0)*+{B}="6";
(36,0)*+{C}="8";
{\ar_{f} "4";"6"};
{\ar_{g} "6";"8"};
{\ar|*+{_{k(g)}} "0";"6"};
{\ar|*+{_{c(f)}} "6";"2"};
{\ar^{\underline{f}} "4";"0"};
{\ar^{\overline{g}} "2";"8"};
{\ar@{=>}_{\underline{\varphi}} (0,8)*{};(0,3)*{}};
{\ar@{=>}^{\overline{\varphi}} (24,8)*{};(24,3)*{}};
{\ar@/_1.80pc/_{0} "4";"8"};
{\ar@{=>}_{\varphi} (12,-2)*{};(12,-6)*{}};
\endxy
\label{CRUDiag2}
\end{eqnarray}
with
\begin{eqnarray*}
(\underline{\varphi}\circ g)\cdot\varphi&=&(\underline{f}\circ\varepsilon_g)\cdot(\underline{f})_I^{\sharp}\\
(f\circ\overline{\varphi})\cdot\varphi&=&(\pi_f\circ\overline{g})\cdot(\overline{g})_I^{\flat},
\end{eqnarray*}
then $\mathrm{Cok}(\underline{f})=0$ if and only if $\mathrm{Ker}(\overline{g})=0$.
\end{lemma}
\begin{proof}
We show only $\mathrm{Cok}(\underline{f})=0$ $\Rightarrow$ $\mathrm{Ker}(\overline{g})=0$, since the other implication can be shown dually. If $\mathrm{Cok}(\underline{f})=0$, i.e. if $\underline{f}$ is fully cofaithful, then we have
\[ \mathrm{Cok}(f)=\mathrm{Cok}(\underline{f}\circ k(g))=\mathrm{Cok}(k(g))=\mathrm{Coim}(g). \]
Thus the following diagram is a coimage factorization, and $\overline{g}$ becomes fully faithful.
\[
\xy
(24,14)*+{\mathrm{Cok}(f)}="2";
(12,0)*+{B}="6";
(36,0)*+{C}="8";
{\ar_{g} "6";"8"};
{\ar^{c(f)} "6";"2"};
{\ar^{\overline{g}} "2";"8"};
{\ar@{=>}^{\overline{\varphi}} (24,8)*{};(24,3)*{}} ;
\endxy
\]
\end{proof}

\begin{definition}\label{final17}
Diagram $(\ref{CRUDiag})$ is said to be 2-exact in $B$, if $\mathrm{Cok}(\underline{f})=0$ $($or equivalently $\mathrm{Ker}(\overline{g})=0$ $)$.
\end{definition}

\begin{remark}
In the notation of Lemma \ref{final16}, the following are equivalent :

{\rm (i)} $(\ref{CRUDiag})$ is 2-exact in $B$.

{\rm (ii)} $\underline{f}$ is fully cofaithful.

{\rm (iii)} $\overline{g}$ is fully faithful.

{\rm (iv)} $c(f)=\mathrm{cok}(k(g))$\quad(i.e. $\mathrm{Cok}(f)=\mathrm{Coim}(g)$).

{\rm (v)} $k(g)=\mathrm{ker}(c(f))$\quad(i.e. $\mathrm{Ker}(g)=\mathrm{Im}(f)$).\end{remark}
\begin{proof}
By the duality, we only show $\mathrm{(i)}\Leftrightarrow\mathrm{(iii)}\Leftrightarrow\mathrm{(v)}$.

$\mathrm{(i)}\Leftrightarrow\mathrm{(iii)}$ follows from Corollary \ref{RevisionLabel1}.

$\mathrm{(iii)}\Rightarrow\mathrm{(v)}$ follows from Proposition \ref{KerFullFaith}.

$\mathrm{(v)}\Rightarrow\mathrm{(iii)}$ follows from Proposition \ref{FirstFactorization}.
\end{proof}

Let us fix the notation for relative (co-)kernels of a complex.
\begin{definition}\label{final18}
For any complex $A_{\bullet}=(A_n,d_n,\delta_n)$ in $\mathbf{S}$, we put\\
{\rm (1)} $[Z^n(A_{\bullet}),z_n^A,\zeta_n^A]:=\mathrm{Ker}(d_n,\delta_{n+1})$.\\
{\rm (2)} $[Q^n(A_{\bullet}),q_n^A,\rho_n^A]:=\mathrm{Cok}(d_{n-1},\delta_{n-1})$.
\end{definition}

\begin{remark}\label{final19}
By the universality of $\mathrm{Ker}(d_n,\delta_{n+1})$ and Lemma \ref{KernelLemma}, there exist $k_n\in\mathbf{S}^1(A_{n-1},Z^n(A_{\bullet}))$, $\nu_{n,1}\in\mathbf{S}^2(k_n\circ z_n,d_{n-1})$ and $\nu_{n,2}\in\mathbf{S}^2(d_{n-2}\circ k_n,0)$ such that
\begin{eqnarray*}
(\nu_{n,1}\circ d_n)\cdot\delta_n&=&(k_n\circ\zeta_n)\cdot(k_n)_I^{\sharp}\\
(d_{n-2}\circ\nu_{n,1})\cdot\delta_{n-1}&=&(\nu_{n,2}\circ z_n)\cdot(z_n)_I^{\flat}.
\end{eqnarray*}
\[
\xy
(0,18)*+{Z^n (A_{\bullet})}="0";
(-36,0)*+{A_{n-2}}="2";
(-12,0)*+{A_{n-1}}="4";
(12,0)*+{A_n}="6";
(36,0)*+{A_{n+1}}="8";
(60,0)*+{A_{n+2}}="10";
{\ar|*+{_{d_{n-2}}} "2";"4"};
{\ar|*+{_{d_{n-1}}} "4";"6"};
{\ar^{k_n} "4";"0"};
{\ar|*+{_{d_n}} "6";"8"};
{\ar^{d_{n+1}} "8";"10"};
{\ar^{z_n} "0";"6"};
{\ar@{=>}^{\zeta_n} (15,3)*{};(18,8)*{}} ;
{\ar@{=>}_{\nu_{n,2}} (-15,3)*{};(-18,8)*{}} ;
{\ar@/_1.65pc/_{0} "2";"6"};
{\ar@/_1.65pc/_{0} "4";"8"};
{\ar@/^0.50pc/^{0} "0";"8"};
{\ar@/^0.50pc/^{0} "2";"0"};
{\ar@{=>}^{\delta_n} (12,-2)*{};(12,-6)*{}} ;
{\ar@{=>}_{\delta_{n-1}} (-12,-2)*{};(-12,-6)*{}} ;
{\ar@{=>}_{\nu_{n,1}} (0,10)*{};(0,4)*{}} ;
\endxy
\]
On the other hand, by the universality of $\mathrm{Ker}(d_n)$, we obtain a factorization diagram
\[
\xy
(0,-18)*+{\mathrm{Ker}(d_n)}="0";
(-36,0)*+{A_{n-2}}="2";
(-12,0)*+{A_{n-1}}="4";
(12,0)*+{A_n}="6";
(36,0)*+{A_{n+1}}="8";
(60,0)*+{A_{n+2}}="10";
{\ar|*+{_{d_{n-2}}} "2";"4"};
{\ar|*+{_{d_{n-1}}} "4";"6"};
{\ar|*+{_{\underline{d}_{n-1}}} "4";"0"};
{\ar|*+{_{d_n}} "6";"8"};
{\ar_{d_{n+1}} "8";"10"};
{\ar_{k(d_n)} "0";"6"};
{\ar@{=>}^{\varepsilon_{d_n}} (15,-3)*{};(18,-8)*{}} ;
{\ar@{=>}_{\underline{\delta}_{n-1}} (-15,-3)*{};(-18,-8)*{}} ;
{\ar@/^1.65pc/^{0} "2";"6"};
{\ar@/^1.65pc/^{0} "4";"8"};
{\ar@/_0.50pc/_{0} "0";"8"};
{\ar@/_0.50pc/_{0} "2";"0"};
{\ar@{=>}_{\delta_n} (12,2)*{};(12,6)*{}} ;
{\ar@{=>}^{\delta_{n-1}} (-12,2)*{};(-12,6)*{}} ;
{\ar@{=>}_{\underline{\delta}_n} (0,-10)*{};(0,-4)*{}} ;
\endxy
\]
which satisfy
\begin{eqnarray*}
(\underline{\delta}_n\circ d_n)\cdot\delta_n&=&(\underline{d}_{n-1}\circ\varepsilon_{d_n})\cdot(\underline{d}_{n-1})_I^{\sharp}\\
(d_{n-2}\circ\underline{\delta}_n)\cdot\delta_{n-1}&=&(\underline{\delta}_{n-1}\circ k(d_n))\cdot(k(d_n))_I^{\flat}
\end{eqnarray*}
By Proposition \ref{KerAndRelKer}, there exists a factorization of $z_n$ through $\mathrm{Ker}(d_n)$
\[
\xy
(-6,10)*+{Z^n (A_{\bullet})}="10";
(-6,-10)*+{\mathrm{Ker}(d_n)}="0";
(12,0)*+{A_n}="6";
(36,0)*+{A_{n+1}}="8";
{\ar^{d_n} "6";"8"};
{\ar_{k(d_n)} "0";"6"};
{\ar_{\underline{z}_n} "10";"0"};
{\ar@{=>}^{\varepsilon_{d_n}} (14,-3)*{};(16,-7)*{}};
{\ar@/_1.00pc/_{0} "0";"8"};
{\ar@/^1.00pc/^{0} "10";"8"};
{\ar^{z_n} "10";"6"};
{\ar@{=>}^{\zeta_n} (14,3)*{};(16,7)*{}};
{\ar@{=>}^{\underline{\zeta}_n} (-2,-4)*{};(2,2)*{}};
\endxy
\]
which satisfies
\[ (\underline{\zeta}_n\circ d_n)\cdot\zeta_n=(\underline{z}_n\circ\varepsilon_{d_n})\cdot(\underline{z}_n)_I^{\sharp}. \]
Moreover $\underline{z}_n$ is fully faithful by Proposition \ref{KerFactRelKer}.
\end{remark}

By the universality of $\mathrm{Ker}(d_n)$, we can show easily the following claim.

\begin{claim}\label{final20}
There exists a unique 2-cell $\widehat{\zeta}_n\in S(k_n\circ z_n,\underline{d}_{n-1})$
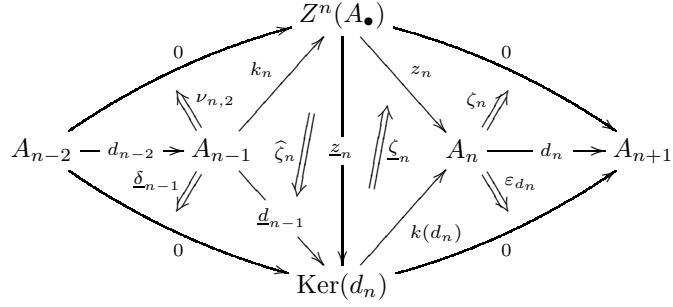
\[
\xy
(0,18)*+{Z^n (A_{\bullet})}="0";
(0,-18)*+{\mathrm{Ker}(d_n)}="20";
(-40,0)*+{A_{n-2}}="2";
(-16,0)*+{A_{n-1}}="4";
(16,0)*+{A_n}="6";
(40,0)*+{A_{n+1}}="8";
{\ar|*+{_{d_{n-2}}} "2";"4"};
{\ar^{k_n} "4";"0"};
{\ar|*+{_{d_n}} "6";"8"};
{\ar^{z_n} "0";"6"};
{\ar|*+{_{\underline{d}_{n-1}}} "4";"20"};
{\ar_{k(d_n)} "20";"6"};
{\ar|*+{_{\underline{z}_n}} "0";"20"};
{\ar@/^0.60pc/^{0} "0";"8"};
{\ar@/^0.60pc/^{0} "2";"0"};
{\ar@/_0.60pc/_{0} "20";"8"};
{\ar@/_0.60pc/_{0} "2";"20"};
{\ar@{=>}^{\zeta_n} (19,3)*{};(22,8)*{}};
{\ar@{=>}_{\nu_{n,2}} (-19,3)*{};(-22,8)*{}};
{\ar@{=>}^{\varepsilon_{d_n}} (19,-3)*{};(22,-8)*{}};
{\ar@{=>}_{\underline{\delta}_{n-1}} (-19,-3)*{};(-22,-8)*{}};
{\ar@{=>}_{\widehat{\zeta}_n} (-4,5)*{};(-6,-6)*{}};
{\ar@{=>}_{\underline{\zeta}_n} (4,-5)*{};(6,6)*{}};
\endxy
\]
such that
\[ (\widehat{\zeta}_n\circ k(d_n))\cdot\underline{\delta}_n=(k_n\circ\underline{\zeta}_n)\cdot\nu_{n,1}. \]
This $\widehat{\zeta}_n$ also satisfies
\[ (d_{n-2}\circ\widehat{\zeta}_n)\cdot\underline{\delta}_{n-1}=(\nu_{n,2}\circ\underline{z}_n)\cdot(\underline{z}_n)_I^{\flat}. \]
\end{claim}

\begin{remark}\label{final21}
Dually, by the universality of the cokernels, we obtain the following two factorization diagrams, where $\overline{q}_n$ is fully cofaithful.
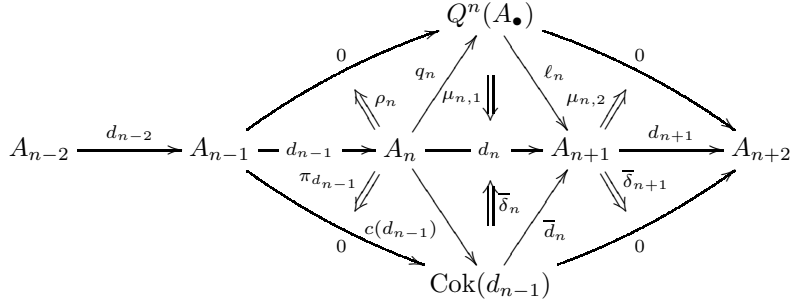
\[
\xy
(24,18)*+{Q^n(A_{\bullet})}="0";
(24,-18)*+{\mathrm{Cok}(d_{n-1})}="20";
(-36,0)*+{A_{n-2}}="2";
(-12,0)*+{A_{n-1}}="4";
(12,0)*+{A_n}="6";
(36,0)*+{A_{n+1}}="8";
(60,0)*+{A_{n+2}}="10";
{\ar^{d_{n-2}} "2";"4"};
{\ar|*+{_{d_{n-1}}} "4";"6"};
{\ar^{q_n} "6";"0"};
{\ar_{c(d_{n-1})} "6";"20"};
{\ar|*+{_{d_n}} "6";"8"};
{\ar^{d_{n+1}} "8";"10"};
{\ar^{\ell_n} "0";"8"};
{\ar_{\overline{d}_n} "20";"8"};
{\ar@{=>}^{\mu_{n,2}} (39,3)*{};(42,8)*{}};
{\ar@{=>}_{\rho_n} (9,3)*{};(6,8)*{}};
{\ar@{=>}^{\overline{\delta}_{n+1}} (39,-3)*{};(42,-8)*{}};
{\ar@{=>}_{\pi_{d_{n-1}}} (9,-3)*{};(6,-8)*{}};
{\ar@/^0.50pc/^{0} "0";"10"};
{\ar@/^0.50pc/^{0} "4";"0"};
{\ar@/_0.50pc/_{0} "20";"10"};
{\ar@/_0.50pc/_{0} "4";"20"};
{\ar@{=>}_{\mu_{n,1}} (24,10)*{};(24,4)*{}};
{\ar@{=>}_{\overline{\delta}_n} (24,-10)*{};(24,-4)*{}};
\endxy
\]
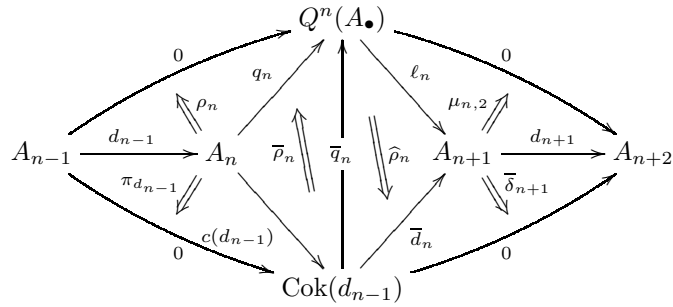
\[
\xy
(24,18)*+{Q^n (A_{\bullet})}="0";
(24,-18)*+{\mathrm{Cok}(d_{n-1})}="20";
(-16,0)*+{A_{n-1}}="4";
(8,0)*+{A_n}="6";
(40,0)*+{A_{n+1}}="8";
(64,0)*+{A_{n+2}}="10";
{\ar^{q_n} "6";"0"};
{\ar_{c(d_{n-1})} "6";"20"};
{\ar^{d_{n+1}} "8";"10"};
{\ar|*+{_{\overline{q}_n}} "20";"0"};
{\ar^{d_{n-1}} "4";"6"};
{\ar^{\ell_n} "0";"8"};
{\ar_{\overline{d}_n} "20";"8"};
{\ar@{=>}^{\mu_{n,2}} (43,3)*{};(46,8)*{}};
{\ar@{=>}_{\rho_n} (5,3)*{};(2,8)*{}};
{\ar@{=>}^{\overline{\delta}_{n+1}} (43,-3)*{};(46,-8)*{}};
{\ar@{=>}_{\pi_{d_{n-1}}} (5,-3)*{};(2,-8)*{}};
{\ar@/^0.60pc/^{0} "0";"10"};
{\ar@/^0.60pc/^{0} "4";"0"};
{\ar@/_0.60pc/_{0} "20";"10"};
{\ar@/_0.60pc/_{0} "4";"20"};
{\ar@{=>}^{\overline{\rho}_n} (20,-5)*{};(18,6)*{}};
{\ar@{=>}^{\widehat{\rho} _n} (28,5)*{};(30,-6)*{}};
\endxy
\]
\end{remark}

We define relative 2-cohomology in the following two ways, which will be shown
to be equivalent later.

\begin{definition}\label{final22}
\begin{eqnarray*}
H_1^n(A_{\bullet})&:=&\mathrm{Cok}(k_n,\nu_{n,2})\\
H_2^n(A_{\bullet})&:=&\mathrm{Ker}(\ell_n,\mu_{n,2})
\end{eqnarray*}
\end{definition}

\begin{lemma}\label{final23}
In the factorization diagram $(\ref{CRUDiag2})$ in Lemma \ref{final16}, if we take the cokernel of $\underline{f}$ and the kernel of $\overline{g}$, then there exist $w\in\mathbf{S}^1(\mathrm{Cok}(\underline{f}),\mathrm{Ker}(\overline{g}))$ and $\omega\in\mathbf{S}^2(c(\underline{f})\circ w\circ k(\overline{g}),k(g)\circ c(f))$ such that
\begin{eqnarray}
(\underline{f}\circ\omega)\cdot(\underline{\varphi}\circ c(f))\cdot\pi_f&=&(\pi_{\underline{f}}\circ w\circ k(\overline{g}))\cdot(w\circ k(\overline{g}))_I^{\flat}
\label{finalequation25-1}\\
(\omega\circ\overline{g})\cdot(k(g)\circ\overline{\varphi})\cdot\varepsilon_g&=&(c(\underline{f})\circ w\circ\varepsilon_{\overline{g}})\cdot(c(\underline{f})\circ w)_I^{\sharp}.
\label{finalequation25-2}
\end{eqnarray}
\[
\xy
(-16,16)*+{\mathrm{Ker}(g)}="0";
(16,16)*+{\mathrm{Cok}(f)}="2";
(-24,0)*+{A}="4";
(0,0)*+{B}="6";
(24,0)*+{C}="8";
(-10,32)*+{\mathrm{Cok}(\underline{f})}="10";
(10,32)*+{\mathrm{Ker}(\overline{g})}="12";
{\ar^{w} "10";"12"};
{\ar_{f} "4";"6"};
{\ar_{g} "6";"8"};
{\ar|*+{_{k(g)}} "0";"6"};
{\ar|*+{_{c(f)}} "6";"2"};
{\ar|*+{_{k(\overline{g})}} "12";"2"};
{\ar|*+{_{c(\underline{f})}} "0";"10"};
{\ar^{\underline{f}} "4";"0"};
{\ar^{\overline{g}} "2";"8"};
{\ar@/_1.80pc/_{0} "4";"8"};
{\ar@/^2.40pc/^{0} "4";"10"};
{\ar@/^2.40pc/^{0} "12";"8"};
{\ar@{=>}_{\underline{\varphi}} (-15,8)*{};(-15,4)*{}};
{\ar@{=>}_{\overline{\varphi}} (15,8)*{};(15,4)*{}};
{\ar@{=>}_{\varphi} (0,-2)*{};(0,-6)*{}};
{\ar@{=>}_{\omega} (0,20)*{};(0,14)*{}};
{\ar@{=>}_{\pi_{\underline{f}}} (-22,18)*{};(-24,20)*{}};
{\ar@{=>}^{\varepsilon_{\overline{g}}} (22,18)*{};(24,20)*{}};
\endxy
\]
Moreover, for any other factorization $(w^{\prime},\omega^{\prime})$ with these
properties, there exists a unique 2-cell $\kappa\in\mathbf{S}^2(w,w^{\prime})$ such that $(c(\underline{f})\circ\kappa\circ k(\overline{g}))\cdot\omega^{\prime}=\omega$.
\end{lemma}
\begin{proof}
Applying Proposition \ref{final5} to 
\begin{equation}
\xy
(-9,6)*+{A}="0";
(-9,-6)*+{\mathrm{Ker}(g)}="2";
(9,6)*+{A}="4";
(9,-6)*+{B}="6";
(11,-7)*+{,}="7";
{\ar@{=}^{\mathrm{id}_A} "0";"4"};
{\ar_{\underline{f}} "0";"2"};
{\ar_{k(g)} "2";"6"};
{\ar^{f} "4";"6"};
{\ar@{=>}_{\overline{\varphi}^{-1}} (0,2)*{};(0,-2)*{}};
\endxy
\label{CRUDiag3}
\end{equation}
we obtain $w_1\in\mathbf{S}^1(\mathrm{Cok}(\underline{f}),\mathrm{Cok}(f))$ and $\omega_1\in\mathbf{S}^2(c(\underline{f})\circ w_1,k(g)\circ c(f))$ which satisfy
\begin{equation}
(\underline{f}\circ\omega_1)\cdot(\underline{\varphi}\circ c(f))\cdot\pi_f=(\pi_{\underline{f}}\circ w_1)\cdot(w_1)_I^{\flat}.
\label{finalm}
\end{equation}
Then $(\omega_1\circ\overline{g})\cdot(k(g)\circ\overline{\varphi})\cdot\varepsilon_g\in\mathbf{S}^2(c(\underline{f})\circ w_1\circ\overline{g},0)$ becomes compatible with $\pi_{\underline{f}}$.
\[
\xy
(-30,0)*+{A}="2";
(-12,0)*+{\mathrm{Ker}(g)}="3";
(12,0)*+{\mathrm{Cok}(\underline{f})}="4";
(30,0)*+{C}="5";
{\ar_{\underline{f}}  "2";"3"};
{\ar_{c(\underline{f})} "3";"4"};
{\ar^{w_1\circ\overline{g}} "4";"5"};
{\ar@/^1.65pc/^{0} "2";"4"};
{\ar@/_1.65pc/_{0} "3";"5"};
{\ar@{=>}^{\pi_{\underline{f}}} (-8,2)*{};(-8,6)*{}};
{\ar@{=>}^{} (8,-2)*{};(8,-6)*{}};
{\ar@{..}@/_0.80pc/^{} (8,-4)*{};(24,-10)*{_{(\omega_1\circ\overline{g})\cdot(k(g)\circ\overline{\varphi})\cdot\varepsilon_g}}};
\endxy
\]
By Lemma \ref{KernelLemma}, there exists a 2-cell $\delta\in S(w_1\circ\overline{g},0)$ such that
\[ (c(\underline{f})\circ\delta)\cdot c(\underline{f})_I^{\sharp}=(\omega_1\circ\overline{g})\cdot(k(g)\circ\overline{\varphi})\cdot\varepsilon_g. \]
So, if we take the cokernels of $k(g)$ and $w_1$, then by Proposition \ref{final9}, we obtain the following diagram:
\[
\xy
(-30,-8)*+{\mathrm{Ker}(g)}="2";
(-10,-8)*+{B}="6";
(-30,-24)*+{\mathrm{Cok}(\underline{f})}="10";
(-10,-24)*+{\mathrm{Cok}(f)}="12";
(30,-8)*+{C}="20";
(30,-24)*+{C}="22";
(10,-8)*+{\mathrm{Coim}(g)}="24";
(10,-24)*+{\mathrm{Cok}(w_1)}="26";
{\ar^{c(k(g))} "6";"24"};
{\ar_{c(w_1)} "12";"26"};
{\ar^{j(g)} "24";"20"};
{\ar_{\exists g^{\dagger}} "26";"22"};
{\ar_{k(g)} "2";"6"};
{\ar@{=}^{\mathrm{id}_C} "20";"22"};
{\ar_{w_1} "10";"12"};
{\ar_{\exists\overline{c}} "24";"26"};
{\ar_{c(\underline{f})} "2";"10"};
{\ar|*+{_{c(f)}} "6";"12"};
{\ar@/^1.80pc/^{g} "6";"20"};
{\ar@/_1.80pc/_{\overline{g}} "12";"22"};
{\ar@{=>}^{\exists\overline{\varphi}_2} (18,-16)*{};(22,-16)*{}};
{\ar@{=>}^{\exists\overline{\varphi}_1} (-2,-16)*{};(2,-16)*{}};
{\ar@{=>}_{\omega_1} (-22,-16)*{};(-18,-16)*{}};
{\ar@{=>}_{\mu_g} (10,-6)*{};(10,-2)*{}};
{\ar@{=>}_{^{\exists}} (10,-26)*{};(10,-30)*{}};
\endxy
\]
Applying Proposition \ref{final8} to $(\ref{CRUDiag3})$, we obtain
\[ [\mathrm{Cok}(w_1),\overline{c},(\overline{c})_I^{\flat}]=\mathrm{Cok}(0\overset{0}{\longrightarrow}\mathrm{Coim}(g)). \]
Thus $\overline{c}$ is an equivalence. Since $j(g)$ is fully faithful, $g^{\dag}$ becomes fully faithful. Thus the following diagram is 2-exact in $\mathrm{Cok}(f)$.
\begin{equation}
\xy
(-20,0)*+{\mathrm{Cok}(\underline{f})}="0";
(0,0)*+{\mathrm{Cok}(f)}="2";
(20,0)*+{C}="4";
{\ar^{w_1} "0";"2"};
{\ar^{\overline{g}} "2";"4"};
{\ar@/_1.80pc/_{0} "0";"4"};
{\ar@{=>}_{\delta} (0,-2)*{};(0,-6)*{}};
\endxy
\label{finaldia2}
\end{equation}
So if we factor $(\ref{finaldia2})$ by $w\in\mathbf{S}^1(\mathrm{Cok}(\underline{f}),\mathrm{Ker}(\overline{g}))$ and $\omega_2\in\mathbf{S}^2(w\circ k(\overline{g}),w_1)$ as in the diagram
\begin{equation}
\xy
(-20,0)*+{\mathrm{Cok}(\underline{f})}="10";
(0,0)*+{\mathrm{Cok}(f)}="12";
(-10,14)*+{\mathrm{Ker}(\overline{g})}="14";
(20,0)*+{C}="22";
{\ar_{w_1} "10";"12"};
{\ar^{w} "10";"14"};
{\ar^{k(\overline{g})} "14";"12"};
{\ar_{\overline{g}} "12";"22"};
{\ar@/_1.80pc/_{0} "10";"22"};
{\ar@/^1.20pc/^{0} "14";"22"};
{\ar@{=>}_{\delta} (0,-2)*{};(0,-6)*{}};
{\ar@{=>}^{\omega_2} (-10,8)*{};(-10,4)*{}};
{\ar@{=>}_{\varepsilon_{\overline{g}}} (4,4)*{};(8,8)*{}};
\endxy
\label{finaldia3}
\end{equation}
which satisfies
\[ (w\circ\varepsilon_{\overline{g}})\cdot w_I^{\sharp}=(\omega_2\circ\overline{g})\cdot\delta, \]
then $w$ becomes fully cofaithful by Lemma \ref{final16}. If we put $\omega:=(c(\underline{f})\circ\omega_2)\cdot\omega_1$, then $(w,\omega)$ satisfies conditions $(\ref{finalequation25-1})$ and $(\ref{finalequation25-2})$.

If $(w^{\prime},\omega^{\prime})$ satisfies
\begin{eqnarray}
(\underline{f}\circ\omega^{\prime})\cdot(\underline{\varphi}\circ c(f))\cdot\pi_{f}&=&(\pi_{\underline{f}}\circ w^{\prime}\circ k(\overline{g}))\cdot(w^{\prime}\circ k(\overline{g}))_I^{\flat}
\label{finalequation28-1}\\
(\omega^{\prime}\circ\overline{g})\cdot(k(g)\circ\overline{\varphi})\cdot\varepsilon_{g}&=&(c(\underline{f})\circ w^{\prime}\circ\varepsilon_{\overline{g}})\cdot(c(\underline{f})\circ w^{\prime})_I^{\sharp},
\nonumber
\end{eqnarray}
then, since both the factorization of $k(g)\circ c(f)$ through $\mathrm{Cok}(\underline{f})$
\begin{eqnarray*}
\omega^{\prime}&:&c(\underline{f})\circ w^{\prime}\circ k(\overline{g})\Longrightarrow k(g)\circ c(f)\\
\omega_1&:&c(\underline{f})\circ w_1\Longrightarrow k(g)\circ c(f)
\end{eqnarray*}
are compatible with $\pi_{\underline{f}}$ and $(\underline{\varphi}\circ c(f))\cdot\pi_f$ by $(\ref{finalm})$ and $(\ref{finalequation28-1})$, there exists $\omega_2^{\prime}\in\mathbf{S}^2(w^{\prime}\circ k(\overline{g}),w_1)$ such that
\[ (c(\underline{f})\circ\omega_2^{\prime})\cdot\omega_1=\omega^{\prime}. \]
Then we can see $\omega_2^{\prime}$ is compatible with $\varepsilon_{\overline{g}}$ and $\delta$. So, comparing this with the factorization $(\ref{finaldia3})$, by the universality of $\mathrm{Ker}(\overline{g})$, we see there exists a unique 2-cell $\kappa\in\mathbf{S}^2(w,w^{\prime})$ such that $(\kappa\circ k(\overline{g}))\cdot\omega_2^{\prime}=\omega_2$. Then $\kappa$ satisfies $(c(\underline{f})\circ\kappa\circ k(\overline{g}))\cdot\omega^{\prime}=\omega$. Uniqueness of such $\kappa$ follows from the fact that $c(\underline{f})$ is cofaithful and $k(\overline{g})$ is faithful.
\end{proof}

\begin{proposition}\label{final25}
In Lemma \ref{final23}, $w$ is an equivalence.
\end{proposition}
\begin{proof}
We showed Lemma \ref{final23} by taking the cokernel first and the kernel second, but we obtain the same $(w,\omega)$ if we take the kernel first and the cokernel second, because of
the symmetricity of the statement $($and the uniqueness of $(w,\omega)$ up to
an equivalence$)$ of Lemma \ref{final23}. As shown in the proof, since $(\ref{finaldia2})$ is 2-exact in $\mathrm{Cok}(\underline{f})$, $w$ becomes fully cofaithful in the factorization $(\ref{finaldia3})$. By the above remark, similarly $w$ can be obtained also by the factorization
\[
\xy
(-19,0)*+{A}="10";
(0,0)*+{\mathrm{Ker}(g)}="12";
(10,14)*+{\mathrm{Cok}(\underline{f})}="14";
(19,0)*+{\mathrm{Ker}(\overline{g})}="22";
{\ar_{} "10";"12"};
{\ar^{} "12";"14"};
{\ar^{w} "14";"22"};
{\ar_{} "12";"22"};
{\ar@{=>}^{} (9.5,8)*{};(9.5,4)*{}};
{\ar@/_1.80pc/_{0} "10";"22"};
{\ar@{=>}_{} (0,-2)*{};(0,-6)*{}};
\endxy
\]
where the bottom row is 2-exact in $\mathrm{Ker}(g)$. So $w$ becomes fully faithful. Thus, $w$ is fully cofaithful and fully faithful, i.e. an equivalence.
\end{proof}

\begin{corollary}\label{final26}
For any complex $A_{\bullet}=(A_n,d_n,\delta_n)$, if we factor it as
\[
\xy
(-16,16)*+{Z^n(A_{\bullet})}="0";
(16,16)*+{Q^n(A_{\bullet})}="2";
(-24,0)*+{A_{n-1}}="4";
(0,0)*+{A_n}="6";
(24,0)*+{A_{n+1}}="8";
(-44,0)*+{A_{n-2}}="14";
(44,0)*+{A_{n+2}}="18";
(-10,32)*+{H^n_1(A_{\bullet})}="10";
(10,32)*+{H^n_2(A_{\bullet})}="12";
{\ar_{d_{n-1}} "4";"6"};
{\ar_{d_n} "6";"8"};
{\ar^{z_n} "0";"6"};
{\ar^{q_n} "6";"2"};
{\ar_{d_{n-2}} "14";"4"};
{\ar_{d_{n+2}} "8";"18"};
{\ar|*+{_{k(\ell_n,\mu_{n,2})}} "12";"2"};
{\ar|*+{_{c(k_n,\nu_{n,2})}} "0";"10"};
{\ar^{k_n} "4";"0"};
{\ar^{\ell_n} "2";"8"};
{\ar@{=>}^{\nu_{n,1}} (-16,8)*{};(-16,3)*{}};
{\ar@{=>}_{\mu_{n,1}} (16,8)*{};(16,3)*{}};
{\ar@/_1.80pc/_{0} "4";"8"};
{\ar@/^1.90pc/^{0} "4";"10"};
{\ar@/^1.90pc/^{0} "12";"8"};
{\ar@/^0.75pc/^{0} "14";"0"};
{\ar@/^0.75pc/^{0} "2";"18"};
{\ar@{=>}_{\delta^A_n} (0,-3)*{};(0,-6)*{}};
{\ar@{=>}^{\nu_{n,2}} (-28,4)*{};(-31,9)*{}};
{\ar@{=>}_{\mu_{n,2}} (28,3)*{};(31,9)*{}};
{\ar@{=>} (-18,19)*{};(-22,21)*{}};
{\ar@{=>} (18,19)*{};(22,21)*{}};
\endxy
\]
$($in the notation of Definition \ref{final18}, Remark \ref{final19} and Remark \ref{final21}$)$, then there exist $w\in\mathbf{S}^1(H_1^n(A_{\bullet}),H_2^n(A_{\bullet}))$ and $\omega\in\mathbf{S}^2(c(k_n,\nu_{n,2})\circ w\circ k(\ell_n,\mu_{n,2}),z_n\circ q_n)$ such that
\begin{eqnarray*}
(k_n\circ\omega)\cdot(\nu_{n,1}\circ q_n)\cdot\rho_n&=&(\pi_{(k_n,\nu_{n,2})}\circ w\circ k(\ell_n,\mu_{n,2}))\cdot(w\circ k(\ell_n,\mu_{n,2}))_I^{\flat}\\
(\omega\circ\ell_n)\cdot(z_n\circ\mu_{n,1})\cdot\zeta_n&=&(c(k_n,\nu_{n,2})\circ w\circ\varepsilon_{(\ell_n,\mu_{n,2})})\cdot(c(k_n,\nu_{n,2})\circ w)_I^{\sharp}.
\end{eqnarray*}
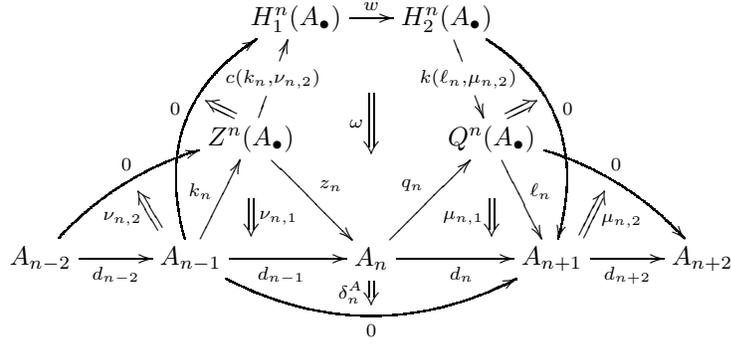
\[
\xy
(-16,16)*+{Z^n(A_{\bullet})}="0";
(16,16)*+{Q^n(A_{\bullet})}="2";
(-24,0)*+{A_{n-1}}="4";
(0,0)*+{A_n}="6";
(24,0)*+{A_{n+1}}="8";
(-44,0)*+{A_{n-2}}="14";
(44,0)*+{A_{n+2}}="18";
(-10,32)*+{H^n _1(A_{\bullet})}="10";
(10,32)*+{H^n _2(A_{\bullet})}="12";
{\ar^{w} "10";"12"};
{\ar_{d_{n-1}} "4";"6"};
{\ar_{d_n} "6";"8"};
{\ar^{z_n} "0";"6"};
{\ar^{q_n} "6";"2"};
{\ar_{d_{n-2}} "14";"4"};
{\ar_{d_{n+2}} "8";"18"};
{\ar|*+{_{k(\ell_n,\mu_{n,2})}} "12";"2"};
{\ar|*+{_{c(k_n,\nu_{n,2})}} "0";"10"};
{\ar^{k_n} "4";"0"};
{\ar^{\ell_n} "2";"8"};
{\ar@/_1.80pc/_{0} "4";"8"};
{\ar@/^1.90pc/^{0} "4";"10"};
{\ar@/^1.90pc/^{0} "12";"8"};
{\ar@/^0.75pc/^{0} "14";"0"};
{\ar@/^0.75pc/^{0} "2";"18"};
{\ar@{=>}^{\nu_{n,1}} (-16,8)*{};(-16,3)*{}};
{\ar@{=>}_{\mu_{n,1}} (16,8)*{};(16,3)*{}};
{\ar@{=>}_{\delta^A_n} (0,-3)*{};(0,-6)*{}};
{\ar@{=>}_{\omega} (0,22)*{};(0,14)*{}};
{\ar@{=>}^{\nu_{n,2}} (-28,4)*{};(-31,9)*{}};
{\ar@{=>}_{\mu_{n,2}} (28,3)*{};(31,9)*{}};
{\ar@{=>} (-18,19)*{};(-22,21)*{}};
{\ar@{=>} (18,19)*{};(22,21)*{}};
\endxy
\]
For any other factorization $(w^{\prime},\omega^{\prime})$ with these conditions, there exists a unique 2-cell $\kappa\in\mathbf{S}^2(\omega,\omega^{\prime})$ such that $(c(k_n,\nu_{n,2})\circ\kappa\circ k(\ell_n,\mu_{n,2}))\cdot\omega^{\prime}=\omega$. Moreover, this $w$ becomes an equivalence.
\end{corollary}

\begin{proof}
For the factorization diagrams
\begin{eqnarray*}
&
\xy
(24,-18)*+{\mathrm{Cok}(d_{n-2})}="20";
(-12,0)*+{A_{n-2}}="4";
(12,0)*+{A_{n-1}}="6";
(36,0)*+{A_n}="8";
(60,0)*+{A_{n+1}}="10";
{\ar|*+{_{d_{n-2}}} "4";"6"};
{\ar|*+{_{c(d_{n-2})}} "6";"20"};
{\ar|*+{_{d_{n-1}}} "6";"8"};
{\ar^{d_n} "8";"10"};
{\ar_{\overline{d}_{n-1}} "20";"8"};
{\ar@/_0.50pc/_{0} "20";"10"};
{\ar@/_0.50pc/_{0} "4";"20"};
{\ar@/^1.65pc/^{0} "4";"8"};
{\ar@/^1.65pc/^{0} "6";"10"};
{\ar@{=>}^{\overline{\delta}_n} (39,-3)*{};(42,-8)*{}} ;
{\ar@{=>}_{\pi_{d_{n-2}}} (9,-3)*{};(6,-8)*{}} ;
{\ar@{=>}^{\delta_{n-1}} (12,2)*{};(12,6)*{}} ;
{\ar@{=>}_{\delta_n} (36,2)*{};(36,6)*{}} ;
{\ar@{=>}_{\overline{\delta}_{n-1}} (24,-9)*{};(24,-4)*{}} ;
\endxy
\\
&
\xy
(0,-18)*+{\mathrm{Ker}(d_{n+1})}="0";
(-36,0)*+{A_{n-1}}="2";
(-12,0)*+{A_n}="4";
(12,0)*+{A_{n+1}}="6";
(36,0)*+{A_{n+2}}="8";
{\ar|*+{_{d_{n-1}}} "2";"4"};
{\ar|*+{_{d_n}} "4";"6"};
{\ar|*+{_{\underline{d}_n}} "4";"0"};
{\ar|*+{_{d_{n+1}}} "6";"8"};
{\ar_{k(d_{n+1})} "0";"6"};
{\ar@/^1.65pc/^{0} "2";"6"};
{\ar@/^1.65pc/^{0} "4";"8"};
{\ar@/_0.50pc/_{0} "0";"8"};
{\ar@/_0.50pc/_{0} "2";"0"};
{\ar@{=>}^{\varepsilon_{d_{n+1}}} (15,-3)*{};(18,-8)*{}};
{\ar@{=>}_{\underline{\delta}_n} (-15,-3)*{};(-18,-8)*{}};
{\ar@{=>}_{\delta_{n+1}} (12,2)*{};(12,6)*{}};
{\ar@{=>}^{\delta_n} (-12,2)*{};(-12,6)*{}};
{\ar@{=>}_{\underline{\delta}_{n+1}} (0,-9)*{};(0,-4)*{}};
\endxy
\end{eqnarray*}
which satisfy
\begin{eqnarray*}
(d_{n-2}\circ\overline{\delta}_{n-1})\cdot\delta_{n-1}&=&(\pi_{d_{n-2}}\circ\overline{d}_{n-1})\cdot(\overline{d}_{n-1})_I^{\flat}\\
(\overline{\delta}_{n-1}\circ d_n)\cdot\delta_n&=&(c(d_{n-2})\circ\overline{\delta}_n)\cdot c(d_{n-2})_I^{\sharp}\\
(\underline{\delta}_{n+1}\circ d_{n+1})\cdot\delta_{n+1}&=&(\underline{d}_n\circ\varepsilon_{d_{n+1}})\cdot(\underline{d}_n)_I^{\sharp}\\
(d_{n-1}\circ\underline{\delta}_{n+1})\cdot\delta_n&=&(\underline{\delta}_n\circ k(d_{n+1}))\cdot k(d_{n+1})_I^{\flat},
\end{eqnarray*}
there exists a unique 2-cell $\delta_n^{\dag}\in\mathbf{S}^2(\overline{d}_{n-1}\circ\underline{d}_n,0)$ such that
\begin{eqnarray*}
(\delta_n^{\dag}\circ k(d_{n+1}))\cdot k(d_{n+1})_I^{\flat}&=&(\overline{d}_{n-1}\circ\underline{\delta}_{n+1})\cdot\overline{\delta}_n\\
(c(d_{n-2})\circ\delta_n^{\dag})\cdot c(d_{n-2})_I^{\sharp}&=&(\overline{\delta}_{n-1}\circ\underline{d}_n)\cdot\underline{\delta}_n.
\end{eqnarray*}
By Proposition \ref{KerAndRelKer}, applying Lemma \ref{final23} and Proposition \ref{final25} to the following diagram, we can obtain Corollary \ref{final26}.
\[
\xy
(-14,14)*+{Z^n(A_{\bullet})}="0";
(14,14)*+{Q^n(A_{\bullet})}="2";
(-28,0)*+{\mathrm{Cok}(d_{n-2})}="4";
(0,0)*+{A_n}="6";
(28,0)*+{\mathrm{Ker}(d_{n+1})}="8";
{\ar_{\overline{d}_{n-1}} "4";"6"};
{\ar_{\underline{d}_n } "6";"8"};
{\ar^{z_n} "0";"6"};
{\ar^{q_n} "6";"2"};
{\ar^{} "4";"0"};
{\ar^{} "2";"8"};
{\ar@{=>}_{} (-14,8)*{};(-14,3)*{}};
{\ar@{=>}_{} (14,8)*{};(14,3)*{}};
{\ar@/_1.80pc/_{0} "4";"8"};
{\ar@{=>}_{\delta^{\dag}_n} (0,-3)*{};(0,-6)*{}};
\endxy
\]
\end{proof}

Thus $H_1^n(A_{\bullet})$ and $H_2^n(A_{\bullet})$ are equivalent. We abbreviate this to $H^n(A_{\bullet})$.

\begin{definition}
A complex $A_{\bullet}$ is said to be relatively 2-exact in $A_n$ if $H^n(A_{\bullet})$ is equivalent to zero.
\end{definition}

\begin{remark}
If the complex is bounded, we consider the relative 2-exactness after adding zeroes as in Remark \ref{RemRemComplex}. For example, a bounded complex
\[
\xy
(-8,0)*+{A}="4";
(8,0)*+{B}="6";
(24,0)*+{C}="8";
{\ar_{f} "4";"6"};
{\ar_{g} "6";"8"};
{\ar@/^1.65pc/^{0} "4";"8"};
{\ar@{=>}_{\varphi} (8,2)*{};(8,6)*{}};
\endxy
\]
is relatively 2-exact in $B$ if and only if
\[
\xy
(-24,0)*+{0}="2";
(-8,0)*+{A}="4";
(8,0)*+{B}="6";
(24,0)*+{C}="8";
(40,0)*+{0}="10";
{\ar_{0} "2";"4"};
{\ar_{f} "4";"6"};
{\ar_{g} "6";"8"};
{\ar_{0} "8";"10"};
{\ar@/_1.65pc/_{0} "2";"6"};
{\ar@/^1.65pc/^{0} "4";"8"};
{\ar@/_1.65pc/_{0} "6";"10"};
{\ar@{=>}_{f^{\flat}_I} (-8,-2)*{};(-8,-6)*{}};
{\ar@{=>}_{g^{\sharp}_I} (24,-2)*{};(24,-6)*{}};
{\ar@{=>}_{\varphi} (8,2)*{};(8,6)*{}};
\endxy
\]
is relatively 2-exact in $B$, and this is equivalent to the 2-exactness in $B$ by Remark \ref{3Rem}.
\end{remark}

\section{Long cohomology sequence in a relatively exact 2-category}\label{Long Cohomology Sequence}

\paragraph{Diagram lemmas {\rm (2)}}\quad\\

\begin{lemma}\label{final28}

Let $A_{\bullet}$ be a complex in $\mathbf{S}$, in which $A_5=0$ and $d_4=0:$
\begin{equation}
\xy
(-32,0)*+{A_1}="1";
(-16,0)*+{A_2}="2";
(0,0)*+{A_3}="3";
(16,0)*+{A_4}="4";
(32,0)*+{0}="5";
{\ar_{d^A_1} "1";"2"};
{\ar_{d^A_2} "2";"3"};
{\ar_{d^A_3} "3";"4"};
{\ar_{0} "4";"5"};
{\ar@/^1.80pc/^{0} "1";"3"};
{\ar@/_1.80pc/_{0} "2";"4"};
{\ar@/^1.80pc/^{0} "3";"5"};
{\ar@{=>}_<<<{\delta_2^A} (-18,2)*{};(-18,6)*{}};
{\ar@{=>}^<<<{\delta_3^A} (0,-2)*{};(0,-6)*{}};
{\ar@{=>}^{(d_3^A)^{\sharp}_I} (18,2)*{};(18,6)*{}};
\endxy
\label{4-3-1}
\end{equation}
Then, $(\ref{4-3-1})$ is relatively 2-exact in $A_3$ and $A_4$ if and only if $\mathrm{Cok}(d_2,\delta_2)=A_4$, i.e. $[Q^3(A_{\bullet}),q_3,\rho_3]=[A_4,d_3,\delta_3]$.
\end{lemma}
\begin{proof}
As in Remark \ref{final21}, we have two factorization diagrams
\begin{eqnarray*}
\xy
(24,18)*+{Q^3 (A_{\bullet})}="0";
(24,-18)*+{\mathrm{Cok}(d_2)}="20";
(-8,0)*+{A_2}="4";
(12,0)*+{A_3}="6";
(36,0)*+{A_4}="8";
{\ar|*+{_{d_2}} "4";"6"};
{\ar^{q_3} "6";"0"};
{\ar_{c(d_2)} "6";"20"};
{\ar|*+{_{d_3}} "6";"8"};
{\ar^{\ell_3} "0";"8"};
{\ar_{\overline{d}_3} "20";"8"};
{\ar@{=>}_{\rho_3} (9,3)*{};(6,8)*{}} ;
{\ar@{=>}_{\pi_{d_2}} (9,-3)*{};(6,-8)*{}} ;
{\ar@/^0.80pc/^{0} "4";"0"};
{\ar@/_0.80pc/_{0} "4";"20"};
{\ar@{=>}_{\mu_{3,1}} (24,10)*{};(24,4)*{}} ;
{\ar@{=>}_{\overline{\delta}_3} (24,-10)*{};(24,-4)*{}};
\endxy
\quad
\xy
(24,18)*+{Q^3 (A_{\bullet})}="0";
(24,-18)*+{\mathrm{Cok}(d_2)}="20";
(-8,0)*+{A_2}="4";
(10,0)*+{A_3}="6";
(38,0)*+{A_4}="8";
{\ar^{q_3} "6";"0"};
{\ar_{c(d_2)} "6";"20"};
{\ar|*+{_{\overline{q}_3}} "20";"0"};
{\ar^{d_2} "4";"6"};
{\ar^{\ell_3} "0";"8"};
{\ar_{\overline{d}_3} "20";"8"};
{\ar@{=>}_{\rho_3} (7,3)*{};(4,8)*{}};
{\ar@{=>}_{\pi_{d_2}} (7,-3)*{};(4,-8)*{}};
{\ar@/^0.80pc/^{0} "4";"0"};
{\ar@/_0.80pc/_{0} "4";"20"};
{\ar@{=>}^{\overline{\rho}_3} (20,-8)*{};(18,4)*{}};
{\ar@{=>}^{\widehat{\rho}_3} (28,8)*{};(30,-4)*{}};
\endxy
\end{eqnarray*}
where $\overline{q}_3$ is fully cofaithful. We have
\begin{eqnarray*}
(\ref{4-3-1})\ \text{is relatively 2-exact in}\ A_4&\Leftrightarrow&\mathrm{Cok}(d_3,\delta_3)=0\\
&\underset{\text{Prop. \ref{KerAndRelKer}}}{\Leftrightarrow}&\mathrm{Cok}(\overline{d}_3)=0\Leftrightarrow\mathrm{Cok}(\overline{q}_3\circ\ell_3)=0\\
&\underset{\text{Prop. \ref{KerFullFaith}}}{\Leftrightarrow}&\mathrm{Cok}(\ell_3)=0\Leftrightarrow\ell_3\ \text{is fully cofaithful}
\end{eqnarray*}
and
\begin{eqnarray*}
(\ref{4-3-1})\ \text{is relatively 2-exact in}\ A_3&\Leftrightarrow&\mathrm{Ker}(\ell_3,(\ell_3)_I^{\sharp})=0\\
&\underset{\text{Rem. \ref{3Rem}}}{\Leftrightarrow}&\mathrm{Ker}(\ell_3)=0\Leftrightarrow\ell_3\ \text{is fully faithful}.
\end{eqnarray*}
Thus, $(\ref{4-3-1})$ is relatively 2-exact in $A_3$ and $A_4$ if and only if $\ell$ is fully cofaithful and fully faithful, i.e. $\ell$ is an equivalence.
\end{proof}

By Remark \ref{3Rem}, we have the following corollary:
\begin{corollary}\label{final29}
Let $(A_n,d_n,\delta_n)$ be a bounded complex in $\mathbf{S}$, as follows$:$
\begin{equation}
\xy
(-32,0)*+{A_1}="1";
(-16,0)*+{A_2}="2";
(0,0)*+{A_3}="3";
(16,0)*+{0}="4";
{\ar_{d^A_1} "1";"2"};
{\ar_{d^A_2} "2";"3"};
{\ar^{0} "3";"4"};
{\ar@/^1.80pc/^{0} "1";"3"};
{\ar@/_1.80pc/_{0} "2";"4"};
{\ar@{=>}_<<<{\delta_2^A} (-18,2)*{};(-18,6)*{}};
{\ar@{=>}^{(d^A_2)^{\sharp}_I} (0,-2)*{};(0,-6)*{}};
\endxy
\label{4-3-2}
\end{equation}
Then, $(\ref{4-3-2})$ is relatively 2-exact in $A_2$ and $A_3$ if and only if $\mathrm{Cok}(d_1^A)=[A_3,d_2^A,\delta_2^A]$.
\end{corollary}

\begin{lemma}\label{FINLem1}
Let $A_{\bullet}$ be a complex. As in Definition \ref{final18}, Remark \ref{final19} and Remark \ref{final21}, take a factorization diagram
\[
\xy
(24,18)*+{Z^{n+1}(A_{\bullet})}="0";
(24,-18)*+{Q^n(A_{\bullet})}="20";
(-12,0)*+{A_{n-1}}="4";
(12,0)*+{A_n}="6";
(36,0)*+{A_{n+1}}="8";
(60,0)*+{A_{n+2}}="10";
{\ar|*+{_{d_{n-1}}} "4";"6"};
{\ar^{k_{n+1}} "6";"0"};
{\ar_{q_n} "6";"20"};
{\ar|*+{_{d_n}} "6";"8"};
{\ar^{d_{n+1}} "8";"10"};
{\ar|*+{_{z_{n+1}}} "0";"8"};
{\ar_{\ell_n} "20";"8"};
{\ar@{=>}^{\zeta_{n+1}} (39,3)*{};(42,8)*{}};
{\ar@{=>}^{\nu_{n+1,2}} (9,3)*{};(6,8)*{}};
{\ar@{=>}^{\mu_{n,2}} (39,-3)*{};(42,-8)*{}};
{\ar@{=>}_{\rho_n} (9,-3)*{};(6,-8)*{}};
{\ar@/^0.60pc/^{0} "0";"10"};
{\ar@/^0.60pc/^{0} "4";"0"};
{\ar@/_0.60pc/_{0} "20";"10"};
{\ar@/_0.60pc/_{0} "4";"20"};
{\ar@{=>}|*+{_{\nu_{n+1,1}}} (24,10)*{};(24,2)*{}};
{\ar@{=>}_{\mu_{n,1}} (24,-10)*{};(24,-4)*{}};
\endxy
\]
which satisfies
\begin{eqnarray*}
(\nu_{n+1,1}\circ d_{n+1})\cdot\delta_{n+1}&=&(k_{n+1}\circ\zeta_{n+1})\cdot(k_{n+1})_I^{\sharp}\\
(d_{n-1}\circ\nu_{n+1,1})\cdot\delta_n&=&(\nu_{n+1,2}\circ z_{n+1})\cdot(z_{n+1})_I^{\flat}\\
(d_{n-1}\circ\mu_{n,1})\cdot\delta_n&=&(\rho_n\circ\ell_n)\cdot(\ell_n)_I^{\flat}\\
(\mu_{n,1}\circ d_{n+1})\cdot\delta_{n+1}&=&(q_n\circ\mu_{n,2})\cdot(q_n)_I^{\sharp}.
\end{eqnarray*}
Then, there exist $x_n\in\mathbf{S}^1(Q^n(A_{\bullet}),Z^{n+1}(A_{\bullet}))$, $\xi_n\in\mathbf{S}^2(x_n\circ z_{n+1},\ell_n)$ and $\eta_n\in\mathbf{S}^2(q_n\circ x_n,k_{n+1})$ such that
\begin{eqnarray}
(\xi_n\circ d_{n+1})\cdot\mu_{n,2}&=&(x_n\circ\zeta_{n+1})\cdot(x_n)_I^{\sharp}\nonumber\\
(q_n\circ\xi_n)\cdot\mu_{n,1}&=&(\eta_n\circ z_{n+1})\cdot\nu_{n+1,1}\nonumber\\
(d_{n-1}\circ\eta_n)\cdot\nu_{n+1,2}&=&(\rho_n\circ x_n)\cdot(x_n)_I^{\flat}.\label{fifil}
\end{eqnarray}
Moreover, for any other $(x_n^{\prime},\xi_n^{\prime},\eta_n^{\prime})$ with these properties, there exists a unique 2-cell $\kappa\in\mathbf{S}^2(x_n,x_n^{\prime})$ such that $(\kappa\circ z_{n+1})\cdot\xi_n^{\prime}=\xi_n$ and $(q_n\circ\kappa)\cdot\eta_n^{\prime}=\eta_n$.
\end{lemma}
\begin{proof}
By the cofaithfulness of $q_n$, we can show $\mu_{n,2}$ is compatible with $\delta_{n+2}$. By the universality of the relative kernel $Z^{n+1}(A_{\bullet})$, there exist $x_n\in\mathbf{S}^1(Q^n(A_{\bullet}),Z^{n+1}(A_{\bullet}))$ and $\xi
_n\in\mathbf{S}^2(x_n\circ z_{n+1},\ell_n)$ such that
\[ (\xi_n\circ d_{n+1})\cdot\mu_{n,2}=(x_n\circ\zeta_{n+1})\cdot(x_n)_I^{\sharp}. \]
Then, both the factorizations
\begin{eqnarray*}
\nu_{n+1,1}&:&k_{n+1}\circ z_{n+1}\Longrightarrow d_n\\
(q_n\circ\xi_n)\cdot\mu_{n,1}&:&q_n\circ x_n\circ z_{n+1}\Longrightarrow d_n
\end{eqnarray*}
are compatible with $\zeta_{n+1}$ and $\delta_{n+1}$. Thus by the universality of relative kernel $Z^{n+1}(A_{\bullet})$, there exists a unique 2-cell $\eta_n\in\mathbf{S}^2(q_n\circ x_n,k_{n+1})$ such that $(q_n\circ\xi_n)\cdot\mu_{n,1}=(\eta_n\circ z_{n+1})\cdot\nu_{n+1,1}$. It can be easily seen that $\eta_n$ also satisfies $(\ref{fifil})$. Uniqueness (up to an equivalence) of $(x_n,\xi_n,\eta_n)$ follows from the universality of the relative kernel $Z^{n+1}(A_{\bullet})$ and the uniqueness of $\eta_n$.
\end{proof}

\begin{lemma}\label{FFFFLem1}
Consider the following complex diagram in $\mathbf{S}$.
\begin{equation}
\xy
(-10,0)*+{A}="4";
(10,0)*+{B}="6";
(30,0)*+{C}="8";
{\ar_{f} "4";"6"};
{\ar_{g} "6";"8"};
{\ar@/^1.65pc/^{0} "4";"8"};
{\ar@{=>}_{\varphi} (10,2)*{};(10,6)*{}};
\endxy
\label{prefindiag1}
\end{equation}
If $(\ref{prefindiag1})$ is 2-exact in $B$ and $g$ is cofaithful, then we have $\mathrm{Cok}(f)=[C,g,\varphi]$.
\end{lemma}
\begin{proof}
If we factor $(\ref{prefindiag1})$ as
\[
\xy
(-19,0)*+{A}="0";
(0,0)*+{B}="2";
(20,10)*+{\mathrm{Cok}(f)}="4";
(20,-10)*+{C}="6";
(24,-11)*+{,}="7";
{\ar^{f} "0";"2"};
{\ar|>>>>>>>*{_{c(f)}} "2";"4"};
{\ar|*+{_g} "2";"6"};
{\ar^{\overline{g}} "4";"6"};
{\ar@/^1.50pc/^{0} "0";"4"};
{\ar@/_1.50pc/_{0} "0";"6"};
{\ar@{=>}_{\pi_f} (0,4)*{};(-1,8)*{}};
{\ar@{=>}_{\varphi}(0,-4)*{};(-1,-8)*{}};
{\ar@{=>}^{\overline{\varphi}} (14,2)*{};(11,-2)*{}};
\endxy
\]
then, since $(\ref{prefindiag1})$ is 2-exact in $B$, $\overline{g}$ becomes fully faithful. On the other hand, since $g$ is cofaithful, $\overline{g}$ is also cofaithful. Thus $\overline{g}$ becomes an equivalence.
\end{proof}

\begin{lemma}\label{FFFFLem2}
Consider the following complex morphism in $\mathbf{S}$.
\[
\xy
(-16,6)*+{A_1}="0";
(-16,-6)*+{A_1}="2";
(0,6)*+{A_2}="4";
(0,-6)*+{B_2}="6";
(16,6)*+{A_3}="8";
(16,-6)*+{B_3}="10";
(32,6)*+{0}="12";
(32,-6)*+{0}="14";
{\ar^{d^A_1} "0";"4"};
{\ar@{=}_{\mathrm{id}} "0";"2"};
{\ar_{d^B_1} "2";"6"};
{\ar|*+{_{f_2}} "4";"6"};
{\ar^{f_3} "8";"10"};
{\ar^{d^A_2} "4";"8"};
{\ar_{d^B_2} "6";"10"};
{\ar^{0} "8";"12"};
{\ar_{0} "10";"14"};
{\ar@/^1.80pc/^{0} "0";"8"};
{\ar@/_1.80pc/_{0} "2";"10"};
{\ar@{=>}_{\lambda} (-8,2)*{};(-8,-2)*{}};
{\ar@{=>}_{\kappa} (8,2)*{};(8,-2)*{}};
{\ar@{=>}_{\delta^A_2} (0,9)*{};(0,12)*{}};
{\ar@{=>}_{\delta^B_2} (0,-9)*{};(0,-12)*{}};
\endxy
\]
If the complexes are relatively 2-exact in $A_2,A_3$ and $B_2,B_3$ respectively, i.e. they satisfy $\mathrm{Cok}(d_1^A)=[A_3,d_2^A,\delta_2^A]$ and $\mathrm{Cok}(d_1^B)=[B_3,d_2^B,\delta_2^B]$ $($see Corollary \ref{final29}$)$, then the following diagram obtained by taking the kernel of $f_2$ becomes 2-exact in $A_3$.
\begin{equation}
\xy
(-24,0)*+{\mathrm{Ker}(f_2)}="4";
(0,0)*+{A_3}="6";
(24,0)*+{B_3}="8";
{\ar_{k(f_2)\circ d_2^A} "4";"6"};
{\ar_{f_3} "6";"8"};
{\ar@/^1.80pc/^{0} "4";"8"};
{\ar@{=>} (0,3)*{};(0,6)*{}};
{\ar@{..}@/_1.20pc/ (0,4)*{};(26,-8)*{_{(k(f_2)\circ\kappa)\cdot(\varepsilon_{f_2} \circ d^B_2 )\cdot(d^B_2 )^{\flat}_I}}};
\endxy
\label{elledia2}
\end{equation}
\end{lemma}
\begin{proof}
By taking the kernels of $\mathrm{id}_{A_1}$ and $f_2$ in the diagram
\[
\xy
(-8,6)*+{A_1}="21";
(8,6)*+{A_2}="22";
(-8,-6)*+{A_1}="31";
(8,-6)*+{B_2}="32";
{\ar^{d_1^A} "21";"22"};
{\ar@{=}_{\mathrm{id}} "21";"31"};
{\ar^{f_2} "22";"32"};
{\ar_{d_1^B} "31";"32"};
{\ar@{=>}_{\lambda} (0,2)*{};(0,-2)*{}};
\endxy
\]
and taking the cokernels of $0_{0,A_1}$ and $k(f_2)$, we obtain the following diagram by Proposition \ref{final9}, where $\theta=k(f_2)_I^{\flat}\cdot(d_1^A)_I^{\flat -1}$:
\begin{equation}
\xy
(-36,6)*+{0}="11";
(-36,-6)*+{\mathrm{Ker}(f_2)}="12";
(-12,6)*+{A_1}="21";
(-12,-6)*+{A_2}="22";
(12,6)*+{A_1}="31";
(12,-6)*+{\mathrm{Coim}(f_2)}="32";
(36,6)*+{A_1}="41";
(36,-6)*+{B_2}="42";
{\ar_{0} "11";"12"};
{\ar^{0} "11";"21"};
{\ar_{k(f_2)} "12";"22"};
{\ar^{d_1^A} "21";"22"};
{\ar@{=}^{\mathrm{id}} "21";"31"};
{\ar_{c(k(f_2))} "22";"32"};
{\ar_{\overline{d}_1^A} "31";"32"};
{\ar@{=}^{\mathrm{id}} "31";"41"};
{\ar_{j(f_2)} "32";"42"};
{\ar^{d_1^B} "41";"42"};
{\ar@/_1.80pc/_{f_2} "22";"42"};
{\ar@{=>}_{\theta} (-26,0)*{};(-22,0)*{}};
{\ar@{=>}_{\lambda_1} (-2,0)*{};(2,0)*{}};
{\ar@{=>}_{\lambda_2} (22,0)*{};(26,0)*{}};
{\ar@{=>}_{\mu_{f_2}} (12,-9)*{};(12,-12)*{}};
\endxy
\label{elledia}
\end{equation}
By taking the cokernels of $0_{0,\mathrm{Ker}(f_2)}$, $d_1^A$ and $d_1^B$ in $(\ref{elledia})$, we obtain the left of the following diagrams, while by Proposition \ref{final15} we obtain the right as a coimage factorization if we take the cokernels of $d_1^A$, $\overline{d}_1^A$ and $d_1^B$ in $(\ref{elledia})$:
\[
\xy
(-24,6)*+{\mathrm{Ker}(f_2)}="11";
(-24,-6)*+{\mathrm{Ker}(f_2)}="12";
(0,6)*+{A_2}="21";
(0,-6)*+{A_3}="22";
(16,6)*+{B_2}="31";
(16,-6)*+{B_3}="32";
{\ar@{=}_{\mathrm{id}} "11";"12"};
{\ar^{k(f_2)} "11";"21"};
{\ar^{k(f_2)\circ d^A_2} "12";"22"};
{\ar_{d_2^A} "21";"22"};
{\ar^{f_2} "21";"31"};
{\ar_{f_3} "22";"32"};
{\ar^{d_2^B} "31";"32"};
{\ar@/^1.80pc/^>>>>>>>>>{0} "11";"31"};
{\ar@/_1.80pc/_>>>>>>>>>{0} "12";"32"};
{\ar@{}|\circlearrowright"11";"22"};
{\ar@{=>}^{\kappa} (6,0)*{};(10,0)*{}};
{\ar@{=>}^{\varepsilon_{f_2}} (0,9)*{};(0,12)*{}};
{\ar@{=>} (0,-9)*{};(0,-12)*{}} ;
{\ar@{.}@/_1.20pc/ (0,-11)*{};(-12,-20)*{_{(k(f_2)\circ\kappa)\cdot(\varepsilon_{f_2}\circ d^B_2)\cdot(d^B_2)^{\flat}_I}}};
\endxy
\quad\xy
(-24,6)*+{A_2}="21";
(-24,-6)*+{A_3}="22";
(0,6)*+{\mathrm{Coim}(f_2)}="31";
(0,-6)*+{\mathrm{Coim}(f_3)}="32";
(20,6)*+{B_2}="41";
(20,-6)*+{B_3}="42";
{\ar_{d_2^A} "21";"22"};
{\ar^{c(k(f_2))} "21";"31"};
{\ar_{c(k(f_3))} "22";"32"};
{\ar_{\overline{d}_2^A} "31";"32"};
{\ar^{j(f_2)} "31";"41"};
{\ar_{j(f_3)} "32";"42"};
{\ar^{d_2^B} "41";"42"};
{\ar@/^1.80pc/^>>>>>>>>>>>{f_2} "21";"41"};
{\ar@/_1.80pc/_>>>>>>>>>>>{f_3} "22";"42"};
{\ar@{=>}^{\kappa_1} (-16,0)*{};(-12,0)*{}};
{\ar@{=>}^{\kappa_2} (8,0)*{};(12,0)*{}};
{\ar@{=>}^{\mu_{f_2}} (0,9)*{};(0,12)*{}};
{\ar@{=>}^{\mu_{f_3}} (0,-9)*{};(0,-12)*{}};
\endxy
\]
On the other hand by Proposition \ref{final8}, if we take the compatible 2-cell $\upsilon=(k(f_2)\circ\kappa_1)\cdot(\pi_{k(f_2)}\circ\overline{d}_2^A)\cdot(\overline{d}_2^A)_I^{\flat}\in\mathbf{S}^2(k(f_2)\circ d_2^A\circ c(k(f_3)),0)$,\[
\xy
(-24,6)*+{\mathrm{Ker}(f_2)}="11";
(-24,-6)*+{\mathrm{Ker}(f_2)}="12";
(0,6)*+{A_2}="21";
(0,-6)*+{A_3}="22";
(24,6)*+{\mathrm{Coim}(f_2)}="31";
(24,-6)*+{\mathrm{Coim}(f_3)}="32";
{\ar@{=}_{\mathrm{id}} "11";"12"};
{\ar^{k(f_2)} "11";"21"};
{\ar_{k(f_2)\circ d^A_2} "12";"22"};
{\ar_{d_2^A} "21";"22"};
{\ar^{c(k(f_2))} "21";"31"};
{\ar_{c(k(f_3))} "22";"32"};
{\ar^{\overline{d}_2^A} "31";"32"};
{\ar@/^2.00pc/^{0} "11";"31"};
{\ar@/_2.00pc/_{0} "12";"32"};
{\ar@{}|\circlearrowright"11";"22"};
{\ar@{=>}_{\kappa_1} (10,0)*{};(14,0)*{}};
{\ar@{=>}_{\pi_{k(f_2)}} (0,9)*{};(0,12)*{}};
{\ar@{=>}_{\upsilon} (0,-9)*{};(0,-12)*{}};
\endxy
\]
then we have $\mathrm{Cok}(k(f_2)\circ d_2^A)=[\mathrm{Coim}(f_3),c(k(f_3)),\upsilon]$. It can be easily shown that $\upsilon$ is compatible with $\mu_{f_3}$ and $(k(f_2)\circ\kappa)\cdot(\varepsilon_{f_2}\circ d_2^B)\cdot(d_2^B)_I^{\flat}$.
\[
\xy
(-24,0)*+{\mathrm{Ker}(f_2)}="2";
(22,-10)*+{B_3}="4";
(0,0)*+{A_3}="6";
(22,10)*+{\mathrm{Coim}(f_3)}="8";
{\ar@/_1.40pc/_<<<<<<<<<<<{0} "2";"4"};
{\ar_{f_3} "6";"4"};
{\ar_{k(f_2)\circ d^A_2} "2";"6"};
{\ar^{_{c(k(f_3))}} "6";"8"};
{\ar^{j(f_3)} "8";"4"};
{\ar@/^1.40pc/^<<<<<<<<<<<{0} "2";"8"};
{\ar@{=>}^{\upsilon} (-1,3)*{};(-2,8)*{}};
{\ar@{=>} (-1,-3)*{};(-3,-8)*{}};
{\ar@{=>}^{\mu_{f_3}} (17,3)*{};(13,-2)*{}};
{\ar@{..}@/^1.60pc/^{} (-3,-7)*{};(-36,-8)*{_{(k(f_2) \circ\kappa)\cdot(\varepsilon_{f_2}\circ d^B_2)\cdot(d^B_2)^{\flat}_I}}};
\endxy
\]
Since $\mathrm{Cok}(k(f_2)\circ d_2^A)=[\mathrm{Coim}(f_3),c(k(f_3)),\upsilon]$ and $j(f_3)$ is fully faithful by Proposition \ref{SecondFactorization}, this means $(\ref{elledia2})$ is 2-exact in $A_3$.
\end{proof}

\begin{lemma}\label{FFFFLem4}
Consider the following complex morphism in $\mathbf{S}$.
\begin{equation}
\xy
(-16,6)*+{A_1}="0";
(-16,-6)*+{B_1}="2";
(0,6)*+{A_2}="4";
(0,-6)*+{B_2}="6";
(16,6)*+{A_3}="8";
(16,-6)*+{B_3}="10";
(-32,-6)*+{0}="12";
{\ar^{d^A_1} "0";"4"};
{\ar_{f_1} "0";"2"};
{\ar_{d^B_1} "2";"6"};
{\ar|*+{_{f_2}} "4";"6"};
{\ar^{f_3} "8";"10"};
{\ar^{d^A_2} "4";"8"};
{\ar_{d^B_2} "6";"10"};
{\ar^{0} "12";"2"};
{\ar@/^1.65pc/^{0} "0";"8"};
{\ar@/_1.65pc/_{0} "2";"10"};
{\ar@{=>}_{\lambda_1} (-8,2)*{};(-8,-2)*{}};
{\ar@{=>}_{\lambda_2} (8,2)*{};(8,-2)*{}};
{\ar@{=>}_{\delta^A _2} (0,9)*{};(0,12)*{}};
{\ar@{=>}_{\delta^B _2} (0,-9)*{};(0,-12)*{}};
\endxy
\label{prefindiag2}
\end{equation}
If the complexes are relatively 2-exact in $A_2$ and $B_1,B_2$ respectively, then the following diagram obtained by taking the kernels
\begin{equation}
\xy
(-24,0)*+{\mathrm{Ker}(f_1)}="4";
(0,0)*+{\mathrm{Ker}(f_2)}="6";
(24,0)*+{\mathrm{Ker}(f_3)}="8";
{\ar_{\underline{d}^A_1} "4";"6"};
{\ar_{\underline{d}^A_2} "6";"8"};
{\ar@/^1.80pc/^{0} "4";"8"};
{\ar@{=>}_{\underline{\delta}^A_2} (0,3)*{};(0,6)*{}};
\endxy
\label{prefindiag3}
\end{equation}
is 2-exact in $\mathrm{Ker}(f_2)$.
\end{lemma}
\begin{proof}
If we decompose $(\ref{prefindiag2})$ into
\[
\xy
(-16,6)*+{A_1}="0";
(-16,-6)*+{B_1}="2";
(4,6)*+{\mathrm{Ker}(d^A_2)}="4";
(4,-6)*+{\mathrm{Ker}(d^B_2)}="6";
(-32,-6)*+{0}="12";
{\ar^{d^{A\dag}_1} "0";"4"};
{\ar_{f_1} "0";"2"};
{\ar_{d^{B\dag}_1} "2";"6"};
{\ar^{\underline{f}_2} "4";"6"};
{\ar^{0} "12";"2"};
{\ar@{=>}_{\lambda^{\dag}_1} (-6,2)*{};(-6,-2)*{}};
\endxy
\quad
\text{and}
\quad
\xy
(-20,6)*+{\mathrm{Ker}(d_2^A)}="0";
(-20,-6)*+{\mathrm{Ker}(d_2^B)}="2";
(0,6)*+{A_2}="4";
(0,-6)*+{B_2}="6";
(16,6)*+{A_3}="8";
(16,-6)*+{B_3}="10";
(18,-7)*+{,}="12";
{\ar^{k(d^A_2)} "0";"4"};
{\ar_{\underline{f}_2} "0";"2"};
{\ar_{k(d^B_2)} "2";"6"};
{\ar|*+{_{f_2}} "4";"6"};
{\ar^{f_3} "8";"10"};
{\ar^{d^A_2} "4";"8"};
{\ar_{d^B_2} "6";"10"};
{\ar@{=>}_{\underline{\lambda}_1} (-10,2)*{};(-10,-2)*{}};
{\ar@{=>}_{\lambda_2} (8,2)*{};(8,-2)*{}};
\endxy
\]
then by (the dual of) Proposition \ref{final8}, we have $\mathrm{Ker}(\underline{d}_2^A)=\mathrm{Ker}(\underline{f}_2)$. Since $d_1^{B\dag}$ is an equivalence by (the dual of) Corollary \ref{final29}, the diagram obtained by taking the kernels of $f_1$ and $\underline{f}_2$
\[
\xy
(-12,6)*+{\mathrm{Ker}(f_1)}="0";
(-12,-6)*+{A_1}="2";
(12,6)*+{\mathrm{Ker}(\underline{f}_2)}="4";
(12,-6)*+{\mathrm{Ker}(d^A_2)}="6";
{\ar^{\underline{d}^{A\dag}_1} "0";"4"};
{\ar_{k(f_1)} "0";"2"};
{\ar_{d^{A\dag}_1} "2";"6"};
{\ar^{k(\underline{f}_2)} "4";"6"};
{\ar@{=>}_{\underline{\lambda}^{\dag}_1} (0,2)*{};(0,-2)*{}};
\endxy
\]
becomes a pullback diagram by (the dual of) Proposition \ref{final12}. Since $d_1^{A\dag}$ is fully cofaithful, $\underline{d}_1^{A\dag}$ becomes also fully cofaithful by Proposition \ref{final13}. This means $(\ref{prefindiag3})$ is 2-exact in $\mathrm{Ker}(f_2)$.
\end{proof}

\begin{lemma}\label{FFFFLem5}

Consider the following complex morphism in $\mathbf{S}$.
\begin{equation}
\xy
(-16,6)*+{A_1}="0";
(-16,-6)*+{A_1}="2";
(0,6)*+{A_2}="4";
(0,-6)*+{A_2}="6";
(16,6)*+{A_3}="8";
(16,-6)*+{B_3}="10";
{\ar^{d^A_1} "0";"4"};
{\ar@{=}_{\mathrm{id}} "0";"2"};
{\ar_{d^B_1} "2";"6"};
{\ar@{=}_{\mathrm{id}} "4";"6"};
{\ar^{f_3} "8";"10"};
{\ar^{d^A_2} "4";"8"};
{\ar_{d^B_2} "6";"10"};
{\ar@/^1.65pc/^{0} "0";"8"};
{\ar@/_1.65pc/_{0} "2";"10"};
{\ar@{=>}_{\lambda_1} (-8,2)*{};(-8,-2)*{}};
{\ar@{=>}_{\lambda_2} (8,2)*{};(8,-2)*{}};
{\ar@{=>}_{\delta^A_2} (0,9)*{};(0,12)*{}};
{\ar@{=>}_{\delta^B_2} (0,-9)*{};(0,-12)*{}};
\endxy
\label{prefindiag4}
\end{equation}
If $f_3$ is faithful and the bottom row is 2-exact in $A_2$, then the top row is also 2-exact in $A_2$.
\end{lemma}

\begin{proof}
By taking the cokernels of $d_1^A$ and $d_1^B$ in $(\ref{prefindiag4})$, we obtain (by Proposition \ref{final9})
\[
\xy
(-10,6)*+{\mathrm{Cok}(d^A_1)}="22";
(10,6)*+{A_3}="23";
(-10,-6)*+{\mathrm{Cok}(d^B_1)}="32";
(10,-6)*+{B_3}="33";
(13,-7)*+{.}="34";
{\ar^{\overline{d}_2^A} "22";"23"};
{\ar@{=}_{\mathrm{id}} "22";"32"};
{\ar^{f_3} "23";"33"};
{\ar_{\overline{d}_2^B} "32";"33"};
{\ar@{=>}_{\overline{\lambda}_2} (0,2)*{};(0,-2)*{}};
\endxy
\]
Since $\overline{d}_2^B$ is fully faithful, by taking the kernels in this diagram, we obtain the following diagram.
\[
\xy
(0,12)*+{0}="12";
(24,12)*+{\mathrm{Ker}(f_3)}="13";
(-24,0)*+{\mathrm{Ker}(\overline{d}^A_2)}="21";
(0,0)*+{\mathrm{Cok}(d^A_1)}="22";
(24,0)*+{A_3}="23";
(-24,-12)*+{0}="31";
(0,-12)*+{\mathrm{Cok}(d^B_1)}="32";
(24,-12)*+{B_3}="33";
{\ar^{0} "12";"13"};
{\ar_{0} "12";"22"};
{\ar^{k(f_3)} "13";"23"};
{\ar^{k(\overline{d}_2^A)} "21";"22"};
{\ar|*+{_{\overline{d}_2^A}} "22";"23"};
{\ar_{0} "21";"31"};
{\ar@{=}_{\mathrm{id}} "22";"32"};
{\ar^{f_3} "23";"33"};
{\ar_{0} "31";"32"};
{\ar_{\overline{d}_2^B} "32";"33"};
{\ar@{=>}_{\exists} (12,8)*{};(12,4)*{}};
{\ar@{=>}_{\exists} (-12,-4)*{};(-12,-8)*{}};
{\ar@{=>}_{\overline{\lambda}_2} (12,-4)*{};(12,-8)*{}};
\endxy
\]
In this diagram, we have
\[ \mathrm{Ker}(\overline{d}_2^A)=\;\mathrm{Ker}(\mathrm{Ker}(\overline{d}_2^A)\longrightarrow 0)\underset{\text{Prop. \ref{final8}}}{=}\mathrm{Ker}(0\longrightarrow \mathrm{Ker}(f_{3}))\underset{\text{Cor. \ref{final4}}}{=}0. \]
This means that the top row in $(\ref{prefindiag4})$ is 2-exact in $A_2$.
\end{proof}


\begin{corollary}\label{FinFinFinCor}

Let
\[
\xy
(-16,0)*+{A_1}="2";
(0,0)*+{A_2}="6";
(16,0)*+{A_3}="10";
{\ar_{d^A_1} "2";"6"};
{\ar_{d^A_2} "6";"10"};
{\ar@/^1.65pc/^{0} "2";"10"};
{\ar@{=>}^{\delta^A_2} (0,2)*{};(0,6)*{}};
\endxy
\quad
\text{and}
\quad
\xy
(-16,0)*+{B_1}="2";
(0,0)*+{A_2}="6";
(16,0)*+{B_3}="10";
{\ar_{d^B_1} "2";"6"};
{\ar_{d^B_2} "6";"10"};
{\ar@/^1.65pc/^{0} "2";"10"};
{\ar@{=>}^{\delta^B_2} (0,2)*{};(0,6)*{}};
\endxy
\]
be two complexes, and assume that there exist 1-cells $f_1,f_3$ and 2-cells $\lambda_1,\lambda_2,\sigma$ as in the following diagram
\[
\xy
(-24,8)*+{B_1}="0";
(-24,-8)*+{A_1}="2";
(0,0)*+{A_2}="4";
(24,8)*+{A_3}="6";
(24,-8)*+{B_3}="8";
(27,-9)*+{,}="10";
{\ar_{f_1} "0";"2"};
{\ar^{d^B_1} "0";"4"};
{\ar_{d^A_1} "2";"4"};
{\ar_{d^B_2} "4";"8"};
{\ar^{d^A_2} "4";"6"};
{\ar^{f_3} "6";"8"};
{\ar@/_1.20pc/_{0} "2";"8"};
{\ar@{=>}_{\sigma} (0,-6)*{};(0,-10)*{}};
{\ar@{=>}^{\lambda_1} (-18,-3)*{};(-14,1)*{}};
{\ar@{=>}^{\lambda_2} (18,3)*{};(14,-1)*{}};
\endxy
\]
where $f_1$ is cofaithful and $f_3$ is faithful. Assume they satisfy
\begin{eqnarray*}
&\mathrm{(d1)}&\ (\lambda_1\circ d_2^B)\cdot\delta_2^B=(f_1\circ\sigma)\cdot(f_1)_I^{\sharp}\\
&\mathrm{(d2)}&\ (d_1^A\circ\lambda_2)\cdot\sigma=(\delta_2^A\circ f_3)\cdot(f_3)_I^{\flat}.
\end{eqnarray*}
Then, if the diagram
\[
\xy
(-16,0)*+{B_1}="2";
(0,0)*+{A_2}="6";
(16,0)*+{B_3}="10";
{\ar_{d^B_1} "2";"6"};
{\ar_{d^B_2} "6";"10"};
{\ar@/^1.65pc/^{0} "2";"10"};
{\ar@{=>}^{\delta^B_2} (0,2)*{};(0,6)*{}};
\endxy
\]
is 2-exact in $A_2$, then the diagram
\[
\xy
(-16,-8)*+{A_1}="2";
(0,-8)*+{A_2}="6";
(16,-8)*+{A_3}="10";
{\ar_{d^A_1} "2";"6"};
{\ar_{d^A_2} "6";"10"};
{\ar@/^1.65pc/^{0} "2";"10"};
{\ar@{=>}^{\delta^A_2} (0,-6)*{};(0,-2)*{}};
\endxy
\]
is also 2-exact in $A_2$.
\end{corollary}
\begin{proof}
This follows if we apply Lemma \ref{FFFFLem5} and its dual to the following diagrams:
\[
\xy
(-16,6)*+{B_1}="0";
(-16,-6)*+{A_1}="2";
(0,6)*+{A_2}="4";
(0,-6)*+{A_2}="6";
(16,6)*+{B_3}="8";
(16,-6)*+{B_3}="10";
(19,-7)*+{,}="12";
{\ar^{d^B_1} "0";"4"};
{\ar_{f_1} "0";"2"};
{\ar_{d^A_1} "2";"6"};
{\ar@{=}^{\mathrm{id}} "4";"6"};
{\ar@{=}^{\mathrm{id}} "8";"10"};
{\ar^{d^B_2} "4";"8"};
{\ar_{d^B_2} "6";"10"};
{\ar@{}|\circlearrowright"4";"10"};
{\ar@/^1.65pc/^{0} "0";"8"};
{\ar@/_1.65pc/_{0} "2";"10"};
{\ar@{=>}_{\lambda^{-1}_1} (-8,2)*{};(-8,-2)*{}};
{\ar@{=>}_{\delta^B_2} (0,8)*{};(0,12)*{}};
{\ar@{=>}_{\sigma} (0,-8)*{};(0,-12)*{}};
\endxy
\quad
\xy
(-16,6)*+{A_1}="0";
(-16,-6)*+{A_1}="2";
(0,6)*+{A_2}="4";
(0,-6)*+{A_2}="6";
(16,6)*+{A_3}="8";
(16,-6)*+{B_3}="10";
{\ar^{d^A_1} "0";"4"};
{\ar@{=}_{\mathrm{id}} "0";"2"};
{\ar_{d^A_1} "2";"6"};
{\ar@{=}_{\mathrm{id}} "4";"6"};
{\ar^{f_3} "8";"10"};
{\ar^{d^A_2} "4";"8"};
{\ar_{d^B_2} "6";"10"};
{\ar@{}|\circlearrowright"0";"6"};
{\ar@/^1.65pc/^{0} "0";"8"};
{\ar@/_1.65pc/_{0} "2";"10"};
{\ar@{=>}_{\lambda_2} (8,2)*{};(8,-2)*{}};
{\ar@{=>}_{\delta^A_2} (0,8)*{};(0,12)*{}};
{\ar@{=>}_{\sigma} (0,-8)*{};(0,-12)*{}};
\endxy
\]
\end{proof}

By Corollary \ref{FinFinFinCor}, it can be shown that the 2-exactness plus compatibility implies the relative 2-exactness (see \cite{RMV} in the case of $\mathrm{SCG}$):

\begin{corollary}
Let $A_{\bullet}=(A_n,d_n,\delta_n)$ be a complex in $\mathbf{S}$. If
\[
\xy
(-16,0)*+{A_{n-1}}="2";
(0,0)*+{A_n}="6";
(16,0)*+{A_{n+1}}="10";
{\ar_{d^A_{n-1}} "2";"6"};
{\ar_{d^A_n} "6";"10"};
{\ar@/^1.65pc/^{0} "2";"10"};
{\ar@{=>}^{\delta^A _n } (0,2)*{};(0,6)*{}};
\endxy
\]
is 2-exact in $A_n$, then $A_{\bullet}$ is relatively 2-exact in $A_n$.
\end{corollary}

\begin{proof}
This follows immediately if we apply Corollary \ref{FinFinFinCor} to the following diagram (see the proof of Corollary \ref{final26}):
\[
\xy
(-26,8)*+{A_{n-1}}="0";
(-26,-8)*+{\mathrm{Cok}(d_{n-2})}="2";
(0,0)*+{A_n}="4";
(26,8)*+{\mathrm{Ker}(d_{n+1})}="6";
(26,-8)*+{A_{n+1}}="8";
{\ar_{c(d_{n-2})} "0";"2"};
{\ar^{d_{n-1}} "0";"4"};
{\ar_{\overline{d}_{n-1}} "2";"4"};
{\ar_{d_n} "4";"8"};
{\ar^{\underline{d}_n} "4";"6"};
{\ar^{k(d_{n+1})} "6";"8"};
{\ar@/_1.20pc/_{0} "2";"8"};
{\ar@{=>}_{\overline{\delta}_n} (0,-6)*{};(0,-10)*{}};
{\ar@{=>}^{\overline{\delta}_{n-1}} (-18,-3)*{};(-14,1)*{}};
{\ar@{=>}^{\underline{\delta}_{n+1}} (18,3)*{};(14,-1)*{}};
\endxy
\]
\end{proof}

\paragraph{Construction of the long cohomology sequence}\quad\\

\begin{definition}
A complex in $\mathbf{S}$
\begin{equation}
\xy
(-16,0)*+{A}="4";
(0,0)*+{B}="6";
(16,0)*+{C}="8";
{\ar_{f} "4";"6"};
{\ar_{g} "6";"8"};
{\ar@/^1.65pc/^{0} "4";"8"};
{\ar@{=>}_{\varphi} (0,2)*{};(0,6)*{}};
\endxy
\label{ExtensionDiagram}
\end{equation}
is called an extension if it is relatively 2-exact in every 0-cell.
\end{definition}

\begin{remark}
By Corollary \ref{final29} (and its dual), $(\ref{ExtensionDiagram})$ is an extension if and only if $\mathrm{Ker}(g)=[A,f,\varphi]$ and $\mathrm{Cok}(f)=[C,g,\varphi]$.
\end{remark}

\begin{definition}
Let $(f_{\bullet},\lambda_{\bullet}):A_{\bullet}\longrightarrow B_{\bullet}$ and $(g_{\bullet},\kappa_{\bullet}):B_{\bullet}\longrightarrow C_{\bullet}$ be complex morphisms and $\varphi_{\bullet}=\{\varphi_n:f_n\circ g_n\Longrightarrow 0\}$ be 2-cells. Then,
\begin{equation}
\xy
(-16,0)*+{A_{\bullet}}="4";
(0,0)*+{B_{\bullet}}="6";
(16,0)*+{C_{\bullet}}="8";
{\ar_{f_{\bullet}} "4";"6"};
{\ar_{g_{\bullet}} "6";"8"};
{\ar@/^1.65pc/^{0} "4";"8"};
{\ar@{=>}_{\varphi_{\bullet}} (0,2)*{};(0,6)*{}};
\endxy
\label{ExtOfCpx}
\end{equation}
is said to be an extension of complexes if it satisfies the following properties$:$\\
{\rm (e1)} For every $n$, the following complex is an extension$:$
\[
\xy
(-8,0)*+{A_n}="4";
(8,0)*+{B_n}="6";
(24,0)*+{C_n}="8";
{\ar_{f_n} "4";"6"};
{\ar_{g_n} "6";"8"};
{\ar@/^1.65pc/^{0} "4";"8"};
{\ar@{=>}_{\varphi_n} (8,2)*{};(8,6)*{}} ;
\endxy
\]
\noindent{\rm (e2)} $\varphi_{\bullet}$ satisfies
\begin{eqnarray*}
&(\lambda_n\circ g_{n+1})\cdot(f_n\circ\kappa_n)\cdot(\varphi_n\circ d_n^C)\cdot(d_n^C)_I^{\flat}=(d_n^A\circ\varphi_{n+1})\cdot(d_n^A)_I^{\sharp}.&\\
&\xy
(-16,7)*+{A_n}="0";
(-16,-7)*+{A_{n+1}}="2";
(0,7)*+{B_n}="4";
(0,-7)*+{B_{n+1}}="6";
(16,7)*+{C_n}="8";
(16,-7)*+{C_{n+1}}="10";
{\ar^{f_n} "0";"4"};
{\ar_{d^A_n} "0";"2"};
{\ar^{f_{n+1}} "2";"6"};
{\ar|*+{_{d^B_n}} "4";"6"};
{\ar^{d^C_n} "8";"10"};
{\ar^{g_n} "4";"8"};
{\ar^{g_{n+1}} "6";"10"};
{\ar@/^1.65pc/^{0} "0";"8"};
{\ar@/_1.65pc/_{0} "2";"10"};
{\ar@{=>}^{\lambda_n} (-11,0)*{};(-5,0)*{}};
{\ar@{=>}^{\kappa_n} (5,0)*{};(11,0)*{}};
{\ar@{=>}_{\varphi_n} (0,9)*{};(0,13)*{}};
{\ar@{=>}_{\varphi_{n+1}} (0,-9)*{};(0,-13)*{}};
\endxy&
\end{eqnarray*}
\end{definition}

Our main theorem is the following:
\begin{theorem}\label{MainTheorem}
For any extension of complexes in $\mathbf{S}$
\[
\xy
(-8,0)*+{A_{\bullet}}="4";
(8,0)*+{B_{\bullet}}="6";
(24,0)*+{C_{\bullet}}="8";
{\ar_{f_{\bullet}} "4";"6"};
{\ar_{g_{\bullet}} "6";"8"};
{\ar@/^1.65pc/^{0} "4";"8"};
{\ar@{=>}_{\varphi_{\bullet}} (8,2)*{};(8,6)*{}};
\endxy
,
\]
we can construct a long 2-exact sequence$:$
\[
\xy
(-48,0)*+{\cdots}="0";
(-30,0)*+{H^n(B_{\bullet})}="2";
(-10,0)*+{H^n(C_{\bullet})}="3";
(10,0)*+{H^{n+1}(A_{\bullet})}="4";
(36,0)*+{H^{n+1}(B_{\bullet})}="5";
(54,0)*+{\cdots}="7";
{\ar "0";"2"};
{\ar "2";"3"};
{\ar "3";"4"};
{\ar "4";"5"};
{\ar "5";"7"};
{\ar@/^1.65pc/^{0} "0";"3"};
{\ar@/_1.65pc/_{0} "2";"4"};
{\ar@/^1.65pc/^{0} "3";"5"};
{\ar@/_1.65pc/_{0} "4";"7"};
{\ar@{=>} (-30,2)*{};(-30,6)*{}};
{\ar@{=>} (-10,-2)*{};(-10,-6)*{}};
{\ar@{=>} (14,2)*{};(14,6)*{}};
{\ar@{=>} (32,-2)*{};(32,-6)*{}};
\endxy
\]
\end{theorem}

\begin{caution}
This sequence is not necessarily a complex. $($See Remark \ref{B4Remark}.$)$
\end{caution}

We prove this theorem in the rest of this section.
\begin{lemma}\label{FINLem2}
In the notation of Lemma \ref{FINLem1}, we have\\
{\rm (1)} $\mathrm{Ker}(x_n)=H^n(A_{\bullet})$,\\
{\rm (2)} $\mathrm{Cok}(x_n)=H^{n+1}(A_{\bullet})$.
\end{lemma}

\begin{proof}
We only show {\rm (1)}, since {\rm (2)} can be shown in the same way. In the notation of Lemma \ref{FINLem1} and Remark \ref{final19}, we can show that the factorization
\[ (x_n\circ\underline{\zeta}_{n+1})\cdot\xi_n:(x_n\circ\underline{z}_{n+1})\circ k(d_{n+1})\Longrightarrow\ell_n \]
is compatible with $\varepsilon_{d_{n+1}}$ and $\mu_{n,2}$.
\[
\xy
(-12,12)*+{Q^n(A_{\bullet})}="10";
(-12,-12)*+{\mathrm{Ker}(d_{n+1})}="0";
(12,0)*+{A_{n+1}}="6";
(36,0)*+{A_{n+2}}="8";
{\ar^{d_{n+1}} "6";"8"};
{\ar_{k(d_{n+1})} "0";"6"};
{\ar_{x_n \circ\underline{z}_{n+1}} "10";"0"};
{\ar@{=>}^{\varepsilon_{d_{n+1}}} (12,-3)*{};(13,-8)*{}};
{\ar@/_1.20pc/_{0} "0";"8"};
{\ar@/^1.20pc/^{0} "10";"8"};
{\ar^{\ell_n} "10";"6"};
{\ar@{=>}^{\mu_{n,2}} (12,3)*{};(13,8)*{}};
{\ar@{=>} (-8,-4)*{};(-4,4)*{}};
{\ar@{..}@/_1.40pc/^{} (-6,0)*{};(-24,-20)*{_{(x_n\circ\underline{\zeta}_{n+1})\cdot\xi_n}}};
\endxy
\]
So, by Proposition \ref{KerAndRelKer}, Proposition \ref{KerFullFaith} and the fact that $\underline{z}_{n+1}$ is fully faithful, we have $H^n(A_{\bullet})=\mathrm{Ker}(\ell_n,\mu_{n,2})=\mathrm{Ker}(x_n\circ\underline{z}_{n+1})=\mathrm{Ker}(x_n).$
\end{proof}

\begin{lemma}\label{FINLem3}
For any extension $(\ref{ExtOfCpx})$ of complexes in $\mathbf{S}$, we can construct a complex morphism
\[
\xy
(-28,8)*+{Q^n(A_{\bullet})}="0";
(-28,-8)*+{Z^{n+1}(A_{\bullet})}="2";
(0,8)*+{Q^n(B_{\bullet})}="4";
(0,-8)*+{Z^{n+1}(B_{\bullet})}="6";
(28,8)*+{Q^n(C_{\bullet})}="8";
(48,8)*+{0}="18";
(28,-8)*+{Z^{n+1}(C_{\bullet})}="10";
(-48,-8)*+{0}="12";
{\ar_{Q^n(f_{\bullet})} "0";"4"};
{\ar_{x^A_n} "0";"2"};
{\ar^{0} "8";"18"};
{\ar_{0} "12";"2"};
{\ar^{Z^{n+1}(f_{\bullet})} "2";"6"};
{\ar_{x^B_n} "4";"6"};
{\ar^{x^C_n} "8";"10"};
{\ar_{Q^n(g_{\bullet})} "4";"8"};
{\ar^{Z^{n+1}(g_{\bullet})} "6";"10"};
{\ar@/^2.00pc/^{0} "0";"8"};
{\ar@/_2.00pc/_{0} "2";"10"};
{\ar@{=>}_{\widetilde{\lambda}_n} (-12,2)*{};(-12,-2)*{}};
{\ar@{=>}^{\widetilde{\kappa}_n} (12,2)*{};(12,-2)*{}};
{\ar@{=>}_{Q^n(\varphi_{\bullet})} (0,11)*{};(0,14)*{}};
{\ar@{=>}^{Z^{n+1}(\varphi_{\bullet})} (0,-11)*{};(0,-14)*{}};
\endxy
\]
where the top line is a complex which is relatively 2-exact in $Q^n(B_{\bullet})$, $Q^n(C_{\bullet})$, and the bottom line is a complex which is relatively 2-exact in $Z^{n+1}(A_{\bullet})$, $Z^{n+1}(B_{\bullet})$.
\end{lemma}
\begin{proof}
If we take the relative cokernels $Q^n(A_{\bullet})$, $Q^n(B_{\bullet})$ and $Q^n(C_{\bullet})$ of the complex diagram
\[
\xy
(-36,12)*+{A_{n-2}}="11";
(-12,12)*+{A_{n-1}}="12";
(12,12)*+{A_n}="13";
(36,12)*+{A_{n+1}}="14";
(-36,0)*+{B_{n-2}}="21";
(-12,0)*+{B_{n-1}}="22";
(12,0)*+{B_n}="23";
(36,0)*+{B_{n+1}}="24";
(-36,-12)*+{C_{n-2}}="31";
(-12,-12)*+{C_{n-1}}="32";
(12,-12)*+{C_n}="33";
(36,-12)*+{C_{n+1}}="34";
(39,-13)*+{,}="28";
{\ar^{d_{n-2}^A} "11";"12"};
{\ar^{d_{n-1}^A} "12";"13"};
{\ar^{d_n^A} "13";"14"};
{\ar_{f_{n-2}} "11";"21"};
{\ar^{f_{n-1}} "12";"22"};
{\ar^{f_n} "13";"23"};
{\ar^{f_{n+1}} "14";"24"};
{\ar_{d_{n-2}^B} "21";"22"};
{\ar_{d_{n-1}^B} "22";"23"};
{\ar_{d_n^B} "23";"24"};
{\ar_{g_{n-2}} "21";"31"};
{\ar^{g_{n-1}} "22";"32"};
{\ar^{g_n} "23";"33"};
{\ar^{g_{n+1}} "24";"34"};
{\ar_{d_{n-2}^C} "31";"32"};
{\ar_{d_{n-1}^C} "32";"33"};
{\ar_{d_n^C} "33";"34"};
{\ar@/^1.65pc/^{0} "11";"13"};
{\ar@/_1.65pc/_{0} "31";"33"};
{\ar@{=>}_{\lambda_{n-2}} (-24,8)*{};(-24,4)*{}};
{\ar@{=>}^{\lambda_{n-1}} (0,8)*{};(0,4)*{}};
{\ar@{=>}^{\lambda_n} (24,8)*{};(24,4)*{}};
{\ar@{=>}_{\delta_{n-1}^A} (-12,14)*{};(-12,18)*{}};
{\ar@{=>}_{\delta_{n-1}^C} (-12,-14)*{};(-12,-18)*{}};
{\ar@{=>}_{\kappa_{n-2}} (-24,-5)*{};(-24,-8)*{}};
{\ar@{=>}^{\kappa_{n-1}} (0,-5)*{};(0,-8)*{}};
{\ar@{=>}^{\kappa_n} (24,-5)*{};(24,-8)*{}};
\endxy
\]
then by (the dual of) Proposition \ref{KerAndRelKer}, Proposition \ref{final5} and Proposition \ref{final6}, we obtain a factorization diagram
\[
\xy
(-36,14)*+{A_{n-1}}="11";
(-12,14)*+{A_n}="12";
(12,14)*+{Q^n(A_{\bullet})}="13";
(36,14)*+{A_{n+1}}="14";
(-36,-0)*+{B_{n-1}}="21";
(-12,-0)*+{B_n}="22";
(12,-0)*+{Q^n(B_{\bullet})}="23";
(36,-0)*+{B_{n+1}}="24";
(-36,-14)*+{C_{n-1}}="31";
(-12,-14)*+{C_n}="32";
(12,-14)*+{Q^n(C_{\bullet})}="33";
(36,-14)*+{C_{n+1}}="34";
{\ar^{d_{n-1}^A} "11";"12"};
{\ar^{q_n^A} "12";"13"};
{\ar^{\ell_n^A} "13";"14"};
{\ar_{f_{n-1}} "11";"21"};
{\ar_{f_n} "12";"22"};
{\ar|*+{_{Q^n(f_{\bullet})}} "13";"23"};
{\ar^{f_{n+1}} "14";"24"};
{\ar_{d_{n-1}^B} "21";"22"};
{\ar_{q_n^B} "22";"23"};
{\ar_{\ell_n^B} "23";"24"};
{\ar_{g_{n-1}} "21";"31"};
{\ar_{g_n} "22";"32"};
{\ar|*+{_{Q^n(g_{\bullet})}} "23";"33"};
{\ar^{g_{n+1}} "24";"34"};
{\ar_{d_{n-1}^C} "31";"32"};
{\ar_{q_n^C} "32";"33"};
{\ar_{\ell_n^C} "33";"34"};
{\ar@{=>}_{\lambda_{n-1}} (-24,9)*{};(-24,5)*{}};
{\ar@{=>}_{\overline{\lambda}_{n-1,1} } (0,9)*{};(0,5)*{}};
{\ar@{=>}^{\overline{\lambda}_{n-1,2}} (24,9)*{};(24,5)*{}};
{\ar@{=>}_{\kappa_{n-1}} (-24,-5)*{};(-24,-9)*{}};
{\ar@{=>}_{\overline{\kappa}_{n-1,1}} (0,-5)*{};(0,-9)*{}};
{\ar@{=>}^{\overline{\kappa}_{n-1,2}} (24,-5)*{};(24,-9)*{}};
\endxy
\]
and a 2-cell $Q^n(\varphi_{\bullet})\in\mathbf{S}^2(Q^n(f_{\bullet})\circ Q^n(g_{\bullet}),0)$, which satisfy compatibility conditions in Proposition \ref{final5} and Proposition \ref{final6}. It is also easy to see by the universality of the relative cokernels that
\[
(\ell_n^A\circ\lambda_{n+1})\cdot(\overline{\lambda}_{n-1,2}\circ d_{n+1}^B)\cdot(Q^n(f_{\bullet})\circ\mu_{n,2}^B)\cdot(Q^n(f_{\bullet}))_I^{\sharp}=(\mu_{n,2}^A\circ f_{n+2})\cdot(f_{n+2})_I^{\flat}.
\]
Now, since
\[
\xy
(-24,0)*+{A_n}="4";
(-8,0)*+{B_n}="6";
(8,0)*+{C_n}="8";
(24,0)*+{0}="10";
{\ar_{f_n} "4";"6"};
{\ar_{g_n} "6";"8"};
{\ar_{0} "8";"10"};
{\ar@/^1.65pc/^{0} "4";"8"};
{\ar@{=>}_{\varphi_n} (-8,2)*{};(-8,6)*{}};
\endxy
\]
is relatively 2-exact in $B_n$ and $C_n$, we have $\mathrm{Cok}(f_n)=[C_n,g_n,\varphi_n]$. So, from $\mathrm{Cok}(f_n)=[C_n,g_n,\varphi_n]$ and $\mathrm{Cok}(f_{n-1})=[C_{n-1},g_{n-1},\varphi_{n-1}]$, by Proposition \ref{FINell} we obtain
\[ \mathrm{Cok}(Q^n(f_{\bullet}))=[Q^n(C_{\bullet}),Q^n(g_{\bullet}),Q^n(\varphi_{\bullet})], \]
i.e. the complex
\[
\xy
(-36,0)*+{Q^n(A_{\bullet})}="0";
(-12,0)*+{Q^n(B_{\bullet})}="4";
(12,0)*+{Q^n(C_{\bullet})}="8";
(36,0)*+{0}="18";
{\ar_{Q^n(f_{\bullet})} "0";"4"};
{\ar^{0} "8";"18"};
{\ar_{Q^n(g_{\bullet})} "4";"8"};
{\ar@/^1.80pc/^{0} "0";"8"};
{\ar@{=>}_{Q^n(\varphi_{\bullet})} (-12,2)*{};(-12,6)*{}};
\endxy
\]
is relatively 2-exact in $Q^n(B_{\bullet})$, $Q^n(C_{\bullet})$. Dually, we obtain a factorization diagram
\[
\xy
(-36,14)*+{A_n}="11";
(-12,14)*+{Z^{n+1}(A_{\bullet})}="12";
(12,14)*+{A_{n+1}}="13";
(36,14)*+{A_{n+2}}="14";
(-36,0)*+{B_n}="21";
(-12,0)*+{Z^{n+1}(B_{\bullet})}="22";
(12,0)*+{B_{n+1}}="23";
(36,0)*+{B_{n+2}}="24";
(-36,-14)*+{C_n}="31";
(-12,-14)*+{Z^{n+1}(C_{\bullet})}="32";
(12,-14)*+{C_{n+1}}="33";
(36,-14)*+{C_{n+2}}="34";
{\ar^{k_{n+1}^A} "11";"12"};
{\ar^{z_{n+1}^A} "12";"13"};
{\ar^{d_{n+1}^A} "13";"14"};
{\ar_{f_n} "11";"21"};
{\ar|*+{_{Z^{n+1}(f_{\bullet})}} "12";"22"};
{\ar^{f_{n+1}} "13";"23"};
{\ar^{f_{n+2}} "14";"24"};
{\ar_{k_{n+1}^B} "21";"22"};
{\ar_{z_{n+1}^B} "22";"23"};
{\ar_{d_{n+1}^B} "23";"24"};
{\ar_{g_n} "21";"31"};
{\ar|*+{_{Z^{n+1}(g_{\bullet})}} "22";"32"};
{\ar^{g_{n+1}} "23";"33"};
{\ar^{g_{n+2}} "24";"34"};
{\ar_{k_{n+1}^C} "31";"32"};
{\ar_{z_{n+1}^C} "32";"33"};
{\ar_{d_{n+1}^C} "33";"34"};
{\ar@{=>}_{\underline{\lambda}_{n+1,2}} (-24,9)*{};(-24,5)*{}};
{\ar@{=>}^{\underline{\lambda}_{n+1,1}} (0,9)*{};(0,5)*{}};
{\ar@{=>}^{\lambda_{n+1}} (24,9)*{};(24,5)*{}};
{\ar@{=>}_{\underline{\kappa}_{n+1,2}} (-24,-5)*{};(-24,-9)*{}};
{\ar@{=>}^{\underline{\kappa}_{n+1,1}} (0,-5)*{};(0,-9)*{}};
{\ar@{=>}^{\kappa_{n+1}} (24,-5)*{};(24,-9)*{}};
\endxy
\]
such that
\begin{eqnarray*}
(z_{n+1}^A\circ\lambda_{n+1})\cdot(\underline{\lambda}_{n+1,1}\circ d_{n+1}^B)\cdot(Z^{n+1}(f_{\bullet})\circ\zeta_{n+1}^B)\cdot Z^{n+1}(f_{\bullet})_I^{\sharp}\\
=(\zeta_{n+1}^A\circ f_{n+2})\cdot(f_{n+2})_I^{\flat}.
\end{eqnarray*}
Then, it can be shown that each of the factorizations
\begin{eqnarray*}
Q^n(f_{\bullet})\circ\xi_n^B&:&Q^n(f_{\bullet})\circ x_n^B\circ z_{n+1}^B\Longrightarrow Q^n(f_{\bullet})\circ\ell_n^B\\
(x_n^A\circ\underline{\lambda}_{n+1,1}^{-1})\cdot(\xi_n^A\circ f_{n+1})\cdot\overline{\lambda}_{n-1,2}&:&x_n^A\circ Z^{n+1}(f_{\bullet})\circ z_{n+1}^B\Longrightarrow Q^n(f_{\bullet})\circ\ell_n^B
\end{eqnarray*}
are compatible with $\zeta_{n+1}^B$ and $(Q^n(f_{\bullet})\circ\mu_{n,2}^B)\cdot(Q^n(f_{\bullet}))_I^{\sharp}$.
\[
\xy
(-17,18)*+{Q^n(A_{\bullet})}="10";
(-17,-18)*+{Z^{n+1}(B_{\bullet})}="0";
(-3,0)*+{B_{n+1}}="6";
(17,0)*+{B_{n+2}}="8";
{\ar^{d_{n+1}^B} "6";"8"};
{\ar_{z^B_{n+1}} "0";"6"};
{\ar@/_0.80pc/_{0} "0";"8"};
{\ar@/^0.80pc/^{0} "10";"8"};
{\ar@/_0.80pc/|*+{_{Q^n(f_{\bullet})\circ x^B_n}} "10";"0"};
{\ar|*+{_{Q^n(f_{\bullet})\circ\ell^B_n}} "10";"6"};
{\ar@{=>} (-2,4)*{};(1,10)*{}};
{\ar@{..}@/_0.80pc/ (0,7)*{};(13,20)*{_{(Q^n(f_{\bullet})\circ\mu^B_{n,2})\cdot(Q^n(f_{\bullet}))^{\sharp}_I}}};
{\ar@{=>}^{\zeta^B_{n+1}} (-2,-4)*{};(1,-9)*{}};
{\ar@{=>} (-13,-6)*{};(-9,2)*{}};
{\ar@{..}@/^0.80pc/ (-11,-2)*{};(1,-20)*{_{Q^n(f_{\bullet})\circ\xi^B_n}}};
\endxy
\xy
(-17,18)*+{Q^n(A_{\bullet})}="10";
(-17,-18)*+{Z^{n+1}(B_{\bullet})}="0";
(-3,0)*+{B_{n+1}}="6";
(17,0)*+{B_{n+2}}="8";
{\ar^{d_{n+1}^B} "6";"8"};
{\ar_{z^B_{n+1}} "0";"6"};
{\ar@/_0.80pc/_{0} "0";"8"};
{\ar@/^0.80pc/^{0} "10";"8"};
{\ar@/_0.80pc/|*+{_{x_n^A\circ Z^{n+1}(f_{\bullet})}} "10";"0"};
{\ar|*+{_{Q^n(f_{\bullet})\circ\ell^B_n}} "10";"6"};
{\ar@{=>}^{\zeta^B_{n+1}} (-2,-4)*{};(1,-9)*{}};
{\ar@{=>} (-2,4)*{};(1,10)*{}};
{\ar@{..}@/_0.80pc/ (0,7)*{};(13,20)*{_{(Q^n(f_{\bullet})\circ\mu^B_{n,2})\cdot(Q^n(f_{\bullet}))^{\sharp}_I}}};
{\ar@{=>}(-13,-6)*{};(-9,2)*{}};
{\ar@{..}@/_0.80pc/ (-11,-2)*{};(1,-26)*{_{(x_n^A \circ\underline{\lambda}^{-1}_{n+1,1})\cdot(\xi^A_n\circ f_{n+1})\cdot\overline{\lambda}_{n-1,2}}}};
\endxy
\]
So, by the universality of the relative kernel, there exists a unique 2-cell $\widetilde{\lambda}_n\in\mathbf{S}^2(Q^n(f_{\bullet})\circ x_n^B,x_n^A\circ Z^{n+1}(f_{\bullet}))$ such that
\[ (\widetilde{\lambda}_n\circ z_{n+1}^B)\cdot(x_n^A\circ\underline{\lambda}_{n+1,1}^{-1})\cdot(\xi_n^A\circ f_{n+1})\cdot\overline{\lambda}_{n-1,2}=Q^n(f_{\bullet})\circ\xi_n^B. \]
This $\widetilde{\lambda}_n$ also satisfies $(q_n^A\circ\widetilde{\lambda}_n)\cdot(\eta_n^A\circ Z^{n+1}(f_{\bullet}))\cdot\underline{\lambda}_{n+1,2}=(\overline{\lambda}_{n-1,1}\circ x_n^B)\cdot(f_n\circ\eta_n^B)$ (see Remark \ref{AlsoRemark}). Similarly, we obtain a 2-cell $\widetilde{\kappa}_n\in\mathbf{S}^2(Q^n(g_{\bullet})\circ x_n^C,x_n^B\circ Z^{n+1}(g_{\bullet}))$ such that
\[ (\widetilde{\kappa}_n\circ z_{n+1}^C)\cdot(x_n^B\circ\underline{\kappa}_{n+1,1}^{-1})\cdot(\xi_n^B\circ g_{n+1})\cdot\overline{\kappa}_{n-1,2}=Q^n(g_{\bullet})\circ\xi_n^C. \]
In the rest, we show the following:
\begin{equation}
(Q^n(f_{\bullet})\circ\widetilde{\kappa}_n)\cdot(\widetilde{\lambda}_n\circ Z^{n+1}(g_{\bullet}))\cdot(x_n^A\circ Z^{n+1}(\varphi_{\bullet}))\cdot(x_n^A)_I^{\sharp}=(Q^n(\varphi_{\bullet})\circ x_n^C)\cdot(x_n^C)_I^{\flat}.
\label{SMCEQ}
\end{equation}
We have the following equalities:
\begin{eqnarray*}
\lefteqn{(Q^n(f_{\bullet})\circ\widetilde{\kappa}_n\circ z_{n+1}^C)\cdot(\widetilde{\lambda}_n\circ Z^{n+1}(g_{\bullet})\circ z_{n+1}^C)}\hspace{0.02cm}\\
&=&(Q^n(f_{\bullet})\circ Q^n(g_{\bullet})\circ\xi_n^C)\cdot(Q^n(f_{\bullet})\circ\overline{\kappa}_{n-1,2}^{-1})\cdot(\overline{\lambda}_{n-1,2}^{-1}\circ g_{n+1})\\
& &\!\!\!\!\!\!\!\cdot\ ((\xi_n^A)^{-1}\circ f_{n+1}\circ g_{n+1})\cdot(x_n^A\circ\underline{\lambda}_{n+1,1}\circ g_{n+1})\cdot(x_n^A\circ Z^{n+1}(f_{\bullet})\circ\underline{\kappa}_{n+1,1}),
\end{eqnarray*}
\begin{eqnarray*}
\lefteqn{(Q^n(\varphi_{\bullet})\circ x_n^C\circ z_{n+1}^C)\cdot(x_n^C\circ z_{n+1}^C)_I^{\flat}\cdot(x_n^A\circ z_{n+1}^A)_I^{\sharp-1}\cdot(x_n^A\circ z_{n+1}^A\circ\varphi_{n+1}^{-1})}\hspace{2cm}\\
&=&(Q^n(f_{\bullet})\circ Q^n(g_{\bullet})\circ\xi_n^C)\cdot(Q^n(f_{\bullet})\circ\overline{\kappa}_{n-1,2}^{-1})\cdot(\overline{\lambda}_{n-1,2}^{-1}\circ g_{n+1})\\
& &\cdot\ ((\xi_n^A)^{-1}\circ f_{n+1}\circ g_{n+1}),
\end{eqnarray*}
\begin{eqnarray*}
(z_{n+1}^C)_I^{\flat-1}&=&(x_n^A\circ z_{n+1}^A)_I^{\sharp-1}\cdot(x_n^A\circ z_{n+1}^A\circ\varphi_{n+1}^{-1})\cdot(x_n^A\circ\underline{\lambda}_{n+1,1}\circ g_{n+1})\\
& &\!\!\!\!\!\!\!\cdot\ (x_n^A\circ Z^{n+1}(f_{\bullet})\circ\underline{\kappa}_{n+1,1})\cdot(x_n^A\circ Z^{n+1}(\varphi_{\bullet})\circ z_{n+1}^C)\cdot((x_n^A)_I^{\sharp}\circ z_{n+1}^C).
\end{eqnarray*}
From these equalities and the faithfulness of $z_{n+1}^C$, we obtain $(\ref{SMCEQ})$.
\end{proof}

\begin{remark}\label{AlsoRemark}
It can be also shown that $\widetilde{\lambda}_n$ in the proof of Lemma \ref{FINLem3} satisfies
\[ (q_n^A\circ\widetilde{\lambda}_n)\cdot(\eta_n^A\circ Z^{n+1}(f_{\bullet}))\cdot\underline{\lambda}_{n+1,2}=(\overline{\lambda}_{n-1,1}\circ x_n^B)\cdot(f_n\circ\eta_n^B). \]
\end{remark}

By Lemma \ref{FINLem2} and Lemma \ref{FINLem3}, Theorem \ref{MainTheorem} is reduced to the following Proposition:
\begin{proposition}
\label{RevisionLabel2}
Consider the following diagram in $\mathbf{S}$, where $(A_{\bullet},d_{\bullet}^A,\delta_{\bullet}^A)$ is a complex which is relatively 2-exact in $A_2$ and $A_3$, and $(B_{\bullet},d_{\bullet}^B,\delta_{\bullet}^B)$ is a complex which is relatively 2-exact in $B_1$ and $B_2$.
\[
\xy
(-16,6)*+{A_1}="0";
(-16,-6)*+{B_1}="2";
(0,6)*+{A_2}="4";
(0,-6)*+{B_2}="6";
(16,6)*+{A_3}="8";
(32,6)*+{0}="18";
(16,-6)*+{B_3}="10";
(-32,-6)*+{0}="12";
{\ar^{d_1^A} "0";"4"};
{\ar_{f_1} "0";"2"};
{\ar^{0} "8";"18"};
{\ar_{0} "12";"2"};
{\ar_{d_1^B} "2";"6"};
{\ar|*+{_{f_2}} "4";"6"};
{\ar^{f_3} "8";"10"};
{\ar^{d_2^A} "4";"8"};
{\ar_{d_2^B} "6";"10"};
{\ar@/^1.65pc/^{0} "0";"8"};
{\ar@/_1.65pc/_{0} "2";"10"};
{\ar@{=>}_{\lambda_1} (-8,2)*{};(-8,-2)*{}};
{\ar@{=>}^{\lambda_2} (8,2)*{};(8,-2)*{}};
{\ar@{=>}_{\delta_2^A} (0,8)*{};(0,12)*{}};
{\ar@{=>}_{\delta_2^B} (0,-8)*{};(0,-12)*{}};
\endxy
\]
Assume $f_{\bullet}:A_{\bullet}\longrightarrow B_{\bullet}$ is a complex morphism. Then there exist $d\in\mathbf{S}^1(\mathrm{Ker}(f_3),\mathrm{Cok}(f_1))$, $\alpha\in\mathbf{S}^2(\underline{d}_2^A\circ d,0)$ and $\beta\in\mathbf{S}^2(d\circ\overline{d}_1^B,0)$ such that the sequence
\begin{equation}
\xy
(-50,0)*+{\mathrm{Ker}(f_1)}="1";
(-30,0)*+{\mathrm{Ker}(f_2)}="2";
(-10,0)*+{\mathrm{Ker}(f_3)}="3";
(10,0)*+{\mathrm{Cok}(f_1)}="4";
(30,0)*+{\mathrm{Cok}(f_2)}="5";
(50,0)*+{\mathrm{Cok}(f_3)}="6";
{\ar_{\underline{d}^A_1} "1";"2"};
{\ar_{\underline{d}^A_2} "2";"3"};
{\ar_{d} "3";"4"};
{\ar^{\overline{d}^B_1} "4";"5"};
{\ar^{\overline{d}^B_2} "5";"6"};
{\ar@/^1.65pc/^{0} "1";"3"};
{\ar@/_1.65pc/_{0} "2";"4"};
{\ar@/^1.65pc/^{0} "3";"5"};
{\ar@/_1.65pc/_{0} "4";"6"};
{\ar@{=>}_{\underline{\delta}^A_2} (-30,2)*{};(-30,6)*{}};
{\ar@{=>}^{\alpha} (-10,-2)*{};(-10,-6)*{}};
{\ar@{=>}_{\beta} (10,2)*{};(10,6)*{}};
{\ar@{=>}_{\overline{\delta}_2^B} (30,-2)*{};(30,-6)*{}};
\endxy
\label{SnakeSeq1}
\end{equation}
is 2-exact in $\mathrm{Ker}(f_2),\mathrm{Ker}(f_3),\mathrm{Cok}(f_1),\mathrm{Cok}(f_2)$.
\end{proposition}

\begin{remark}
\label{B4Remark}
This sequence does not necessarily become a complex. Indeed, for a relatively exact 2-category $\mathbf{S}$, the following are shown to be equivalent by an easy diagrammatic argument:

{\rm (i)} Any $(\ref{SnakeSeq1})$ obtained in Proposition \ref{RevisionLabel2} becomes a complex.

{\rm (ii)} For any $f\in \mathbf{S}^1(A,B)$,
\begin{equation}
\xy
(-32,0)*+{\mathrm{Ker}(f)}="1";
(-16,0)*+{A}="2";
(0,0)*+{B}="3";
(16,0)*+{\mathrm{Cok}(f)}="4";
{\ar_{k(f)} "1";"2"};
{\ar_{f} "2";"3"};
{\ar_{c(f)} "3";"4"};
{\ar@/^1.65pc/^{0} "1";"3"};
{\ar@/_1.65pc/_{0} "2";"4"};
{\ar@{=>}_<<<{\varepsilon_f} (-16,2)*{};(-16,6)*{}};
{\ar@{=>}^<<<{\pi_f} (0,-2)*{};(0,-6)*{}};
\endxy
\label{SnakeSeq2}
\end{equation}
is a complex.

(Indeed, if $(\ref{SnakeSeq2})$ is a complex for each of $f_1,f_2$ and $f_3$, then $(\ref{SnakeSeq1})$ becomes a complex.)

Thus if $\mathbf{S}$ satisfies {\rm (ii)}, then the long cohomology sequence in Theorem \ref{MainTheorem} becomes a complex. But this assumption is a bit too strong, since it is not satisfied by $\mathrm{SCG}$. This is pointed out by the referee.

\end{remark}

\begin{proof}[ of Proposition \ref{RevisionLabel2}]
Put $\mathrm{Ker}(d_2^A\circ f_3)=[K,k,\zeta]$. If we take the kernel of the diagram
\begin{eqnarray}
\xy
(-16,7)*+{A_1}="0";
(-16,-7)*+{0}="2";
(0,7)*+{A_2}="4";
(0,-7)*+{B_3}="6";
(16,7)*+{A_3}="8";
(16,-7)*+{B_3}="10";
(32,7)*+{0}="12";
(32,-7)*+{0}="14";
{\ar^{d_1^A} "0";"4"};
{\ar_{0} "0";"2"};
{\ar_{0} "2";"6"};
{\ar|*+{_{d_2^A\circ f_3}} "4";"6"};
{\ar^{f_3} "8";"10"};
{\ar^{d_2^A} "4";"8"};
{\ar@{=}_{\mathrm{id}} "6";"10"};
{\ar^{0} "8";"12"};
{\ar_{0} "10";"14"};
{\ar@{}|\circlearrowright"4";"10"};
{\ar@/^1.65pc/^{0} "0";"8"};
{\ar@{=>}_{\xi_0} (-8,2)*{};(-8,-2)*{}};
{\ar@{=>}_{\delta_2^A } (0,9)*{};(0,13)*{}};
\endxy
\end{eqnarray}
where $\xi_0:=(\delta_2^A\circ f_3)\cdot(f_3)_I^{\flat}$, then by Proposition \ref{final6} we obtain a diagram
\begin{equation}
\xy
(-18,6)*+{A_1}="0";
(-18,-6)*+{A_1}="2";
(0,6)*+{K}="4";
(0,-6)*+{A_2}="6";
(18,6)*+{\mathrm{Ker}(f_3)}="8";
(18,-6)*+{A_3}="10";
(36,-6)*+{0}="14";
{\ar^{k_1} "0";"4"};
{\ar@{=}_{\mathrm{id}} "0";"2"};
{\ar_{d_1^A} "2";"6"};
{\ar_{k} "4";"6"};
{\ar^{k(f_3)} "8";"10"};
{\ar^{k_2} "4";"8"};
{\ar_{d_2^A} "6";"10"};
{\ar_{0} "10";"14"};
{\ar@/^1.65pc/^{0} "0";"8"};
{\ar@/_1.65pc/_{0} "2";"10"};
{\ar@{=>}_{\xi_1} (-9,2)*{};(-9,-2)*{}};
{\ar@{=>}_{\xi_2} (9,2)*{};(9,-2)*{}};
{\ar@{=>}_{\alpha_2} (0,8)*{};(0,12)*{}};
{\ar@{=>}_{\delta_2^A} (0,-8)*{};(0,-12)*{}};
\endxy
\end{equation}
which satisfies
\begin{eqnarray*}
(k_2\circ\varepsilon_{f_3})\cdot(k_2)_I^{\sharp}&=&(\xi_2\circ f_3)\cdot\zeta\\
(\xi_1\circ d_2^A\circ f_3)\cdot\xi_0&=&(k_1\circ\zeta)\cdot(k_1)_I^{\sharp}\\
(k_1\circ\xi_2)\cdot(\xi_1\circ d_2^A)\cdot\delta_2^A&=&(\alpha_2\circ k(f_3))\cdot k(f_3)_I^{\flat}.
\end{eqnarray*}
By Lemma \ref{FFFFLem4},
\[
\xy
(-8,0)*+{A_1}="4";
(8,0)*+{K}="6";
(24,0)*+{\mathrm{Ker}(f_3)}="8";
{\ar_{k_1} "4";"6"};
{\ar_{k_2} "6";"8"};
{\ar@/^1.65pc/^{0} "4";"8"};
{\ar@{=>}_{\alpha_2} (8,2)*{};(8,6)*{}};
\endxy
\]
is 2-exact in $K$. On the other hand, by (the dual of) Proposition \ref{final12},
\[
\xy
(-10,6)*+{K}="0";
(-10,-6)*+{A_2}="2";
(10,6)*+{\mathrm{Ker}(f_3)}="4";
(10,-6)*+{A_3}="6";
{\ar^{k_2} "0";"4"};
{\ar_{k} "0";"2"};
{\ar_{d_2^A} "2";"6"};
{\ar^{f_3} "4";"6"};
{\ar@{=>}_{\xi_2} (0,2)*{};(0,-2)*{}};
\endxy
\]
is a pullback diagram, and $k_2$ becomes cofaithful since $d_2^A$ is cofaithful. Thus, we have $\mathrm{Cok}(k_1)=[\mathrm{Ker}(f_3),k_2,\alpha_2]$ by Lemma \ref{FFFFLem1}. Dually, if we put $\mathrm{Cok}(f_1\circ d_1^B)=[Q,q,\rho]$, then we obtain the following diagram
\[
\xy
(-27,14)*+{0}="10";
(-9,14)*+{A_1}="11";
(9,14)*+{A_1}="12";
(27,14)*+{0}="13";
(-27,0)*+{0}="20";
(-9,0)*+{B_1}="21";
(9,0)*+{B_2}="22";
(27,0)*+{B_3}="23";
(-9,-14)*+{\mathrm{Cok}(f_1)}="31";
(9,-14)*+{Q}="32";
(27,-14)*+{B_3}="33";
{\ar^{0} "10";"11"};
{\ar@{=}^{\mathrm{id}} "11";"12"};
{\ar^{0} "12";"13"};
{\ar_{f_1} "11";"21"};
{\ar|*+{_{f_1\circ d_1^B}} "12";"22"};
{\ar^{0} "13";"23"};
{\ar^{0} "20";"21"};
{\ar|*+{_{d_1^B}} "21";"22"};
{\ar|*+{_{d_2^B}} "22";"23"};
{\ar_{c(f_1)} "21";"31"};
{\ar|*+{_{q}} "22";"32"};
{\ar@{=}^{\mathrm{id}} "23";"33"};
{\ar_{q_1} "31";"32"};
{\ar_{q_2} "32";"33"};
{\ar@/_1.65pc/_{0} "31";"33"};
{\ar@{}|\circlearrowright"11";"22"};
{\ar@{=>}^{\eta_0} (18,9)*{};(18,5)*{}};
{\ar@{=>}_{\eta_1} (0,-5)*{};(0,-9)*{}};
{\ar@{=>}_{\eta_2} (16,-5)*{};(16,-9)*{}};
{\ar@{=>}_{\beta_2} (8,-16)*{};(8,-20)*{}};
\endxy
\]
which satisfies
\begin{eqnarray*}
\eta_0&=&(f_1)_I^{\sharp-1}\cdot(f_1\circ\delta_2^{B-1})\\
\rho&=&(f_1\circ\eta_1)\cdot(\pi_{f_1}\circ q_1)\circ(q_1)_I^{\flat}\\
\mathrm{id}_0&=&\eta_0\cdot(f_1\circ d_1^B\circ\eta_2)\cdot(\rho\circ q_2)\cdot(q_2)_I^{\flat}\\
\delta_2^B&=&(d_1^B\circ\eta_2)\cdot(\eta_1\circ q_2)\cdot(c(f_1)\circ\beta_2)\cdot c(f_1)_I^{\sharp},
\end{eqnarray*}
and we have $\mathrm{Ker}(q_2)=[\mathrm{Cok}(f_1),q_1,\beta_2]$. (The \textquotedblleft un-duality\textquotedblright\ in appearance is simply because of the direction of the 2-cells.) Thus, we obtain complex morphisms:
\[
\xy
(-32,14)*+{A_1}="01";
(-32,0)*+{K}="02";
(-32,-14)*+{\mathrm{Ker}(f_3)}="03";
(-16,14)*+{A_1}="11";
(-16,0)*+{A_2}="12";
(-16,-14)*+{A_3}="13";
(0,14)*+{B_1}="21";
(0,0)*+{B_2}="22";
(0,-14)*+{B_3}="23";
(16,14)*+{\mathrm{Cok}(f_1)}="31";
(16,0)*+{Q}="32";
(16,-14)*+{B_3}="33";
{\ar_{k_1} "01";"02"};
{\ar_{k_2} "02";"03"};
{\ar@{=}^{\mathrm{id}} "01";"11"};
{\ar|*+{_k} "02";"12"};
{\ar_{k(f_3)} "03";"13"};
{\ar|*+{_{d_1^A}} "11";"12"};
{\ar|*+{_{d_2^A}} "12";"13"};
{\ar^{f_1} "11";"21"};
{\ar|*+{_{f_2 }} "12";"22"};
{\ar_{f_3} "13";"23"};
{\ar|*+{_{d_1^B}} "21";"22"};
{\ar|*+{_{d_2^B}} "22";"23"};
{\ar^{c(f_1)} "21";"31"};
{\ar|*+{_{q}} "22";"32"};
{\ar@{=}_{\mathrm{id}} "23";"33"};
{\ar^{q_1} "31";"32"};
{\ar^{q_2} "32";"33"};
{\ar@{=>}^{\xi_1} (-27,7)*{};(-21,7)*{}};
{\ar@{=>}^{\xi_2} (-27,-7)*{};(-21,-7)*{}};
{\ar@{=>}^{\lambda_1} (-11,7)*{};(-5,7)*{}};
{\ar@{=>}^{\lambda_2} (-11,-7)*{};(-5,-7)*{}};
{\ar@{=>}^{\eta_1} (5,7)*{};(11,7)*{}};
{\ar@{=>}^{\eta_2} (5,-7)*{};(11,-7)*{}};
\endxy
\]

If we put
\begin{eqnarray*}
c&:=&k\circ f_2\circ q\nonumber\\
\alpha_K&:=&(\xi_1\circ f_2\circ q)\cdot(\lambda_1\circ q)\cdot\rho\\
\beta_Q&:=&(k\circ f_2\circ\eta_2^{-1})\cdot(k\circ\lambda_2^{-1})\cdot\zeta,
\end{eqnarray*}
then, it can be shown that the following diagram is a complex.
\[
\xy
(-24,0)*+{A_1}="2";
(-8,0)*+{K}="3";
(8,0)*+{Q}="4";
(24,0)*+{B_3}="5";
{\ar_{k_1} "2";"3"};
{\ar_{c} "3";"4"};
{\ar^{q_2} "4";"5"};
{\ar@/^1.65pc/^{0} "2";"4"};
{\ar@/_1.65pc/_{0} "3";"5"};
{\ar@{=>}^{\alpha_K} (-8,2)*{};(-8,6)*{}};
{\ar@{=>}^{\beta_Q} (8,-2)*{};(8,-6)*{}};
\endxy
\]
Since $\mathrm{Cok}(k_1)=[\mathrm{Ker}(f_3),k_2,\alpha_2]$ as already shown, we have a factorization diagram
\[
\xy
(0,18)*+{\mathrm{Ker}(f_3)}="0";
(-36,0)*+{A_1}="2";
(-12,0)*+{K}="4";
(12,0)*+{Q}="6";
(36,0)*+{B_3}="8";
{\ar_{k_1} "2";"4"};
{\ar_{c} "4";"6"};
{\ar^{k_2} "4";"0"};
{\ar_{q_2} "6";"8"};
{\ar^{\overline{c}} "0";"6"};
{\ar@{=>}^{\overline{\beta}_Q} (14,2)*{};(18,8)*{}};
{\ar@{=>}_{\alpha_2} (-15,3)*{};(-18,8)*{}};
{\ar@/_1.65pc/_{0} "2";"6"};
{\ar@/_1.65pc/_{0} "4";"8"};
{\ar@/^0.50pc/^{0} "0";"8"};
{\ar@/^0.50pc/^{0} "2";"0"};
{\ar@{=>}^{\beta_Q} (12,-2)*{};(12,-6)*{}};
{\ar@{=>}_{\alpha_K} (-12,-2)*{};(-12,-6)*{}};
{\ar@{=>}_{\overline{\alpha}_K} (0,10)*{};(0,4)*{}};
\endxy
\]
which satisfies
\begin{eqnarray*}
(k_1\circ\overline{\alpha}_K)\cdot\alpha_K&=&(\alpha_2\circ\overline{c})\cdot(\overline{c})_I^{\flat}\\
(\overline{\alpha}_K\circ q_2)\cdot\beta_Q&=&(k_2\circ\overline{\beta}_Q)\cdot(k_2)_I^{\sharp}.
\end{eqnarray*}
Similarly, since $\mathrm{Ker}(q_2)=[\mathrm{Cok}(f_1),q_1,\beta_2]$, we have a
factorization diagram
\[
\xy
(0,-18)*+{\mathrm{Cok}(f_1)}="0";
(-36,0)*+{A_1}="2";
(-12,0)*+{K}="4";
(12,0)*+{Q}="6";
(36,0)*+{B_3}="8";
{\ar^{k_1} "2";"4"};
{\ar^{c} "4";"6"};
{\ar_{\underline{c}} "4";"0"};
{\ar^{q_2} "6";"8"};
{\ar_{q_1} "0";"6"};
{\ar@{=>}^{\beta_2} (14,-2)*{};(18,-8)*{}};
{\ar@{=>}_{\underline{\alpha} _K} (-15,-3)*{};(-18,-8)*{}};
{\ar@/^1.65pc/^{0} "2";"6"};
{\ar@/^1.65pc/^{0} "4";"8"};
{\ar@/_0.50pc/_{0} "0";"8"};
{\ar@/_0.50pc/_{0} "2";"0"};
{\ar@{=>}_{\beta_Q} (12,2)*{};(12,6)*{}};
{\ar@{=>}^{\alpha_K} (-12,2)*{};(-12,6)*{}};
{\ar@{=>}_{\underline{\beta} _Q} (0,-10)*{};(0,-4)*{}};
\endxy
\]
which satisfies
\begin{eqnarray*}
(\underline{\beta}_Q\circ q_2)\cdot\beta_Q&=&(\underline{c}\circ\beta_2)\cdot(\underline{c})_I^{\sharp}\\
(k_1\circ\underline{\beta}_Q)\cdot\alpha_K&=&(\underline{\alpha}_K\circ q_1)\cdot(q_1)_I^{\flat}.
\end{eqnarray*}
Then, there exist $d\in\mathbf{S}^1(\mathrm{Ker}(f_3),\mathrm{Cok}(f_1))$, $\alpha^{\dag}\in\mathbf{S}^2(k_2\circ d,\underline{c})$ and $\beta^{\dag}\in\mathbf{S}^2(d\circ q_1,\overline{c})$
such that
\begin{eqnarray*}
(k_1\circ\alpha^{\dag})\cdot\underline{\alpha}_K&=&(\alpha_2\circ d)\cdot d_I^{\flat}\\
(\beta^{\dag}\circ q_2)\cdot\overline{\beta}_Q&=&(d\circ\beta_2)\cdot d_I^{\sharp}\\
(k_2\circ\beta^{\dag})\cdot\overline{\alpha}_K&=&(\alpha^{\dag}\circ q_1)\cdot\underline{\beta}_Q
\end{eqnarray*}
(note that $\mathrm{Cok}(k_1)=[\mathrm{Ker}(f_3),k_2,\alpha_2]$ and $\mathrm{Ker}(q_2)=[\mathrm{Cok}(f_1),q_1,\beta_2]$ (cf. Lemma \ref{FINLem1})):
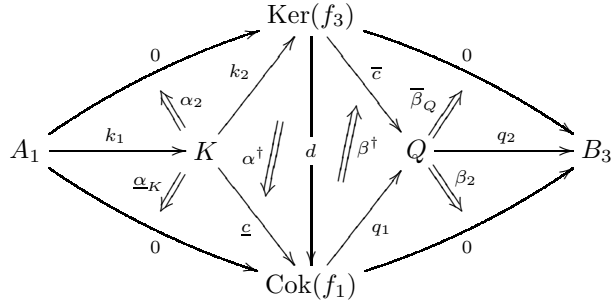
\[
\xy
(0,18)*+{\mathrm{Ker}(f_3)}="0";
(0,-18)*+{\mathrm{Cok}(f_1)}="20";
(-38,0)*+{A_1}="4";
(-14,0)*+{K}="6";
(14,0)*+{Q}="8";
(38,0)*+{B_3}="10";
{\ar^{k_2} "6";"0"};
{\ar_{\underline{c}} "6";"20"};
{\ar^{q_2} "8";"10"};
{\ar|*+{_{d}} "0";"20"};
{\ar^{k_1} "4";"6"};
{\ar^{\overline{c}} "0";"8"};
{\ar_{q_1} "20";"8"};
{\ar@/^0.50pc/^{0} "0";"10"};
{\ar@/^0.50pc/^{0} "4";"0"};
{\ar@/_0.50pc/_{0} "20";"10"};
{\ar@/_0.50pc/_{0} "4";"20"};
{\ar@{=>}^{\overline{\beta}_Q} (16,2)*{};(20,8)*{}};
{\ar@{=>}_{\alpha_2} (-17,3)*{};(-20,8)*{}};
{\ar@{=>}^{\beta_2} (16,-2)*{};(20,-8)*{}};
{\ar@{=>}_{\underline{\alpha}_K} (-17,-3)*{};(-20,-8)*{}};
{\ar@{=>}_{\alpha^{\dag}} (-4,4)*{};(-6,-6)*{}};
{\ar@{=>}_{\beta^{\dag}} (4,-4)*{};(6,6)*{}};
\endxy
\]
Applying (the dual of) Proposition \ref{final9} to the diagram
\[
\xy
(-20,14)*+{\mathrm{Ker}(f_1)}="01";
(0,14)*+{\mathrm{Ker}(f_2)}="02";
(20,14)*+{\mathrm{Ker}(f_3)}="03";
(-20,0)*+{A_1}="11";
(0,0)*+{A_2}="12";
(20,0)*+{A_3}="13";
(-20,-14)*+{0}="21";
(0,-14)*+{B_3}="22";
(20,-14)*+{B_3}="23";
(23,-15)*+{,}="24";
{\ar^{\underline{d}_1^A} "01";"02"};
{\ar^{\underline{d}_2^A} "02";"03"};
{\ar_{k(f_1)} "01";"11"};
{\ar|*+{_{k(f_2)}} "02";"12"};
{\ar^{k(f_3)} "03";"13"};
{\ar|*+{_{d_1^A}} "11";"12"};
{\ar|*+{_{d_2^A}} "12";"13"};
{\ar_{0} "11";"21"};
{\ar|*+{_{d_2^A\circ f_3}} "12";"22"};
{\ar^{f_3} "13";"23"};
{\ar_{0} "21";"22"};
{\ar@{=}_{\mathrm{id}} "22";"23"};
{\ar@{}|\circlearrowright"12";"23"};
{\ar@{=>}_{\underline{\lambda}_1} (-10,9)*{};(-10,5)*{}};
{\ar@{=>}_{\underline{\lambda}_2} (10,9)*{};(10,5)*{}};
{\ar@{=>}_{\xi_0} (-10,-5)*{};(-10,-9)*{}};
\endxy
\]
we see that there exist $k^{\prime}\in\mathbf{S}^1(\mathrm{Ker}(f_2),K)$, $\xi_1^{\prime}\in\mathbf{S}^2(\underline{d}_1^A\circ k^{\prime},k(f_1)\circ k_1)$, $\xi_2^{\prime}\in\mathbf{S}^2(\underline{d}_2^A,k^{\prime}\circ k_2)$ and $\xi\in\mathbf{S}^2(k^{\prime}\circ k,k(f_2))$ 
such that
\begin{eqnarray*}
\underline{\delta}_2^A&=&(\underline{d}_1^A\circ\xi_2^{\prime})\cdot(\xi_1^{\prime}\circ k_2)\cdot(k(f_1)\circ\alpha_2)\cdot k(f_1)_I^{\sharp}\\
\underline{\lambda}_2&=&(\xi_2^{\prime}\circ k(f_3))\cdot(k^{\prime}\circ\xi_2)\cdot(\xi\circ d_2^A)\\
(\underline{d}_1^A\circ\xi)\cdot\underline{\lambda}_1&=&(\xi_1^{\prime}\circ k)\cdot(k(f_1)\circ\xi_1).
\end{eqnarray*}
Similarly, there exist $q^{\prime}\in\mathbf{S}^1(Q,\mathrm{Cok}(f_2))$, $\eta_1^{\prime}\in\mathbf{S}^2(q_1\circ q^{\prime},\overline{d}_1^B)$, $\eta_2^{\prime}\in\mathbf{S}^2(q_2\circ c(f_3),q^{\prime}\circ\overline{d}_2^B)$ and $\eta\in\mathbf{S}^2(q\circ q^{\prime},c(f_2))$ such that
\begin{eqnarray*}
(\beta_2\circ c(f_3))\cdot c(f_3)_I^{\flat}&=&(q_1\circ\eta_2^{\prime})\cdot(\eta_1^{\prime}\circ\overline{d}_2^B)\cdot\overline{\delta}_2^B\\
(\eta_1\circ q^{\prime})\cdot(c(f_1)\circ\eta_1^{\prime})&=&(d_1^B\circ\eta)\cdot\overline{\lambda}_1\\
(\eta_2\circ c(f_3))\cdot(q\circ\eta_2^{\prime})\cdot(\eta\circ\overline{d}_2^B)&=&\overline{\lambda}_2.
\end{eqnarray*}
If we put
\[ \alpha_0:=(\underline{d}_2^A\circ\beta^{\dag})\cdot(\xi_2^{\prime}\circ\overline{c})\cdot(k^{\prime}\circ\overline{\alpha}_K)\cdot(\xi\circ f_2\circ q)\cdot(\varepsilon_{f_2}\circ q)\cdot q_I^{\flat}, \]
then it can be shown that $\alpha_0:\underline{d}_2^A\circ d\circ q_1\Longrightarrow0$ is compatible with $\beta_2$.
\[
\xy
(-36,0)*+{\mathrm{Ker}(f_2)}="2";
(-12,0)*+{\mathrm{Cok}(f_1)}="3";
(12,0)*+{Q}="4";
(36,0)*+{B_3}="5";
{\ar_{\underline{d}_2^A\circ d}  "2";"3"};
{\ar_{q_1} "3";"4"};
{\ar^{q_2} "4";"5"};
{\ar@/^1.65pc/^{0} "2";"4"};
{\ar@/_1.65pc/_{0} "3";"5"};
{\ar@{=>}^{\alpha_0} (-12,2)*{};(-12,6)*{}};
{\ar@{=>}^{\beta_2} (12,-2)*{};(12,-6)*{}};
\endxy
\]
So by Lemma \ref{KernelLemma}, there exists $\alpha\in\mathbf{S}^2(\underline{d}_2^A\circ d,0)$ such that
\[ (\alpha\circ q_1)\cdot(q_1)_I^{\flat}=\alpha_0. \]
Dually, if we put
\[ \beta_0:=(\alpha^{\dag}\circ\overline{d}_1^B)\cdot(\underline{c}\circ\eta_1^{\prime -1})\cdot(\underline{\beta}_Q\circ q^{\prime})\cdot(k\circ f_2\circ\eta)\cdot(k\circ\pi_{f_2})\cdot k_I^{\sharp}, \]
then $\beta_0:k_2\circ d\circ\overline{d}_1^B\Longrightarrow 0$ is compatible with $\alpha_2$, and there exists $\beta\in\mathbf{S}^2(d\circ\overline{d}_1^B,0)$ such that
\[ (k_2\circ\beta)\cdot(k_2)_I^{\sharp}=\beta_0. \]
\[
\xy
(-50,0)*+{\mathrm{Ker}(f_1)}="1";
(-30,0)*+{\mathrm{Ker}(f_2)}="2";
(-10,0)*+{\mathrm{Ker}(f_3)}="3";
(10,0)*+{\mathrm{Cok}(f_1)}="4";
(30,0)*+{\mathrm{Cok}(f_2)}="5";
(50,0)*+{\mathrm{Cok}(f_3)}="6";
{\ar_{\underline{d}^A_1} "1";"2"};
{\ar_{\underline{d}^A_2} "2";"3"};
{\ar_{d} "3";"4"};
{\ar^{\overline{d}^B_1} "4";"5"};
{\ar^{\overline{d}^B_2} "5";"6"};
{\ar@/^1.65pc/^{0} "1";"3"};
{\ar@/_1.65pc/_{0} "2";"4"};
{\ar@/^1.65pc/^{0} "3";"5"};
{\ar@/_1.65pc/_{0} "4";"6"};
{\ar@{=>}_{\underline{\delta}^A_2} (-30,2)*{};(-30,6)*{}};
{\ar@{=>}^{\alpha} (-10,-2)*{};(-10,-6)*{}};
{\ar@{=>}_{\beta} (10,2)*{};(10,6)*{}};
{\ar@{=>}_{\overline{\delta}_2^B} (30,-2)*{};(30,-6)*{}};
\endxy
\]
In the rest, we show that this is 2-exact in $\mathrm{Ker}(f_2),\mathrm{Ker}(f_3),\mathrm{Cok}(f_1),\mathrm{Cok}(f_2)$. We show only the 2-exactness in $\mathrm{Ker}(f_2)$ and $\mathrm{Ker}(f_3)$, since the rest can be shown dually. The 2-exactness in $\mathrm{Ker}(f_2)$ follows immediately from Lemma \ref{FFFFLem4}. So, we show the 2-exactness in $\mathrm{Ker}(f_3)$. Since we have $\mathrm{Cok}(d_1^A)=[A_3,d_2^A,\delta_2^A]$ and $\mathrm{Cok}(f_1\circ d_2^B)=[Q,q,\rho]$, there exists a factorization $(\ell,\varpi_1)$
\begin{equation}
\xy
(-16,6)*+{A_1}="0";
(-16,-6)*+{A_1}="2";
(0,6)*+{A_2}="4";
(0,-6)*+{B_2}="6";
(16,6)*+{A_3}="8";
(16,-6)*+{Q}="10";
(32,6)*+{0}="12";
(32,-6)*+{0}="14";
{\ar^{d^A_1} "0";"4"};
{\ar@{=}_{\mathrm{id}} "0";"2"};
{\ar_{f_1\circ d^B_1} "2";"6"};
{\ar_{f_2} "4";"6"};
{\ar^{\ell} "8";"10"};
{\ar^{d^A_2} "4";"8"};
{\ar_{q} "6";"10"};
{\ar^{0} "8";"12"};
{\ar_{0} "10";"14"};
{\ar@/^2.00pc/^{0} "0";"8"};
{\ar@/_2.00pc/_{0} "2";"10"};
{\ar@{=>}_{\lambda_1} (-8,2)*{};(-8,-2)*{}};
{\ar@{=>}_{\varpi_1} (8,2)*{};(8,-2)*{}};
{\ar@{=>}_{\delta^A_2 } (0,9)*{};(0,13)*{}};
{\ar@{=>}_{\rho} (0,-9)*{};(0,-13)*{}};
\endxy
\label{myan1}
\end{equation}
such that
\[ (d_1^A\circ\varpi_1)\cdot(\lambda_1\circ q)\cdot\rho=(\delta_2^A\circ\ell)\cdot\ell_I^{\flat}. \]
Applying Lemma \ref{FFFFLem2} to diagram $(\ref{myan1})$, we see that the following diagram becomes 2-exact in $A_3$:
\begin{equation}
\xy
(-24,0)*+{\mathrm{Ker}(f_2)}="4";
(0,0)*+{A_3}="6";
(24,0)*+{Q}="8";
{\ar_{k(f_2)\circ d^A_2} "4";"6"};
{\ar_{\ell} "6";"8"};
{\ar@/^1.65pc/^{0} "4";"8"};
{\ar@{=>}_{} (0,2)*{};(0,6)*{}};
{\ar@{..}@/_0.80pc/ (0,4)*{};(24,8)*{_{(k(f_2)\circ\varpi_1)\cdot(\varepsilon_{f_2}\circ q)\cdot q^{\flat}_I}}};
\endxy
\label{myandia}
\end{equation}
Then it can be shown that $(\varpi_1\circ q_2)\cdot(f_2\circ\eta_2^{-1}):d_2^A\circ\ell\circ q_2\Longrightarrow f_2\circ d_2^B$ is compatible with $\delta_2^A$ and $(\lambda_1\circ d_2^B)\cdot(f_1\circ d_2^B)\cdot(f_1)_I^{\sharp}$.
So, comparing the following two factorizations
\[
\xy
(-14,0)*+{A_1}="2";
(18,-14)*+{B_3}="4";
(4,0)*+{A_2}="6";
(18,14)*+{A_3}="8";
{\ar|*+{_{f_2\circ d_2^B}} "6";"4"};
{\ar^{d^A_1} "2";"6"};
{\ar^{d^A_2} "6";"8"};
{\ar@/_1.20pc/_{0} "2";"4"};
{\ar@/^0.80pc/|*+{_{\ell\circ q_2}} "8";"4"};
{\ar@/^1.20pc/^{0} "2";"8"};
{\ar@{=>}_{\delta^A_2} (2,3)*{};(0,9)*{}};
{\ar@{=>} (2,-3)*{};(0,-9)*{}};
{\ar@{=>} (16,5)*{};(14,-2)*{}};
{\ar@{..}@/^1.60pc/ (1,-6)*{};(-6,-19)*{_{(\lambda_1\circ d^B_2)\cdot(f_1\circ\delta^B_2)\cdot(f_1)^{\sharp}_I}}};
{\ar@{..}@/^1.60pc/ (15,2)*{};(10,20)*{_{(\varpi_1\circ q_2)\cdot(f_2\circ\eta^{-1}_2)}}};
\endxy
\quad
\xy
(-14,0)*+{A_1}="2";
(18,-14)*+{B_3}="4";
(21,-15)*+{,}="10";
(4,0)*+{A_2}="6";
(18,14)*+{A_3}="8";
{\ar|*+{_{f_2\circ d_2^B}} "6";"4"};
{\ar^{d^A_1} "2";"6"};
{\ar^{d^A_2} "6";"8"};
{\ar@/_1.20pc/_{0} "2";"4"};
{\ar@/^0.80pc/|*+{_{f_3}} "8";"4"};
{\ar@/^1.20pc/^{0} "2";"8"};
{\ar@{=>}_{\delta^A_2} (2,3)*{};(0,9)*{}};
{\ar@{=>} (2,-3)*{};(0,-9)*{}};
{\ar@{=>}^{\lambda_2} (16,5)*{};(14,-2)*{}};
{\ar@{..}@/^1.60pc/ (1,-6)*{};(-6,-19)*{_{(\lambda_1\circ d^B_2)\cdot(f_1\circ\delta^B_2)\cdot(f_1)^{\sharp}_I}}};
\endxy
\]
we see there exists a unique 2-cell $\varpi_2\in\mathbf{S}^2(\ell\circ q_2,f_3)$ such that
\[ (d_2^A\circ\varpi_2)\cdot\lambda_2=(\varpi_1\circ q_2)\cdot(f_2\circ\eta_2^{-1}). \]
Then it can be shown that each of the two factorizations
\begin{eqnarray*}
&\mathrm{(1)}&\ (\alpha^{\dag}\circ q_1)\cdot\underline{\beta}_Q:k_2\circ d\circ q_1\Longrightarrow c\\
&\mathrm{(2)}&\ (\xi_2\circ\ell)\cdot(k\circ\varpi_1):k_2\circ k(f_3)\circ\ell\Longrightarrow k\circ f_2\circ q=c
\end{eqnarray*}
is compatible with $\alpha_2$ and $\alpha_K$.
\[
\xy
(-14,0)*+{A_1}="2";
(18,-14)*+{Q}="4";
(4,0)*+{K}="6";
(18,14)*+{\mathrm{Ker}(f_3)}="8";
{\ar_{c} "6";"4"};
{\ar^{k_1} "2";"6"};
{\ar^{k_2} "6";"8"};
{\ar@/_1.20pc/_{0} "2";"4"};
{\ar@/^0.80pc/|*+{_{d\circ q_1}} "8";"4"};
{\ar@/^1.20pc/^{0} "2";"8"};
{\ar@{=>}_{\alpha_2} (2,3)*{};(0,9)*{}};
{\ar@{=>}_{\alpha_K} (2,-3)*{};(0,-9)*{}};
{\ar@{=>} (16,5)*{};(14,-2)*{}};
{\ar@{..}@/^1.60pc/ (15,2)*{};(10,20)*{_{(\alpha^{\dag}\circ q_1)\cdot\underline{\beta}_Q}}};
\endxy
\quad
\xy
(-14,0)*+{A_1}="2";
(18,-14)*+{Q}="4";
(4,0)*+{K}="6";
(18,14)*+{\mathrm{Ker}(f_3)}="8";
{\ar_{c} "6";"4"};
{\ar^{k_1} "2";"6"};
{\ar^{k_2} "6";"8"};
{\ar@/_1.20pc/_{0} "2";"4"};
{\ar@/^0.80pc/|*+{_{k(f_3)\circ\ell}} "8";"4"};
{\ar@/^1.20pc/^{0} "2";"8"};
{\ar@{=>}_{\alpha_2} (2,3)*{};(0,9)*{}};
{\ar@{=>}_{\alpha_K} (2,-3)*{};(0,-9)*{}};
{\ar@{=>}^{} (16,5)*{};(14,-2)*{}} ;
{\ar@{..}@/^1.60pc/^{} (15,2)*{};(10,20)*{_{(\xi_2\circ\ell)\cdot(k\circ\varpi_1)}}};
\endxy
\]
So there exists a unique 2-cell $\varpi_3\in\mathbf{S}^2(d\circ q_1,k(f_3)\circ\ell)$ such that
\[ (k_2\circ\varpi_3)\cdot(\xi_2\circ\ell)\cdot(k\circ\varpi_1)=(\alpha^{\dag}\circ q_1)\cdot\underline{\beta}_Q \]
(recall that $\mathrm{Cok}(k_1)=[\mathrm{Ker}(f_3),k_2,\alpha_2]$). Then we have $(\varpi_3\circ q_2)\cdot(k(f_3)\circ\varpi_2)\cdot\varepsilon_{f_3}=(d\circ\beta_2)\cdot d_I^{\sharp}$.
\begin{equation}
\xy
(-12,14)*+{\mathrm{Ker}(f_3)}="11";
(12,14)*+{\mathrm{Cok}(f_1)}="12";
(-12,0)*+{A_3}="21";
(12,0)*+{Q}="22";
(-12,-14)*+{B_3}="31";
(12,-14)*+{B_3}="32";
{\ar^{d} "11";"12"};
{\ar|*+{_{k(f_3)}} "11";"21"};
{\ar|*+{_{q_1}} "12";"22"};
{\ar_{\ell} "21";"22"};
{\ar|*+{_{f_3}} "21";"31"};
{\ar|*+{_{q_2}} "22";"32"};
{\ar@{=}_{\mathrm{id}} "31";"32"};
{\ar@/^1.80pc/^{0} "12";"32"};
{\ar@/_1.80pc/_{0} "11";"31"};
{\ar@{=>}_{\varpi_3} (0,9)*{};(0,5)*{}};
{\ar@{=>}_{\varpi_2} (0,-5)*{};(0,-9)*{}};
{\ar@{=>}_{\varepsilon_{f_3}} (-14,0)*{};(-18,0)*{}};
{\ar@{=>}_{\beta_2} (14,0)*{};(18,0)*{}};
\endxy
\label{myan5}
\end{equation}
By taking kernels of $d,\ell$ and $\mathrm{id}_{B_3}$ in $(\ref{myan5})$, we obtain the following diagram.
\[
\xy
(-24,14)*+{\mathrm{Ker}(d)}="01";
(0,14)*+{\mathrm{Ker}(f_3)}="02";
(24,14)*+{\mathrm{Cok}(f_1)}="03";
(-24,0)*+{\mathrm{Ker}(\ell)}="11";
(0,0)*+{A_3}="12";
(24,0)*+{Q}="13";
(-24,-14)*+{0}="21";
(0,-14)*+{B_3}="22";
(24,-14)*+{B_3}="23";
{\ar^{k(d) } "01";"02"};
{\ar^{d} "02";"03"};
{\ar_{\underline{k(f_3)}} "01";"11"};
{\ar|*+{_{k(f_3)}} "02";"12"};
{\ar^{q_1} "03";"13"};
{\ar_{k(\ell)} "11";"12"};
{\ar|*+{_{\ell}} "12";"13"};
{\ar_{0} "11";"21"};
{\ar|*+{_{f_3}} "12";"22"};
{\ar^{q_2} "13";"23"};
{\ar_{0} "21";"22"};
{\ar@{=}_{\mathrm{id}} "22";"23"};
{\ar@{=>}_{\underline{\varpi}_3} (-12,9)*{};(-12,5)*{}};
{\ar@{=>}^{\varpi_3} (12,9)*{};(12,5)*{}};
{\ar@{=>}_{\underline{\varpi}_2} (-12,-5)*{};(-12,-9)*{}};
{\ar@{=>}_{\varpi_2} (12,-5)*{};(12,-9)*{}};
\endxy
\]
Since $\mathrm{Ker}(0:\mathrm{Ker}(\ell)\longrightarrow0 )=\mathrm{Ker}(d)$ by (the dual of) Proposition \ref{final8}, so $\underline{k(f_3)}$ becomes an equivalence. On the other hand, the following is a complex morphism, where $s:=k(f_2)\circ d_2^A$.
\begin{equation}
\xy
(-20,6)*+{\mathrm{Ker}(f_2)}="0";
(-20,-6)*+{\mathrm{Ker}(f_2)}="2";
(0,6)*+{\mathrm{Ker}(f_3)}="4";
(0,-6)*+{A_3}="6";
(20,6)*+{\mathrm{Cok}(f_1)}="8";
(20,-6)*+{Q}="10";
{\ar^{\underline{d}^A_2} "0";"4"};
{\ar@{=}_{\mathrm{id}} "0";"2"};
{\ar_{s} "2";"6"};
{\ar^{k(f_3)} "4";"6"};
{\ar^{q_1} "8";"10"};
{\ar^{d} "4";"8"};
{\ar_{\ell} "6";"10"};
{\ar@/^1.65pc/^{0} "0";"8"};
{\ar@/_1.65pc/_{0} "2";"10"};
{\ar@{=>}_{\underline{\lambda}_2} (-10,2)*{};(-10,-2)*{}};
{\ar@{=>}^{\varpi_3} (10,2)*{};(10,-2)*{}};
{\ar@{=>}_{\alpha} (0,8)*{};(0,12)*{}};
{\ar@{=>} (0,-8)*{};(0,-12)*{}};
{\ar@{..}@/^0.50pc/ (0,-10)*{};(24,-16)*{_{(k(f_2)\circ\varpi_1)\cdot(\varepsilon_{f_2}\circ q)\cdot q^{\flat}_I}}};
\endxy
\label{myan6}
\end{equation}
Thus by taking kernels of $d$ and $\ell$ in diagram $(\ref{myan6})$, we obtain the following factorization by (the dual of) Proposition \ref{final9}.
\[
\xy
(-24,6)*+{\mathrm{Ker}(f_2)}="0";
(-24,-6)*+{\mathrm{Ker}(f_2)}="2";
(0,6)*+{\mathrm{Ker}(d)}="4";
(0,-6)*+{\mathrm{Ker}(\ell)}="6";
(24,6)*+{\mathrm{Ker}(f_3)}="8";
(24,-6)*+{A_3}="10";
{\ar^{\exists(\underline{d}^A_2)^{\dag}} "0";"4"};
{\ar@{=}_{\mathrm{id}} "0";"2"};
{\ar_{\exists\underline{s}} "2";"6"};
{\ar^{\underline{k(f_3)}} "4";"6"};
{\ar^{k(f_3)} "8";"10"};
{\ar^{k(d)} "4";"8"};
{\ar_{k(\ell)} "6";"10"};
{\ar@/^1.65pc/^{\underline{d}^A_2} "0";"8"};
{\ar@/_1.65pc/_{s} "2";"10"};
{\ar@{=>}_{\exists} (-12,2)*{};(-12,-2)*{}};
{\ar@{=>}^{\underline{\varpi}_3} (12,2)*{};(12,-2)*{}};
{\ar@{=>}_{\exists\underline{\alpha}} (0,8)*{};(0,12)*{}};
{\ar@{=>}_{\exists} (0,-8)*{};(0,-12)*{}};
\endxy
\]
Since $(\ref{myandia})$ is 2-exact in $A_3$, so $\underline{s}$ becomes fully cofaithful. Since $\underline{k(f_3)}$ is an equivalence, this means $(\underline{d}_2^A)^{\dag}$ is fully cofaithful, and
\[
\xy
(-24,0)*+{\mathrm{Ker}(f_2)}="4";
(0,0)*+{\mathrm{Ker}(f_3)}="6";
(24,0)*+{\mathrm{Cok}(f_1)}="8";
{\ar_{\underline{d}_2^A} "4";"6"};
{\ar_{d} "6";"8"};
{\ar@/^1.65pc/^{0} "4";"8"};
{\ar@{=>}_{\alpha} (0,2)*{};(0,6)*{}};
\endxy
\]
becomes 2-exact in $\mathrm{Ker}(f_3)$.
\end{proof}

\end{document}